\documentclass{amsart}
\usepackage{latexsym}
\newcommand{\oldversion}[1]{}
\newcommand{\newversion}[1]{#1}

\newcommand{\id}{{\operatorname{id}}}
\newcommand{\Mor}{\operatorname{Mor}}

\newcommand{\Spec}{\operatorname{Sp}}
\newcommand{\Ad}{{\operatorname{Ad}}}

\newcommand{\Graph}{\operatorname{Graph}}

\newcommand{\Hil}{{\mathcal H}}
\newcommand{\skrK}{{\mathcal K}}
\newcommand{\skrL}{{\mathcal L}}

\newcommand{\Dom}{{\mathcal D}}
\newenvironment{pf}{\vs

\noindent{\it Proof.}}{\qed}
\newcommand{\hA}{{\widehat A}}
\newcommand{\ha}{{\widehat a}}
\newcommand{\hb}{{\widehat b}}

\newcommand{\hF}{{\widehat F}}
\newcommand{\hh}{{\widehat h}}
\newcommand{\hJ}{{\widehat J}}

\newcommand{\hQ}{{\widehat Q}}
\newcommand{\hR}{{\widehat R}}
\newcommand{\hS}{{\widehat S}}
\newcommand{\hW}{{\widehat W}}

\newcommand{\hDelta}{{\widehat\Delta}}
\newcommand{\hdelta}{{\widehat\delta}}
\newcommand{\heta}{{\widehat\eta}}
\newcommand{\hkappa}{{\widehat\kappa}}

\newcommand{\hrho}{{\widehat\rho}}
\newcommand{\hsigma}{{\widehat\sigma}}
\newcommand{\htau}{{\widehat\tau}}
\newcommand{\hlambda}{{\widehat\lambda}}
\newcommand{\hgamma}{{\widehat\gamma}}

\newcommand{\Btilde}{{\widetilde B}}
\newcommand{\alphatilde}{{\widetilde \alpha}}

\newcommand{\compl}{{\Bbb C}}
\newcommand{\natu}{{\Bbb N}}
\newcommand{\real}{{\Bbb R}}

\newcommand{\inte}{{\Bbb Z}}

\newcommand{\skrD}{{\mathcal D}}
\newcommand{\skrR}{{\mathcal R}}
\newcommand{\rf}[1]{{\rm (\ref{#1})}}
\newcommand{\reef}[1]{{\rm \ref{#1}}}

\newcommand{\Nideal}{\mathfrak{N}}
\newcommand{\Mideal}{\mathfrak{M}}

\newcommand{\Vs}[1]{\rule{0mm}{#1mm}}
\newcommand{\vs}{\vspace{2mm}}
\newcommand{\comp}{\raisebox{.5mm}{\tiny o}}

\newcommand{\set}[2]{\left\{#1:#2\right\}}
\newcommand{\seq}[2]{\ntuple{#1}_{#2}}
\newcommand{\ntuple}[1]{\left\{#1\right\}}
\newcommand{\CLS}[1]{{\mbox{CLS}}\left\{#1\right\}}

\newcommand{\norm}[1]{\|#1\|}

\numberwithin{equation}{section}
\newtheorem{Thm}{Theorem}[section]
\newtheorem{Def}[Thm]{Definition}
\newtheorem{Lem}[Thm]{Lemma}
\newtheorem{Prop}[Thm]{Proposition}
\newtheorem{Cor}[Thm]{Corollary}
\newtheorem{rem}[Thm]{Remark}
\newenvironment{Rem}{\begin{rem}\em}{\end{rem}}

\newcommand{\IS}[2]{\left(#1\,\vline\,#2\right)}
\newcommand{\ITS}[3]{\left(#1\,\vline\,#2\,\vline\,#3\right)}

\newcommand{\xxx}[1]{}

\begin{document}
\title{A C$^*$-algebraic framework for quantum groups}
\author{T.~MASUDA}
\address[T.~MASUDA]{Institute of Mathematics, University of Tsukuba, 
1-1, Tennodai 1 chome,
Tsukuba-shi, Ibaraki, 305-8571 Japan.}
\email[T.~MASUDA]{tetsuya@math.tsukuba.ac.jp}
\author{Y.~NAKAGAMI}
\address[Y.~NAKAGAMI]{Department of Mathematical and Physical 
Science, Faculty of Science, Japan
Women's University, 8-1, Mejirodai 2 chome, Bunkyo-ku, Tokyo, 112-8681 Japan.}
\email[Y.~NAKAGAMI]{nakagami@fc.jwu.ac.jp}
\author{S. L.~WORONOWICZ}
\address[S. L.~WORONOWICZ]{Department of Mathematical Methods in 
Physics, Faculty of Physics,
University of Warsaw, Ho\.{z}a 74, 00-682 Warszawa, Poland.}
\email[S. L.~WORONOWICZ]{slworono@fuw.edu.pl}
\date{}

\begin{abstract} We develop a general framework to deal with the 
unitary representations of quantum
groups using the language of
$C^*$-algebras. Using this framework, we prove that the duality holds 
in a general context. This
extends the framework of the duality theorem using the language of 
von Neumann algebras previously
developed by Masuda and Nakagami.
\end{abstract}

\maketitle

\section{Introduction}
When we look at quantum groups, we still have the shortage of
examples which are studied in detail from all \oldversion{the aspects of
mathematics, in particular the non-compact ones.} \newversion{points 
of view. This
applies mostly to noncompact groups.} As long as we deal with the 
examples of the
quantum groups of the compact type, the existence of the Haar 
measure, the (relative)
invariance of the Haar measure with respect to the scaling automorphism group
etc. are all clear even from \oldversion{the} \newversion{a} very 
general point of
view. Due to the Peter-Weyl theorem for quantum groups of classical type, we
can characterize the object by making use of the set of all matrix
elements $W=\{w_{ij}\}_{ij}$ of all possible finite dimensional
(unitary) representations  without using functional analysis.
\oldversion{Though} \newversion{However}, even in the compact case, the
C$^*$-algebraic framework is essential to guarantee the existence of the Haar
measure, from which \oldversion{we have developed} almost all the theory of
representations \newversion{has been developed}~\cite{w3}. \vs

However, once we become interested in quantum groups of
non-compact type, we are faced \newversion{for instance} with a 
serious necessity for
the convergence of the infinite sums of matrix elements of infinite
dimensional unitary representations\oldversion{, for instance}. Furthermore,
the Haar measures correspond to unbounded linear functionals and
the elements which generate the coordinate ring of the quantum
group are unbounded linear operators, which cause analytical
difficulties. Therefore we are obliged to work with some
appropriate notions and the basic tools of functional analysis, in
particular the language of operator algebras. Despite these
difficulties, there are some quite explicitly studied quantum
groups like quantum Lorentz groups, $E_q(2)$, quantum `$az+b$'
group, quantum `$ax+b$'\linebreak group, $\widetilde{SU}_q(1,1)$ etc. of
non-compact type~\cite{pw}, \cite{wz}, \cite{w2}, \cite{w6},
\cite{w7}, \cite{kor}, \cite{kk} with infinite dimensional
irreducible unitary representations. Furthermore we are still
discovering new explicitly computable quantum groups of
non-compact type. All of them should be investigated from the
point of view of algebra, analysis and geometry. \vs \xxx{}

Now, the history of the operator algebraic study of the classical
locally compact groups (with the duality theorem) goes back to
Stinespring who transferred the regular representation to what we
nowadays call the {\it multiplicative unitary operator} for the
purpose of generalizing the Tannaka-Krein's duality to \oldversion{the}
unimodular locally compact groups. Then it was Kac and Takesaki
who independently used this unitary operator to extend the duality
theorem for the category which is much wider than that of locally
compact groups. Immediately after that, it was Tatsuuma who
recognized the importance of this unitary operator and named it
the "Kac-Takesaki operator". After that, Enock-Schwartz and
Kac-Vainermann completed the theory of generalized locally compact
groups, which they called ``Kac algebras'' and introduced the
``fundamental unitary'' playing the role of Kac-Takesaki operator.
The details concerning the Kac algebras including the history of
the development of the subject (before the appearance of \oldversion{the}
quantum groups) \oldversion{are} \newversion{can be} found in their 
book~\cite{eSch}.
It is also remarked that in the development of those theoretical frameworks
to deal with the general theory of locally compact groups
including the duality theorem, all the achievements were made on
the basis of von Neumann algebras by the heavy use of the
Tomita-Takesaki theory. The only approach based on the language of
C$^*$-algebras was proposed by Vallin~\cite{val1}, in which the
theory was not developed in full generality due to the lack of the
reasonable integration theory for C$^*$-algebras. It is also
emphasized that the condition which we now call the {\it strong
$($right$)$ invariance} had been used from the early stage for the
development of the subject. \vs

In the middle \oldversion{of} eighties, after the famous discoveries of
Leningrad school, Drinfeld and Jimbo introduced quantum groups as
deformations of enveloping algebras of semisimple Lie algebras. At
the beginning these objects were of purely algebraic nature. In
particular only finite-dimensional representations and no
$*$-involutions were considered. \vs

A little later the third author introduced a noncommutative
analogue of the algebras of functions on compact groups. Then the
$*$-involution was an essential part of the considered structure.
It allowed to distinguish compact and noncompact quantum groups
and to consider unitary representations of quantum groups. For
noncompact groups unitary representations can be infinite
dimensional \oldversion{that require} \newversion{which requires}
the use \oldversion{of concepts}
of analytical methods. One of the aims of this publication is to draw attention
to the importance of the functional analytical methods from this
point of view. \vs

It was the paper of Baaj and Skandalis~\cite{bs} who first gave a
reasonable framework to deal with such problems on the basis of
the multiplicative unitary operator. Their approach avoided the
discussion of \oldversion{the} Haar measures \oldversion{and} 
\newversion{while} the
duality and the fact that the dual of the dual is isomorphic to the 
original object
were a trivial consequence of the definitions. Their theory in large
part worked with C$^*$-algebras. \vs

In~\cite{mn}, the first and the second authors gave an alternative
formulation of the theory including the discussion of the Haar
measures with the duality being a deep theorem. That paper was
based on the straightforward use of the Tomita-Takesaki theory to
deal with the Haar measure \oldversion{in a general context} following the
above mentioned theory of Kac algebras. We were obliged to
sacrifice the topological nature of the underlying quantum groups
which we believe to be important nowadays. This is why we need a
reasonable framework to deal with the quantum groups as a
generalization of locally compact groups. For this purpose, the
right choice of the language and the axioms are important. \vs

Here we would like to list a number of requirements that any good
axiomatization of the theory of locally compact  quantum groups
should satisfy. In our opinion the quantum theory should be
developed in the way parallel to the classical one. In particular
this rule determines the choice of basic concepts. Classically the
groups are locally compact spaces endowed with a group rule. By
the famous Gelfand-Naimark theorem locally compact spaces are in
one to one correspondence with commutative C$^*$-algebras (unital
if the considered spaces are compact). Following the general
philosophy, passage to \oldversion{the} quantum spaces consists in using
noncommutative algebras. Consequently the underlying spaces of
quantum groups should be described by noncommutative (in general)
C$^*$-algebras. The use of von Neumann algebras seems not to be
appropriate. They correspond (in commutative case) to measure
spaces. Of course, one is able to develop the theory of groups with
underlying spaces being measure spaces,
but this theory is less natural and more complicated than the
theory of locally compact groups. \vs

The next reason why one should work with the language of
C$^*$-algebras is the more natural description of homomorphisms of
(quantum) groups. Assume that we have a locally compact subgroup
$H$ (of lower dimension) of a locally compact group $G$.
\oldversion{Working within the C$^*$-algebra language we describe the
inclusion $H\subset G$ be the restriction morphism $C_0(G)\ni a\rightarrow
a|_H\in C_0(H)$}
\newversion{Within the C$^*$-algebra language, the inclusion $H\subset
G$ is described by the restriction map $C_0(G)\ni a\rightarrow a|_H\in
C_{\rm bounded}(H)$, which is a morphism  (in the sense of
\cite{w11}; cf. Definition \ref{defmor}) from $C_0(G)$ into $C_0(H)$}.
In the von Neumann algebra setting the corresponding restriction map
$L^\infty(G) \rightarrow L^\infty(H)$ does not make sense.   \vs

Now we pass to the second ingredient of the concept of a group
which is the binary operation called the group rule. Let $G$ be a
group. Classically the group rule is a continuous mapping from
$G\times G$ onto $G$. Passing to the quantum theory we replace
cartesian product of spaces by tensor product of algebras
(alternative choice would be the free product of algebras
but then all the classical groups would be lost). Let $\alpha$ be
a C$^*$-norm on the algebraic tensor product $A\otimes_{\rm alg}B$
of C$^*$-algebras $A$ and $B$. By definition, tensor product
$A\otimes_\alpha B$ is defined as the completion of $A\otimes_{\rm
alg}B$ with respect to the norm $\alpha$. \vs

Let $A$ be a C$^*$-algebra related to a quantum group. The
correspondence bet\-ween spaces and algebras is a contravariant
functor. It changes the direction of morphisms. Therefore the
group rule corresponds to a morphism $\delta$ from $A$ into
$A\otimes_\alpha A$. This is the so called {\it comultiplication}
(we use small letter $\delta$, because the 
\oldversion{notation}\newversion{symbol}
$\Delta$ is reserved to denote the modular operator in Tomita-Takesaki
theory). In this context one has to pay some attention to the
notion of morphism of C$^*$-algebras. The usual $*$-algebra
homomorphisms are not appropriate to use in the context of quantum
spaces. They correspond (in the commutative case) to the proper
continuous mappings of underlying locally compact spaces. On the
other hand the group rule is not a proper mapping from $G\times G$
into $G$ if $G$ is not compact. Therefore one has to modify the
notion of morphism creating a new category of C$^*$-algebras.
Restricted to commutative algebras the category should be
(anti)-equivalent (via Gelfand Naimark theorem) to the category of
locally compact spaces and (all) continuous mappings. Such a
category of C$^*$-algebras was proposed in~\cite{w0} (see
also~\cite{val1}, \cite{w1}, \cite{w11}). For reader's convenience
we recalled the definition of morphisms of C$^*$-algebras in
Appendix \reef{A} (Definition~\ref{defmor}).\vs\vs

In general we have more than one C$^*$-norm on $A\otimes_{\rm
alg}A$. Among them there is the minimal one that corresponds to
the spatial tensor product of C$^*$-algebras denoted simply by
$A\otimes A$. Due to the minimality, for any C$^*$-norm $\alpha$
we have canonical epimorphism $\Psi_\alpha:A\otimes_\alpha
A\rightarrow A\otimes A$. Replacing $\delta$ by
$\Psi_\alpha\comp\delta$ we obtain a morphism from $A$ into
$A\otimes A$. Therefore the choice of spatial tensor product as
the operation corresponding to the cartesian product of spaces
produces the theory with larger (than with other choices) class of
examples. It is also worth mentioning that in most known examples
the algebra $A$ is nuclear and all tensor products
$A\otimes_\alpha A$ coincide.\vs

Summing up: The basic concepts of the theory of locally compact
quantum groups should be a C$^*$-algebra (say $A$) together with a
morphism $\delta\in\Mor(A,A\otimes A)$. In the last expression
'$\otimes$' and `$\Mor$' denote the spatial tensor product of
C$^*$-algebras and the set of morphisms in the sense of
Definition~\ref{defmor}. \vs

If $A$ is commutative then by Gelfand Naimark theory
$A=C_\infty(G)$, where $G$ is a locally compact topological space
and $\delta(a)(g,g')=a(gg')$, where
\begin{equation}
\label{dzialaniegr} G\times G\ni(g,g')\longmapsto gg'\in G
\end{equation}
is a continuous mapping.\vs

Now we shall discuss the most desired system of axioms. Quantum
groups with commutative $A$ should coincide with the usual groups.
In this case our axioms should reduce to the axioms or {\it simple
theorems} of the theory of locally compact topological groups.
This requirement is satisfied by axioms of the coassociativity of
$\delta$, the properness and the cancellation law (see Definitions
\ref{defProp} and \ref{defCan} in the next section). The
coassociativity clearly  corresponds to the associativity of the
group rule \rf{dzialaniegr}: $(g_1g_2)g_3=g_1(g_2g_3)$ for any
$g_1,g_2,g_3\in G$. The properness means that the mappings
$G\times G\ni(g,g')\mapsto(g,gg')\in G\times G$ and $G\times
G\ni(g,g')\mapsto(gg',g')\in G\times G$ are proper i.e: $(g,gg')$
and $(gg',g')$ tend to infinity when $(g,g')$ tends to infinity.
The cancellation property states that for any $g,g_1,g_2\in G$ we
have:
\[
\begin{array}{c}
(gg_1=gg_2)\Rightarrow(g_1=g_2),\\
(g_1g=g_2g)\Rightarrow(g_1=g_2).\Vs{5}
\end{array}
\]
A good candidate for an axiom is the statement of Theorem
\ref{newAx}. Combined with the properness and cancellation law it
reduces (for commutative $A$) to the postulate saying that the
relations $g=hh'$ and $h=gg'$ define a homeomorphism
\[
G\times G\ni(g,g')\longmapsto(h,h')\in G\times G
\]
acting on $G\times G$. Unfortunately at present we are not able to
base our theory on axioms of the form similar to the statement
mentioned above.\vs

By simple theorems of the theory of locally compact topological
groups (the phrase used in the previous paragraph) we mean
statements that may be derived from axioms by short and easy
reasoning. In principle they should be found (explicitly or
implicitly) on first pages of handbooks devoted to the subject.
The Haar theorem stating the existence of the right (or left)
invariant measure does not belong to the category of simple
theorems. In a good theory of locally compact quantum groups, the
existence of the Haar weight should be one of advanced theorems.
\vs

This is why the existing axiomatizations of the theory of locally
compact groups are not sa\-tisfactory. In all the approaches
(including that of Kustermans and Vaes \cite{Vaes} and the one
used in the present paper) the existence of the Haar measure is
included as one of the axiom of the theory. We  strongly believe
that this is a temporary situation and that in the future we shall
have a  theory that will satisfy the requirements  listed above.
To find a satisfactory axiomatization it is important to
understand well the logical relations between the basic concepts
and statements of the theory. This is the main purpose of the
present paper. \vs

\oldversion{Under these observations, our joint research began by aiming
to develop a reasonable C$^*$-algebraic framework for the quantum
groups. We are going to improve our previous approach~\cite{mn}
based on the language of von Neumann algebras. Our first results
were announced by the third author at Fields Institute and by the
first author at Banach Center in 1995. Then the first author also
announced our preliminary results in an expository
article~\cite{m}. \vs

Taking the above discussions into account, we aim to present a
reasonable framework, as general as possible, to formulate (within
the formalism of C$^*$-algebras) the objects which we recognize as
the locally compact quantum groups.}

\newversion{Taking the above discussions into account, we aim to present a
reasonable framework, as general as possible, to formulate (within the
formalism of C$^*$-algebras) the objects which we recognize as \oldversion{the}
locally compact quantum groups. We are going to improve our previous
approach~\cite{mn} based on the language of von Neumann algebras. Our
first results were announced by the third author at Fields Institute and
by the first author at Banach Center in 1995. Then the first author also
announced our preliminary results in an expository article~\cite{m}.} On
the way\oldversion{ to develop our present version}, the third author
modified the theory of multiplicative unitaries~\cite{w4} replacing the
(semi)-regularity condition of Baaj and Skandalis by a new one better
suited to our aims. While the present version of our paper was being
developed, Kustermans and Vaes~\cite{Vaes} found an alternative approach
to the locally compact quantum groups starting from a completely
different set of axioms \newversion{and} proving the duality theorem. 
However, in
both cases, the existence of the Haar measure(s) is an
\oldversion{unavoidable} \newversion{inevitable} assumption. This  means
that the both are still on the way to be completed. The existence of the
Haar measure should be derived from a more primitive set of axioms. \vs

Now we have a final comment on the strong right invariance which
is one of our axioms. It is the only axiom that relates the
antipode with remaining structure (comultiplication, Haar measure).
On the other hand, Kustermans and Vaes~\cite{Vaes} proposed to
assume the existence of both the right and the left Haar measures.
Then the existence of the antipode is \oldversion{an uneasy} \newversion{a
sophisticated} theorem. It seems that, in both cases, those conditions are
essentially related to the {\it manageability} which plays an 
important role in both
approaches. \vs

In~\cite{mn} we discussed the properties of the special one
parameter automorphism groups which we now call the scaling
groups. This is the symmetry first recognized in~\cite{mncmp}
which describes an obstruction of our object to being a locally
compact group or a Kac algebra. We also remark that there are many
improvements of the proofs as well as the basic axioms comparing
with what is given in~\cite{mn}:
\begin{enumerate}
\item 
The main definition \ref{def1.2} is simplified: the commutativity
of the coproduct with the unitary antipode as well as with the
scaling group; and the commutativity of the left and the right
Haar weights; are not assumed.
\item 
In case of weights (contrary to the states) the commutativity of
the two weights $\varphi$ and $\psi$ does not imply that the
Connes' Radon-Nikodym cocycle $(D\psi:D\varphi)_t$ is invariant
under the group of modular automorphism
$\seq{\sigma_t^\varphi}{t\in\real}$.
\item 
There is a gap in the proof of Lemma 2.14 in~\cite{mn}.  The
assertion is true by Theorem \ref{th6.10} of the present paper.
\item
The invariance of the Haar weight with respect to the scaling
automorphisms is replaced by relative invariance.
\end{enumerate}
\vs

The relative invariance (instead of the invariance) of the Haar
weight under the scaling group appeared first in the paper of
Kustermans and Vaes~\cite{Vaes} as a result of their axioms.
Recently, Van Daele~\cite{Fons} found that the quantum groups
\oldversion{$ax+b$ and $az+b$} \newversion{`$ax+b$' and `$az+b$'}
introduced in~\cite{w6,w7} have Haar measures which are relatively
invariant but not invariant with respect to the scaling automorphisms. To
include such examples in our framework we assume the relative invariance
of the Haar weight. This modification was very easy, it resulted in a
number of minor changes. \vs

\oldversion{Now we summarize the organization of this paper. In Section 1,
we give our basic definitions and the axioms of weighted Hopf
C$^*$-algebra which we discuss throughout this paper. Section 2 is
devoted to an investigation of matrix elements of the Kac-Takesaki
operator. In Section 3 we prove that the operator is a
multiplicative unitary and that it is manageable in the sense
of~\cite{w4}. In Section 4, we discuss the Radon-Nikodym
derivative of the left Haar weight with respect to the right one.
Section 5 is devoted to the construction of the dual modular
operator. In Section 6 we introduce the dual Haar weight. To this
end  the dual GNS-mapping is defined. In Section 7, we introduce
the dual weighted Hopf C$^*$-algebra and prove the duality.
Section 8 is devoted to the quantum codouble construction (which
is interpreted to be a dual version to the Drinfeld's quantum
double) to deal with examples like the quantum Lorentz group
$SL_q(2,{\Bbb C})$.}

\newversion{Now we summarize the organization of this paper. In Section
\ref{ex1}, we give our basic definitions and the axioms of weighted 
Hopf C$^*$-algebra which we
discuss throughout this paper. Then using the GNS construction we 
introduce the Hilbert space and
reveal the content  of the strong right invariance. Section \ref{ex2} 
is devoted to an
investigation of matrix elements of the Kac-Takesaki operator. We 
prove the theorems describing the
{\it coalgebraic} properties of scaling group and unitary antipode. 
In Section \ref{sek3} we show
that the Kac-Takesaki operator is a multiplicative unitary and that 
it is manageable in the sense
of~\cite{w4}. Then the general theory of multiplicative unitaries 
enable us to introduce the basic
objects related to the dual weighted Hopf C$^*$-algebra. In Section 
\ref{sec4}, we discuss the
Radon-Nikodym derivative of the left Haar weight with respect to the 
right one. This derivative is
determined by two operators $\rho$ and $\gamma$. Later in Section 
\ref{sek6} we show that $\rho$ is
affiliated with the considered C$^*$-algebra and $\gamma$ is a 
multiple of unity
related to the constant describing the relative invariance of the 
Haar weight. In
the last part of Section \ref{sec4} we prove the uniqueness of the 
antipode and the
Haar weight. Section \ref{sek5} is devoted to the properties of an operator
$\hDelta$ which plays an important role in our analysis. In Section 
\ref{sek6} we
introduce the dual GNS-mapping $\heta$. Operator $\hDelta$ introduced in the
previous section turns out to be the modular operator related to 
$\heta$. At the end
of this section we discuss the behaviour of the (right invariant) 
Haar weight with
respect to left shifts. In Section \ref{sek7}, we introduce the dual 
Haar weight.
This way the construction of dual weighted Hopf C$^*$-algebra is 
achieved and the
duality is shown. In our work, we deal all the time with the ``regular''
representations and the ``regular'' dual in order to have faithful 
dual Haar weight.
This means that our multiplicative unitary operator is the 
``regular'' bicharacter.
However, we are also able to discuss the universal dual with the universal
bicharacter. This topic will be discussed in a separate publication. Section
\ref{sek8} is devoted to the quantum codouble construction which is a 
dual version
to the Drinfeld's quantum double. It helps to deal with examples like 
quantum Lorentz
groups.} \vs

\newversion{Then we have six Appendices ranging from \ref{A} to
\ref{AnaGene}. This is where we collect the technicalities used in our
paper. A few fundamental facts on the multiplier algebra and the
definition of the set
$\Mor(A,B)$ of morphisms are given in Appendix \ref{A}. Appendix
\ref{GNSmap} contains the concepts of GNS mapping, its commutant (and the
corresponding bicommutant theorem) and (exact) vector presentation of a
GNS mapping. These concepts provide a convenient language for the
Tomita-Takeski theory outlined in Appendix \ref{TomTak}. Appendix
\ref{B} is devoted to weights and their properties (such as strict
faithfulness) on a separable C$^*$-algebras. Tensor products of the
GNS-mappings are discussed in Appendix \ref{Tensor}. In the last Appendix
\ref{AnaGene} we present the basic properties of analytic generators of
one parameter groups acting on Banach spaces.}
\vs

Some remarks on the terminology used in this paper: First of all,
we decided to use the name ``weighted Hopf C$^*$-algebra'' for the
basic object of this paper. We have to acknowledge that the
terminology ``Hopf C$^*$-algebra'' is already used in the paper of
Vallin~\cite{val1} and also in the paper of Baaj and
Skandalis~\cite{bs} in a different context. The word ``weighted''
means that we deal with the object with a Haar weight. Then, the
one parameter family of automorphism which comes from the ``radial
part'' of the antipode is called the scaling group. (In \oldversion{our
previous publication}~\cite{mn}, we called it the deformation
automorphism group.) We hope that our choice of the terminology
will not cause any confusion to the readers. \vs

\newversion{
Another change which we have to acknowledge is the choice of the particular
form of Kac-Takesaki operator $W$. In general, the Kac-Takesaki operator is
introduced by choosing a certain linear mapping acting on the tensor 
product of the
function algebra by itself and  pushing it down to the Hilbert space 
level. There
are several possibilities. In~\cite{mn} we used the mapping
\begin{equation}
\label{defW}
  a\otimes b \longmapsto \delta (b)(a\otimes 1).
\end{equation}
This choice forced us to work with the left invariant Haar weight and with
the strong left invariance axiom. In the present paper, we replace \rf{defW} by
the mapping
\[
a\otimes b \longmapsto \delta (a)(1\otimes b).
\]
Consequently we work with the right invariant Haar weight and
the strong right invariance. This choice is consistent with
the notational convention of Baaj and Skandalis~\cite{bs} and the framework
developed in~\cite{w4}.}\vs

\par
{\bf Acknowledgement}. The first two authors are grateful to the
Department of Mathematical Methods in Physics of the University of
Warsaw, and the last two authors are grateful to the Department
of Applied Mathematics of Fukuoka University and the Department of
Mathematics of Norwegian University of Science and Technology. The
first author was partially supported by the Grant-in-Aid for
Scientific Research (No. 11440052) and the second author was
partially supported by the Grant-in-Aid for Scientific Research
(No. 12640216), the Ministry of Education, Science, Sports and
Culture, Japan. The third author was supported by Polish KBN,
grants No 2 P301 020 07, No 2 P0A3 030 14 and No 2 P03A 040 22,
by Foundation for Polish Science and by Japan Society for Promoting Science.

\section{Basic definitions and axioms}
\label{ex1}

The \newversion{main} aim of this section is \oldversion{to give our 
definition}
\newversion{to introduce the notion} of a weighted Hopf
$C^*$-algebra\oldversion{which we consider to be the basis of the whole
discussions of this paper. As we have already mentioned in the 
introduction, we use
the terminology "Hopf $C^*$-algebra"  even though the same 
terminology had already
been used in a different context}.\vs

\oldversion{We also have some remarks on the notational convention 
that the unit 1
of the $C^*$-algebra $A$ always means the unit element of the 
multiplier algebra
$M(A)$. Then $\id$ is the notation for the identity mapping.  For the inner
product the physicist convention is used, i.e. $\IS{y}{x}$ is linear in
$x$ and conjugate linear in $y$. Throughout the discussions of this 
paper, the $C^*$-algebras and the
Hilbert spaces are always {\em assumed to be separable} without 
specially mentioned. Then all the
tensor products of the $C^*$-algebras which appear in this paper are 
{\em minimal tensor products}
without specially mentioned. We shall use compositions of the form 
$\varphi\comp\kappa$, where
$\varphi\in A^*$ and $\kappa$ is a closed operator acting on $A$. We write
$\varphi\comp\kappa\in A^*$ if there exists a $\psi\in A^*$ such that 
$\psi(a)=\varphi(\kappa(a))$
for all $a\in\Dom(\kappa)$. Clearly $\psi$ is unique and will be  denoted by
$\varphi\comp\kappa$.\vs}

\newversion{Some remarks on used notation. The unit $1=1_A$
of a $C^*$-algebra $A$ always means the unit element of the multiplier algebra
$M(A)$. Then $\id$ denotes the identity mapping.  For the inner
product the physicist convention is used, i.e. $\IS{y}{x}$ is linear in
$x$ and conjugate linear in $y$. Throughout this paper, all 
$C^*$-algebras and all
Hilbert spaces are {\em assumed to be separable}.
Moreover all tensor products of $C^*$-algebras are {\em minimal 
tensor products}
unless specially mentioned. We shall use compositions of the
form $\varphi\comp\kappa$, where
$\varphi\in A^*$ and $\kappa$ is a closed operator acting on $A$. We write
$\varphi\comp\kappa\in A^*$ if there exists a $\psi\in A^*$ such that 
$\psi(a)=\varphi(\kappa(a))$
for all $a\in\Dom(\kappa)$. Clearly $\psi$ is unique and is denoted by
$\varphi\comp\kappa$.\vs}

\oldversion{Before going into the description of our basic axiom, we 
need to present
some technical lemmas concerning the multipliers and the weights on the
$C^*$-algebras. We refer these technical matters to Appendices \reef{A} and
\reef{B}.}

\newversion{To understand well the forthcoming definitions and results the
reader is advised to look at Appendices \reef{A} and \reef{B}, where the basic
notions such as multipliers, morphisms and weights are explained.\vs}

\par Now, we are ready to give our basic definition of the weighted Hopf
$C^*$-algebra.\vs

\begin{Def}
\label{def1.1}
A pair $(A,\delta )$ where $A$ is a $C^*$-algebra and
$\delta\in\Mor (A,A\otimes A)$ satisfying the coassociativity 
condition $(\delta\otimes
\id_A)\comp\delta =(\id_A\otimes\delta)\comp\delta$ is called a {\em
$C^*$-bialgebra}.
\end{Def}

The mapping $\delta$ will be called a comultiplication or a 
coproduct.  In what follows we denote
by $1=1_A$ the unit of the multiplier algebra $M(A)$. The relation 
$\delta\in\Mor (A,A\otimes A)$
does not mean that $\delta(a)\in A\otimes A$ for $a\in A$. We know only that
$\delta(a)\in M(A\otimes A)$. It means that $\delta(a)(b\otimes c)\in 
A\otimes A$ and
$(b\otimes c)\delta(a)\in A\otimes A$ for any $a,b,c\in A$. Therefore 
the following definitions are
relevant:

\begin{Def}
\label{defProp}
A C$^*$-bialgebra $(A,\delta )$ is said to be proper if 
$\delta(a)(1\otimes b)$ and
$(b\otimes 1)\delta(a)$ belong to $A\otimes A$ for any $a,b,c\in A$.
\end{Def}

\begin{Def}
\label{defCan}
Let  $(A,\delta )$ be a proper C$^*$-bialgebra. We say that 
$(A,\delta )$ has the cancellation pro\-perty
if the linear spans of the sets $\set{\delta(a)\Vs{3}(1\otimes 
b)}{a,b\in A}$ and
$\set{(b\otimes 1)\Vs{3}\delta(a)}{a,b\in A}$ are dense in $A\otimes A$.
\end{Def}

Let $(A,\delta )$ be a proper C$^*$-bialgebra. Then for any 
continuous linear functionals
$\varphi,\psi\in A^*$ and any $a\in A$, the convolution products
\begin{equation}
\label{conpro}
\begin{array}{c}
\varphi*a:=(\id\otimes\varphi)\delta(a),\\
a*\psi:=(\psi\otimes\id)\delta(a).\Vs{5}
\end{array}
\end{equation}
belong to $A$. Indeed, by Proposition \reef{AP2}, we may assume that 
$\varphi=b\varphi'$ and $\psi=\psi'c$, where $\varphi',\psi'\in
A^*$ and $b,c\in A$. Then
\[
\begin{array}{c}
\varphi*a=(\id\otimes\varphi')\left\{\delta(a)(1\otimes b)\Vs{3}\right\},\\
a*\psi=(\psi'\otimes\id)\left\{(c\otimes 1)\delta(a)\Vs{3}\right\}\Vs{5}
\end{array}
\]
and our statement follows immediately from the properness of the 
considered bialgebra. We shall also
use the convolution product of functionals:
\[
\varphi*\psi:=(\varphi\otimes\psi)\comp\delta\in A^*.
\]
One can easily verify that
\[
(\varphi*\psi)(a)=\varphi(\psi*a)=\psi(a*\varphi).
\]
Due to the coassociativity of the comultiplication, the convolution 
product obey the following
associativity laws:
\[
\begin{array}{c}
(\psi*\varphi)*a=\psi*(\varphi*a),\\
(a*\psi)*\varphi=a*(\psi*\varphi),\Vs{4}\\
(\varphi*a)*\psi=\varphi*(a*\psi),\Vs{4}\\
(\varphi*\chi)*\psi=\varphi*(\chi*\psi).\Vs{4}
\end{array}
\]
In these formulae $a\in A$ and $\varphi,\chi,\psi\in A^*$.\vs

Throughout the paper, the concept of a right invariant weight on a 
C$^*$-bialgebra plays an essential
role. In brief, $h$ is right invariant if $(h\otimes\id)\delta(a)=h(a)\,1$.

\begin{Def}
\label{hinv}
Let $(A,\delta)$ be a C$^*$-bialgebra and $h$ be a locally finite, 
lower semicontinuous weight on
$A$. We say that $h$ is right invariant if for any $a\in A_+$ such 
that $h(a)<\infty$ and any
$\varphi\in A^*_+$ we have
\[
h(\varphi*a)=\varphi(1)h(a).
\]
\end{Def}

We shall use standard notation:  For any weight $h$ we set:
\[
\begin{array}{r@{\;}c@{\;}l}
\Nideal_h&=&\set{a\in A}{h(a^*a)<\infty},\\
\Mideal_h&=&\Vs{5}\set{a\in A_+}{h(a)<\infty}^{\rm linear\ span}\\
&=&\Vs{5}\set{a^*b}{a,b\in\Nideal_h}^{\rm linear\ span}
\end{array}
\]
It is known that $h$ extends uniquely to a linear functional on $\Mideal_h$.\vs

Assume that $h$ is right invariant. For any $\varphi\in A^*_+$ the 
mapping $a\mapsto \varphi*a=
(\id_A\otimes\varphi)(\delta (a))$ is completely positive. Hence, by 
the Kadison inequality
\[
\{\varphi*a\}^*\{\varphi*a\}
         \leq \varphi(1)\{\varphi*(a^*a)\}\;.
\]
Now, using the right invariance we see that
\begin{equation}
\label{Kadison}
h((\varphi*a)^*(\varphi*a))\leq \varphi(1)^2 h(a^*a)<\infty \;.
\end{equation}
for any $a\in\Nideal_h$. Remembering that any $\varphi\in A^*$ can be 
written by a linear combination of
four positive functionals we conclude that $\varphi*a\in\Nideal_h$ 
for any $a\in\Nideal_h$
and $\varphi\in A^*$. Therefore for any $a,b\in\Nideal_h$ and any 
$\varphi,\psi\in A^*$ the products
$(\varphi^**a)^*b=(\varphi*a^*)b$ and $a^*(\psi*b)$ belong to 
$\Mideal_h$ and expressions
$h((\varphi*a^*)b)$ and
$h(a^*(\psi*b))$ make sense. We shall use them in our main definition.\vs

The object which we study in this paper is defined as follows.

\begin{Def}
\label{def1.2}
Let $(A,\delta )$ with $A$ separable be a proper C$^*$-bialgebra with 
the cancellation property. We
say that $(A,\delta)$ is a {\em weighted Hopf $C^*$-algebra} if there 
exist a closed densely defined linear map $\kappa$
acting on $A$ $($called {\em antipode}$)$ and a locally finite strictly
faithful lower semicontinuous right invariant weight $h$ on $A$ 
$($called {\em Haar weight}$)$ such
that
  \begin{enumerate}
  \item
The operator $\kappa$ admits the following polar decomposition$:$
\begin{equation}
\label{poldec}
\kappa=R\comp\tau_{i/2}
\end{equation}
where $\tau_{i/2}$ is the analytic generator of a one parameter group
$\tau=\seq{\tau_t}{t\in\real}$ of automorphisms of the $C^*$-algebra 
$A$ $($called {\em scaling group}$)$
and $R$ is an involutive antiautomorphism of $A$ $($called {\em 
unitary antipode}$)$ commuting with
automorphisms $\tau_t$ for all $t\in\real$. \vs

  \item 
  The weight $h$ is relatively invariant under $\seq{\tau_t}{t\in\real}$
  $:$ there exists a positive scalar $\lambda$ such that
  $h\comp\tau_t=\lambda^{t}h$ for all $t\in\real$.\vs

\item
$(${\em Strong right invariance}$)$
For any $\varphi\in A^*$ with $\varphi\comp\kappa\in
A^*$ we have
  \[
  h((\varphi*a^*)b) =h(a^*((\varphi\comp\kappa)*b))
  \]
  for all $a,\ b\in \Nideal_h$.
  \end{enumerate}
\end{Def}

It should be noted that the scaling constant $\lambda$ is 1 for the 
most examples of quantum groups
constructed so far. The first examples with $\lambda\neq 1$ appeared 
in \cite{w6,w7} (cf.
\cite{Fons}). This possibility was foreseen for the first time in 
\cite{Vaes}.\vs

It is understood that the group $\tau$
is pointwise norm continuous: for any
$a\in A$, $\norm{\tau_t(a)-a}\rightarrow 0$ when $t\rightarrow 0$.
We refer to Appendix \reef{AnaGene} for the definition of the 
analytic generator. Clearly, $\kappa^2=\tau_i$. Performing
the holomorphic continuation in the formula 
$\tau_t(ab)=\tau_t(a)\tau_t(b)$ we see that $\Dom(\tau_{i/2})$ is a
subalgebra in $A$ and $\tau_{i/2}(ab)=\tau_{i/2}(a)\tau_{i/2}(b)$ for any
$a,b\in\Dom(\tau_{i/2})$. Similarly performing the holomorphic continuation
in the formula $\tau_t(a^*)=\tau_t(a)^*$ we see that 
$\Dom(\tau_{i/2})^*=\Dom(\tau_{-i/2})$ and
$\tau_{-i/2}(a^*)=\tau_{i/2}(a)^*$ for all $a\in\Dom(\tau_{i/2})$. Therefore
$(\tau_{i/2}\Dom(\tau_{i/2}))^*=\Dom(\tau_{i/2})$ and 
$\tau_{i/2}(\tau_{i/2}(a)^*)^*=a$ for any
$a\in\Dom(\tau_{i/2})$. Remembering that $R$ is an involutive 
antiautomorphism of $A$ commuting with $\tau_t$ we
conclude that
\[
\begin{array}{c}
\kappa(ab)=\kappa(b)\kappa(a),\\
\kappa(\kappa(a)^*)^*=a\Vs{5}
\end{array}
\]
for any $a,b\in\Dom(\kappa)$. In other words the mapping 
$a\mapsto\kappa(a)^*$ is a conjugate linear multiplicative
involution acting on $\Dom(\kappa)$.

\begin{Rem}\label{rmk.Uniqueness} The polar decomposition \rf{poldec} 
is unique. The scaling group and the  unitary
antipode are determined by $\kappa$.
Indeed $\tau_i=\kappa^2$ and using Corollary $\reef{uniquegr}$ we 
obtain the uniqueness of $\tau$. In particular $\tau_{i/2}$ is unique 
and the uniqueness of
$R$ follows immediately from \rf{poldec}.
\end{Rem}

\begin{Rem}\label{rmk1.4} We have to point out that the automorphism 
$\tau_t$ was denoted in
\cite{mn} by $\tau_{t/2}$ and called the {\em deformation 
automorphism}. However, in this paper, we call it the scaling
automorphism following the terminology of \cite{w4}.
\end{Rem}

\begin{Rem}\label{rmk1.6} The strong right invariance is a condition 
which establishes the connection between
comultiplication and antipode$:$ it replaces the well known formula
\begin{equation}
\label{1941}
{\rm m}(\id\otimes\kappa)\delta(a)={\rm m}(\kappa\otimes\id)\delta(a)=e(a)1
\end{equation}
of the usual theory of Hopf algebras. This means the following.
\oldversion{If $A$ is a Hopf algebra, the mapping $W :A\otimes 
A\rightarrow A\otimes
A$ defined by $a\otimes c\mapsto\delta (a)(1\otimes c)$ for $a,\;c\in A$ is a
bijective linear mapping acting on the vector space $A\otimes A$ due to the
existence of the explicit inverse $W^{-1}:b\otimes d\mapsto
(\id_A\otimes\kappa)(\delta (b))(1\otimes d)$. The proof that those 
mappings are
inverse to each other uses the above formula. On the Hilbert space 
level $W$ should
be a unitary operator.  Therefore we require the condition $W^*=W^{-1}$, i.e.
\[
\IS{W(a\otimes c)}{b\otimes d} =\IS{a\otimes c}{W^{-1}(b\otimes d)}.
\]
In other words
\[
(h\otimes h)((1\otimes c^*)\delta(a^*)(b\otimes d)) =(h\otimes 
h)((a^*\otimes c^*)(\id\otimes\kappa)(\delta(b))(1\otimes
d)).
\]
  Then by putting $\varphi(x) =h(c^*xd)$, we get the strong right invariance.}
\newversion{If $A$ is a Hopf algebra, then the linear mapping $W :A\otimes
A\rightarrow A\otimes A$ defined by $a\otimes c\mapsto\delta 
(a)(1\otimes c)$ (where
$a,\;c\in A$) is bijective due to the existence of the explicit inverse
$V:b\otimes d\mapsto (\id_A\otimes\kappa)(\delta (b))(1\otimes d)$. The proof
that $V=W^{-1}$ uses \rf{1941}. On the Hilbert space level $W$ should unitary.
Therefore we require the equality $W^*=V$, i.e.
\[
(h\otimes h)((1\otimes c^*)\delta(a^*)(b\otimes d)) =
(h\otimes h)((a^*\otimes c^*)(\id\otimes\kappa)(\delta(b))(1\otimes
d)).
\]
This is the strong right invariance with $\varphi(x) =h(c^*xd)$.}

\end{Rem}

Some of our main assertions in this paper are summarized as follows:

\begin{Thm}
\label{th1.5}
  Let $(A,\delta)$ be a weighted Hopf $C^*$-algebra. Then
\begin{enumerate}
\item The unitary antipode $R$, the scaling group 
$\seq{\tau_t}{t\in\real}$ and the Haar weight
$h$ are uniquely determined.
\item  The coproduct $\delta$ commute with the unitary antipode and 
the scaling group $:$
\[
\delta\comp R=\sigma\comp(R\otimes R)\comp\delta, \quad
\delta\comp\tau_t=(\tau_t\otimes\tau_t)\comp\delta, \quad t\in\real,
\]
where $\sigma$ is the flip
automorphism on $A\otimes A$.
\item There exists a pointwise norm continuous one parameter 
$*$-automorphism group
$\seq{\sigma_t}{t\in\real}$ of $A$ such that
\[
h(\sigma_{i/2}(a)\sigma_{i/2}(a)^*)=h(a^*a)
\]
for any $a\in\Dom(\sigma_{i/2})$. The group 
$\seq{\sigma_t}{t\in\real}$ is the modular
automorphism related to the Haar weight $h$. Clearly, the above 
formula is the KMS-condition for $h$.
\item For any $t\in\real$ we have
\[
\delta\comp\sigma_t=(\sigma_t\otimes\tau_t)\comp\delta.
\]
\item There exists a strictly positive element $\rho$ affiliated with 
the $C^*$-algebra $A$ such that for any $a\in A_+$
  with $h(a)<\infty$ and any $\varphi\in A^*_+$ we have
\[
h(a*\varphi)=\varphi(\rho)h(a),
\]
where by definition $\varphi(\rho)={\displaystyle 
\lim_{n\rightarrow\infty}}\varphi(\frac{n\rho}{n1+\rho})$. Moreover 
we have $:$
$\delta(\rho)=\rho\otimes\rho$, $\tau_t(\rho)=\rho$ and $R(\rho)=\rho^{-1}$.
\item The composition $h^L=h\comp R$ is the left Haar weight on $A$. 
The Connes' Radon-Nikodym cocycle is of the form
\[
\left(Dh^L:Dh\right)_t=\lambda^{-it^2/2}\rho^{it},
\]
where $\lambda$ is the constant appearing in Definition 
\reef{def1.2}. Moreover $\sigma_t(\rho)=\lambda^{-t}\rho$.

\end{enumerate}
\end{Thm}

The strong right invariance plays an essential role in our theory. We 
end this Section with a number of equivalent
formulations of this condition. By the way we introduce Hilbert space 
objects that will be used throughout of the paper.
The most important one is the Kac-Takesaki operator $W$ (in other 
words the right regular representation of
the considered quantum group). In this (and the next) section the 
matrix elements of $W$ are defined and
their properties are investigated. The operator itself will be 
introduced in Section \reef{sek3}.\vs

Till the end of this section we shall deal with a proper 
C$^*$-bialgebra $(A,\delta)$ with the cancellation property
equipped with a faithful locally finite, lower semicontinuous 
right-invariant weight $h$ and with a closed
densely defined linear map $\kappa$ satisfying the first two 
conditions of Definition \reef{def1.2}. Therefore our
results will be independent of the strong right invariance condition and the
strict faithfulness of the weight $h$.\vs

  Using the GNS construction we obtain a Hilbert space $\Hil$ and a 
representation $\pi$ of $A$ acting on $\Hil$. Since $h$ is
faithful, so is $\pi$ and we shall identify $A$ with its $\pi$-image: 
$A\subset\skrL(\Hil)$. The corresponding GNS map
will be denoted by $\eta$. It is a densely defined closed linear 
mapping from $A$ into $\Hil$ with the domain
$\Dom(\eta)=\Nideal_h=\set{a\in A}{h(a^*a)<\infty}$ such that
\[
\IS{\eta(c)}{a\eta(b)}=h(c^*ab), \quad a\in A, \quad b,c\in\Dom(\eta),
\]
The GNS mapping $\eta$ plays a crucial role in the entire paper. It is closed
with respect to the strong topology on $A$ and the norm topology in $\Hil$. The
separability of $A$ implies that of the Hilbert space $\Hil$. See 
Appendices \reef{GNSmap} and \reef{B} for the details. \vs

We start with some technicalities which we use throughout this paper. 
For any $x,
y\in\Hil$, $\omega_{x,y}$ will denote the functional on $\skrL(\Hil)$ such that
\[
\omega_{x,y}(a)=\IS{x}{ay}
\]
for $a\in\skrL(\Hil)$. It is clear that the set of all such 
functionals is linearly dense in
the predual $\skrL(\Hil)_*$ of $\skrL(\Hil)$. In general, if $A$ is a
$C^*$-algebra acting on $\Hil$ then $A$ is endow with the $\sigma$-weak
topology coming from $\skrL(\Hil )$. In what follows $A_*$ will denote the set
of all $\sigma$-weakly continuous functionals on $A$. It is known that
\[
A_*=\set{\varphi |_A}{\varphi\in\skrL(\Hil)_*} .
\]
Consequently we set $A_{*+}=A_*\cap A_+$.
\vs

\oldversion{As it was shown after \rf{Kadison},}
\newversion{We already know (see the paragraph containing formula 
\rf{Kadison}) that}
\begin{equation}
\label{pr2.0}
\varphi*a\in\Dom(\eta)
\end{equation}
for any $a\in\Dom(\eta)$ and any $\varphi\in A^*$.

\begin{Prop}
\label{pr2.1}
For any $\varphi\in A^*$, there exists $W_{\varphi}\in{\mathcal
L}(\Hil)$ such that
\begin{equation}
\label{Ftran}
W_{\varphi}\eta (a)=\eta (\varphi *a)
\end{equation}
for any $a\in\Dom(\eta)$.  The mapping $A^*\ni \varphi\mapsto 
W_{\varphi}\in{\mathcal
L}(\Hil)$ is injective, linear and bounded.  Moreover,
\begin{equation}
\label{wproduct} W_{\varphi}W_{\psi}=W_{\varphi*\psi}
\end{equation}
for all $\varphi,\ \psi\in A^*$.
\end{Prop}

\begin{pf} Formula \rf{Ftran} defines $W_\varphi$ on a dense linear 
subset of $\Hil$. We have
to show that $W_\varphi$ is bounded. If $\varphi$ is positive then 
using \rf{Kadison} we get
$\norm{\eta(\varphi*a)}\leq\varphi(1)\norm{\eta(a)}$ and 
$\norm{W_\varphi}\leq\varphi(1)$. For general
$\varphi$ we have $\norm{W_\varphi}\leq 4\norm{\varphi}$. This is 
because any continuous linear
functional on $A$ is a linear combination of four positive 
functionals with norms not larger than the
norm of the original functional. \vs

If $W_\varphi=0$ then $\varphi*a=0$ for all $a\in A$. Hence 
$(\id\otimes\varphi)\left((b\otimes
1)\delta(a)\Vs{3}\right)=b(\varphi*a)=0$ for all $a,b\in A$ and using 
the cancellation property we get
$\varphi=0$. It shows that the mapping $\varphi\mapsto W_\varphi$ is 
injective. \vs

Formula (\reef{wproduct}) follows immediately from the associativity 
of the convolution product.
\end{pf}

\begin{Prop}\label{pr2.2} For every $x,\;y\in\Hil$, there exists a 
unique element $W(x,y)\in
A$ such that
\begin{equation}\label{wzorpr2.2}
\varphi (W(x,y))=\IS{x}{W_{\varphi}y}
\end{equation}
for any $\varphi\in A^*$. Furthermore, $W(x,y)$ is sesquilinear in 
$(x,y)$ and $\norm{W(x,y)}\leq 4\norm{x} \norm{y}$.
\end{Prop}

\begin{pf} The uniqueness, the sesquilinearity and the estimate 
follow from the formula. It is enough
to prove the existence of operators $W(x,y)$ for a dense set of 
vectors $x,\ y$. We shall use the commutant $\eta'$ of
the GNS map $\eta$ (see Appendix \reef{GNSmap}). For any 
$a\in\Dom(\eta)$, $a'\in\Dom(\eta')$ and
$x\in\Hil$ we put :
\[
W({a'}^*x,\eta(a)):=a*\omega_{x,\eta'(a')}.
\]
Then for any $\varphi\in A^*$ we have
\[
\begin{array}{r@{\;}c@{\;}l}
\varphi(W({a'}^*x,\eta(a)))&=&\varphi(a*\omega_{x,\eta'(a')}) 
=\IS{x}{(\varphi*a)\eta'(a')}\\
&=&\IS{x}{a'\eta(\varphi*a)}=\IS{{a'}^*x}{W_{\varphi}\eta(a)}\Vs{5},
\end{array}
\]
where the second equality due to the cyclic property of the 
convolution product. This proves the
assertion.
\end{pf}

\begin{Prop}\label{pr2.3} The set $\set{W(x,y)}{x,\;y\in\Hil}$ is 
linearly dense in $A$.
\end{Prop}

\begin{pf} Suppose $\varphi\in A^*$ satisfies $\varphi (W(x,y))=0$ 
for all $x,\; y\in \Hil$.
Then
$\IS{x}{\eta (\varphi *a)}=\IS{x}{W_{\varphi}\eta(a)}=\varphi 
(W(x,\eta (a)))=0$ and $\varphi *a=0$ for all
$a\in\Dom(\eta)$. Therefore$W_\varphi=0$ and hence $\varphi =0$. This 
proves the assertion.
\end{pf}\vs

\begin{Thm}
\label{SRIC} Let $(A,\delta)$ be a proper C$^*$-bialgebra with the 
cancellation property, $h$ be a
faithful locally finite, lower semicontinuous right invariant weight 
on $A$ and $\kappa$ be a
closed densely defined linear map acting on $A$. Assume that the 
first two conditions of Definition
\reef{def1.2} are satisfied. Then the following statements are 
equivalent $:$\vs

$1.$ The strong right invariance $($Condition $3$ of Definition 
\reef{def1.2}$)$.\vs

$2.$ For all $\varphi\in A^*$ such that $\varphi\comp\kappa\in A^*$ we have
\begin{equation}
W_{\varphi^*}=\left(W_{\varphi\comp\kappa}\right)^*. \label{wadjoint}
\end{equation}

$3.$ Formula \rf{wadjoint} holds for all $\varphi$ from a weakly${}^*$ dense
subset $\Dom_0$ of $A^*$ such that
$\psi\comp\tau_t\in\Dom_0$ for all $\psi\in\Dom_0$ and $t\in\real$.\vs

$4.$ For any $x,\;y\in\Hil$, the element $W(x,y)$ is in the domain of $\kappa$
and $\kappa (W(x,y))=W(y,x)^*$.

\end{Thm}

\begin{pf}
Assume that the strong right invariance holds. Then for any $\varphi\in A^*$
with $\varphi\comp\kappa\in A^*$ and any $a,b\in\Dom(\eta)$ we have
\[
\begin{array}{r@{\;}c@{\;}l}
\IS{W_{\varphi^*}\eta(a)}{\eta(b)}&=&\IS{\eta(\varphi^* 
*a)}{\eta(b)}=h((\varphi*a^*)b)\\
&=&h(a^*((\varphi\comp\kappa)*b))
=\IS{\eta(a)}{\eta((\varphi\comp\kappa)*b)}= 
\IS{\eta(a)}{W_{\varphi\comp\kappa}\eta(b)},\Vs{5}
\end{array}
\]
where we used the equality $\varphi *a^*=(\varphi^**a)^*$ for the 
second equality. Remembering that the range of $\eta$ is dense in 
$\Hil$, we obtain relation (\reef{wadjoint}). It shows that Statement
1 implies Statement 2.\vs

Clearly Statement 2 implies Statement 3. Assume that Statement 3 
holds. Then for any $x,y\in\Hil$ and any
$\varphi\in\Dom_0$ we have
\[
\begin{array}{r@{\;}c@{\;}l} (\varphi\comp\kappa)(W(x,y))
&=&(x|W_{\varphi\comp\kappa}y)=(W_{\varphi^*}x|y) \\ 
&=&\overline{(y|W_{\varphi^*}x)}
=\overline{\varphi^*(W(y,x))}=\varphi(W(y,x)^*)\Vs{5}.
\end{array}
\]
Inserting $\varphi\comp\tau_t$ instead of $\varphi$, multiplying the 
both sides by $e^{-t^2}$, integrating over
$t\in\real$ and remembering that $\kappa$ commutes with $\tau_t$, we obtain
\[
(\varphi\comp\kappa)\left(\Vs{3}\skrR_\tau(W(x,y)\right)=\varphi\left(\Vs{3}\skrR_\tau 
(W(y,x)^*)\right),
\]
where $\skrR_\tau$ is a linear operator acting on $A$ introduced by the formula
\[
\skrR_\tau(a)=\frac{1}{\sqrt{\pi}}\int_{-\infty}^{\infty}e^{-t^2}\tau_t(a)\,dt
\]
for any $a\in A$. It is known (cf. Appendix \reef{AnaGene}) that the 
range of $\skrR_\tau$ is contained in
$\Dom(\tau_{i/2})=\Dom(\kappa)$. Therefore
$(\varphi\comp\kappa)\left(\Vs{3}\skrR_\tau(W(x,y)\right)=\varphi\left(\kappa\left(\Vs{3}\skrR_\tau(W(x,y)\right)\right)$.
Remembering that $\varphi$ runs over a weakly$^*$ dense subset 
$\Dom_0\subset A^*$ we obtain
\[
\kappa\left(\Vs{3}\skrR_\tau(W(x,y))\right)=\skrR_\tau (W(y,x)^*).
\]
Theorem \reef{removeReg} shows now that $W(x,y)\in\Dom(\kappa)$ and 
that $\kappa(W(x,y))=W(y,x)^*$. This shows that Statement 3 implies 
Statement 4.\vs

To complete the proof we have to show that Statement 4 implies the 
strong right invariance. Let $a,b\in\Nideal_h$.
Inserting $x=\eta(a)$ and $y=\eta(b)$ in Statement 4, we see that 
$W(\eta(a),\eta(b))\in\Dom(\kappa)$ and
$\kappa (W(\eta(a),\eta(b)))=W(\eta(b),\eta(a))^*$. Therefore for any 
$\varphi\in A^*$ with $\varphi\comp\kappa\in A^*$
we have
\[
\begin{array}{r@{\;}c@{\;}l}
h((\varphi*a^*)b)&=&\IS{\eta(\varphi^**a)}{\eta(b)}
  =\overline{\IS{\eta(b)}{W_{\varphi^*}\eta(a)}}
  =\overline{\varphi^*(W(\eta(b),\eta(a)))}\\ \Vs{5}
&=&\varphi\left(W(\eta(b),\eta(a))\Vs{3}^*\right)
  =\varphi\left(\kappa (W(\eta (a),\eta (b)))\right) \\
&=&(\varphi\comp\kappa)\left(W(\eta(a),\eta(b))\Vs{3}\right)
=\IS{\eta(a)}{W_{\varphi\comp\kappa}\eta(b)}\\ \Vs{5}
&=&\IS{\eta(a)}{\eta((\varphi\comp\kappa)*b)}=h(a^*((\varphi\comp\kappa)*b)).
\end{array}
\]
\end{pf}

\section{Matrix elements of Kac-Takesaki operator}
\label{ex2}

In this Section we derive the simplest consequences of axioms listed 
in Definition \reef{def1.2}. Let $(A,\delta)$ be a
weighted Hopf C$^*$-algebra. We shall use the notation and results of the
previous Section. In particular all \oldversion{the}
four Statements of Theorem \reef{SRIC} hold. \vs

We assumed that the weight $h$ satisfies the relative invariance 
condition: $h\comp\tau_t=\lambda^t h$ for
$t\in\real$. Therefore $\Dom(\eta)$ is $\seq{\tau_t}{t\in\real}$ 
invariant and there exists a positive
self-adjoint operator $Q$ on $\Hil$ with $\ker (Q)=\{0\}$ such that
\begin{equation}
Q^{2it}\eta (a):=\lambda^{-t/2}\eta (\tau_t(a)), \quad 
a\in\Dom(\eta). \label{Qdef}
\end{equation}
In \cite{mn} the operator $Q^2$ \oldversion{is} \newversion{was} 
denoted by $H$. For
any $a\in A$ and $b\in\Dom(\eta)$ we have
\[
Q^{2it}a\eta(b)= Q^{2it}\eta(ab)=\lambda^{-t/2}\eta(\tau_t(ab))
=\lambda^{-t/2}\tau_t(a)\eta(\tau_t(b)) =\tau_t(a)Q^{2it}\eta(b).
\]
Therefore
\begin{equation}
\tau_t(a)=Q^{2it}aQ^{-2it} \label{scaling}
\end{equation}
for all $t\in \real$ and $a\in A$.\vs

  By Lemma \reef{lem:strictfaithful}, the strict faithfulness of $h$ 
implies that
the mapping $\eta (a)\mapsto\eta(a^*)$ (where $a$ runs over 
$\Dom(\eta )\cap\Dom(\eta )^*$) is a
closable conjugate linear involution. We denote by $S$ its closure. 
Let $F:=S^*$. Then there exist an involutive antiunitary $J$
and a strictly positive self-adjoint operator $\Delta$ such that 
$S=J\Delta^{1/2}$ and $F=J\Delta^{-1/2}$ by the fact that
$S^2=1$ on $\Dom(\Delta^{1/2})$.\vs

Taking into account \rf{Qdef} we see that the set 
$\set{\eta(a)}{a\in\Dom(\eta)\cap\Dom(\eta)^*}$
is a $\seq{Q^{it}}{t\in\real}$ - invariant core of $S$. Moreover for 
any $x=\eta(a)$ in this set we have
\[
Q^{2it}Sx=Q^{2it}\eta(a^*)=\lambda^{-t/2}\eta(\tau_t(a^*)),
\]
\[
SQ^{2it}x=\lambda^{-t/2}S\eta(\tau_t(a))=\lambda^{-t/2}\eta(\tau_t(a)^*)
\]
and the right hand sides of the above equalities coincide. Therefore 
$Q^{2it}SQ^{-2it}=S$ for any $t\in\real$. By the
uniqueness of the polar decomposition, we have
\begin{equation}
Q^{2it}\Delta Q^{-2it}=\Delta, \label{qdelta}
\end{equation}
\[
Q^{2it}J Q^{-2it}=J.
\]
The first equation means that $\Delta$ strongly commutes with $Q$. 
The second one is equivalent to
the formula
\begin{equation}
JQJ=Q^{-1}. \label{Q}
\end{equation}
The operators $Q$ and $\Delta$ play an essential role in all our 
considerations.

\begin{Lem}\label{lm2.4} Let $x\in\Dom(\Delta^{1/2})$ and 
$y\in\Dom(\Delta^{-1/2})$.
Then
$W(J\Delta^{-1/2}y,J\Delta^{1/2}x)=W(y,x)^*$.
\end{Lem}

\begin{pf} Since $\eta(\Dom(\eta)\cap\Dom(\eta)^*)$ is a core for 
$\Delta^{1/2}$, it
is sufficient to prove the statement for $x=\eta(a)$ where 
$a\in\Dom(\eta)\cap{\mathcal
D}(\eta)^*$. For any $\varphi\in A^*$, we have
\[
\begin{array}{r@{\;}c@{\;}l}
\varphi 
(W(J\Delta^{-1/2}y,J\Delta^{1/2}x))&=&(J\Delta^{-1/2}y|W_{\varphi}\eta 
(a^*))\\
    &=&(J\Delta^{-1/2}y|\eta (\varphi 
*a^*))=(J\Delta^{-1/2}y|J\Delta^{1/2}\eta (\varphi^**a))\Vs{5}\\
    &=&(\eta (\varphi^**a)|y)=\overline{(y|W_{\varphi^*}\eta (a))} \Vs{5} \\
    &=&\overline{\varphi^*(W(y,\eta (a)))}=\varphi (W(y,\eta (a))^*).\Vs{5}
\end{array}
\]

\end{pf}\vs

Combining Statement 4 of Theorem \reef{SRIC} with Lemma \reef{lm2.4}, we get
\begin{equation}
\kappa(W(x,y))=W(J\Delta^{-1/2}y,J\Delta^{1/2}x) \label{kappaw}
\end{equation}
for all $x\in\Dom(\Delta^{1/2})$ and $y\in\Dom(\Delta^{-1/2})$.

\begin{Prop}\label{pr2.8} For any $x,\; y\in\Hil$, we have
\begin{enumerate}
\item $\tau_t(W(x,y))=W(\Delta^{it}x,\Delta^{it}y)$ for
$t\in\real$.
\item $R(W(x,y))=W(Jy,Jx)$.
\end{enumerate}
\end{Prop}

\begin{Rem}
This proposition shows that the unitary antipode and the scaling 
group are uniquely determined by the Haar weight.
So is $\kappa$. Unfortunately the proof
of the uniqueness of the Haar weight requires the prior knowledge of 
the uniqueness of $\kappa$. See the proof of Theorem \reef{th7.3}.
\end{Rem}

\begin{pf} Let $x\in\Dom(\Delta)$ and $y\in\Dom(\Delta^{-1})$.  Then 
using twice
(\reef{kappaw}), we get
\[
\begin{array}{r@{\;}c@{\;}l} W(\Delta x,\Delta^{-1}y)
  &=&W(J\Delta^{-1/2}J\Delta^{1/2}x,J\Delta^{1/2}J\Delta^{-1/2}y) \\
  &=&\kappa (W(J\Delta^{-1/2}y,J\Delta^{1/2}x))\Vs{5} \\
  &=&\kappa^2(W(x,y))=\tau_i(W(x,y))\Vs{5}
\end{array}
\]
It shows that $W(x,y)\in\Dom(\tau_i)$ and
\[
\tau_i(W(x,y))=W(\Delta x,\Delta^{-1}y).
\]
Using now Theorem \reef{drugie}, we obtain the first statement. 
Moreover, for any $x\in\Dom(\Delta^{1/2} )$ and 
$y\in\Dom(\Delta^{-1/2})$, the element $W(x,y)\in{\mathcal
D}(\tau_{i/2})$ and
\[
\tau_{i/2}(W(x,y))=W(\Delta^{1/2}x,\Delta^{-1/2}y).
\]
Combining this with (\reef{kappaw}), we see that
$R(W(\Delta^{1/2}x,\Delta^{-1/2}y))=W(J\Delta^{-1/2}y,J\Delta^{1/2}x)$ 
for $x\in\Dom(\Delta^{1/2})$ and
$y\in\Dom(\Delta^{-1/2})$. Now the continuity of $W$ shows that the 
second statement holds in full generality.
\end{pf}\vs

By the above proposition the set $\set{W(x,y)}{x,y\in\Hil}$ is $R$ 
and $\seq{\tau_t}{t\in\real}$ invariant. By
Proposition \reef{pr2.3} this set is linearly dense in $A$.  Using Theorem
\reef{core} we obtain

\begin{Cor}
\label{pr2.6}
The linear span of the set $\set{\Vs{3}W(x,y)}{x,y\in\Hil}$ is a core 
for $\kappa$.
\end{Cor}

\begin{Prop}
\label{pr2.9}
For any $\varphi\in A^*$, we have
\begin{enumerate}
\item $W_{\varphi\comp R}=J(W_{\varphi})^*J$.
\item $W_{\varphi\comp\tau_t}=\Delta^{-it}W_{\varphi}\Delta^{it}$.
\end{enumerate}
\end{Prop}

\begin{pf} Ad 1. For $x,\;y\in\Hil$ and $\varphi\in A^*$, we have
\[
\begin{array}{r@{\;}c@{\;}l} (x|W_{\varphi\comp R}y)&=&\varphi\comp 
R(W(x,y))=\varphi (W(Jy,Jx))\\
&=&(Jy|W_{\varphi}Jx)=(JW_{\varphi}Jx|y)\;\Vs{5},
\end{array}
\]
where we used Proposition \reef{pr2.8} for the second equality. This 
proves the equality
$J(W_{\varphi})^*J=W_{\varphi\comp R}$.

Ad 2. Similarly, we have
\[
\begin{array}{r@{\;}c@{\;}l} 
(x|W_{\varphi\comp\tau_t}y)&=&\varphi\comp\tau_t(W(x,y))\\
&=&\varphi(W(\Delta^{it}x,\Delta^{it}y)) 
=(\Delta^{it}x|W_{\varphi}\Delta^{it}y)\Vs{5}.
\end{array}
\]
This proves the equality 
$\Delta^{-it}W_{\varphi}\Delta^{it}=W_{\varphi\comp\tau_t}$.
\end{pf}

\begin{Prop}\label{pr2.10}
\begin{enumerate}
\item $\delta\comp R=\sigma\comp (R\otimes R)\comp\delta$, where 
$\sigma$ is the flip automorphism on
$A\otimes A$.
\item $\delta\comp\tau_t=(\tau_t\otimes\tau_t)\comp\delta$.
\end{enumerate}
\end{Prop}

\begin{pf} Ad 1. By Proposition \reef{pr2.9}, for any $\varphi,\ 
\psi\in A^*$ we have
\[
\begin{array}{r@{\;}c@{\;}l} W_{(\varphi *\psi )\comp
R}&=&J(W_{\varphi*\psi})^*J=J(W_{\varphi}W_{\psi})^*J\\
&=&J(W_{\psi})^*JJ(W_{\varphi})^*J=W_{\psi\comp R}W_{\varphi\comp R}\Vs{5}.
\end{array}
\]
By Proposition \reef{pr2.1} , $W_{\varphi}=0$ implies $\varphi =0$. Therefore
$(\varphi *\psi)\comp R=(\psi\comp R)*(\varphi\comp R)$ and hence we obtain
$(\varphi\otimes\psi )\comp\delta\comp R=(\varphi\otimes\psi )
\comp\sigma\comp (R\otimes R)\comp\delta$ for $\varphi , \;\psi\in 
A^*$.  This proves Assertion 1.

Ad 2. By Proposition \reef{pr2.9}, we have
\[
\begin{array}{r@{\;}c@{\;}l}
W_{(\varphi\comp\tau_t)*(\psi\comp\tau_t)}&=&W_{\varphi\comp\tau_t}W_{\psi\comp\tau_t}
=(\Delta^{-it}W_{\varphi}\Delta^{it})(\Delta^{-it}W_{\psi}\Delta^{it})\\
&=&\Delta^{-it}W_{\varphi}W_{\psi}\Delta^{it} 
=\Delta^{-it}W_{\varphi*\psi}\Delta^{it}
=W_{(\varphi*\psi)\comp\tau_t}\Vs{5}.
\end{array}
\]
We then have
$(\varphi\comp\tau_t)*(\psi\comp\tau_t)=(\varphi*\psi)\comp\tau_t$ 
and hence we obtain
\[
\begin{array}{rcl} (\varphi\otimes\psi)((\tau_t\otimes\tau_t)\comp\delta)
&=&((\varphi\comp\tau_t)\otimes(\psi\comp\tau_t))\comp\delta\\
&=&(\varphi\comp\tau_t)*(\psi\comp\tau_t)=(\varphi*\psi)\comp\tau_t
=(\varphi\otimes\psi)\comp\delta\comp\tau_t\Vs{5}
\end{array}
\]
for all $\varphi,\ \psi\in A^*$.  This proves Assertion 2.
\end{pf}

\begin{Cor}\label{cor2.11}
$\delta\comp\kappa=\sigma\comp (\kappa\otimes\kappa)\comp\delta$ on 
$\Dom(\kappa)$.
\end{Cor}

\begin{pf} By Assertion 2 of Proposition \reef{pr2.10}, 
$a\in\Dom(\tau_{i/2})$ if and only if
$\delta(a)\in\Dom(\tau_{i/2}\otimes\tau_{i/2})$.  Since 
$\Dom(\kappa)=\Dom(\tau_{i/2})$ and
$\Dom(\kappa\otimes\kappa)=\Dom(\tau_{i/2}\otimes\tau_{i/2})$, our 
statement is immediate from
Proposition \reef{pr2.10}.
\end{pf}\vs

\begin{Lem}
\label{lm2.13}
\begin{enumerate}
\item $Q^{2it}W_{\varphi}Q^{-2it}=W_{\varphi\comp\tau_{-t}}$ for 
$\varphi\in A^*$ and $t\in \real$.
\item For any $x,\; y\in \Hil$, we have
$\tau_t(W(x,y))=W(Q^{2it}x,Q^{2it}y)$ for $t\in \real$.
\end{enumerate}
\end{Lem}

\begin{pf} Ad 1. Let $\varphi\in A^*$ and $a\in\Dom(\eta)$.  Then
\[
\begin{array}{r@{\;}c@{\;}l}
\tau_{-t}(\varphi *\tau_t(a))
    &=&\tau_{-t}\comp(\id_A\otimes\varphi )\comp\delta (\tau_t(a))  \\
    &=&(\id_A\otimes\varphi\comp\tau_t)\comp\delta (a)\Vs{5}
    =(\varphi\comp\tau_t)*a
\end{array}
\]
and applying $\eta$ to the both sides, we get
\[
Q^{-2it}W_{\varphi}Q^{2it}\eta (a)=W_{\varphi\comp\tau_t}\eta (a) \;.
\]

Ad 2. For $x,\;y\in \Hil$ and $\varphi\in A^*$, we have
\[
\begin{array}{r@{\;}c@{\;}l}
\varphi(\tau_t(W(x,y)))&=&(x|W_{\varphi\comp\tau_t}y)\\
    &=&(x|Q^{-2it}W_{\varphi}Q^{2it}y)
    =\varphi(W(Q^{2it}x,Q^{2it}y))\Vs{5}
\end{array}
\]
by Assertion 1.  This completes the proof.
\end{pf}

\begin{Prop}\label{pr2.14} Let $\varphi,\ \psi\in A^*$.  Assume that 
$W_{\varphi}^*=W_{\psi}$. Then
$\psi=\varphi^*\comp\kappa$.
\end{Prop}

\begin{pf} We have to show that $\psi(a)=\varphi^*(\kappa(a))$ for 
any $a\in\Dom(\kappa)$. By Corollary \reef{pr2.6} we may
assume that $a=W(x,y)$, where $x,\ y\in\Hil$. Now we compute
\[
\begin{array}{r@{\;}c@{\;}l}
\psi(W(x,y))&=&(x|W_{\psi}y)=(W_{\varphi}x|y)=\overline{(y|W_{\varphi}x)}=\overline{\varphi(W(y,x))}\\
&=&\varphi^*(W(y,x)^*)=\varphi^*(\kappa(W(x,y)))\Vs{5}.
\end{array}
\]
\end{pf}

\section{Pentagonal equation and manageability} \label{sek3}

In this section, we construct the Kac-Takesaki operator
$W$ and prove that this operator is a multiplicative unitary.  We also prove
that $W$ is manageable in the sense of~\cite{w4}.

We start with a lemma of Dini type.
\par
\begin{Lem}
\label{lm3.1}
Let $\seq{a_n}{n\in \natu }$ be an increasing sequence of elements in $A_+$.
\begin{enumerate}
\item If $a\in A$ and for any
$\psi\in A_+^*$, the sequence $\seq{\psi (a_n)}{n\in \natu }$ 
converges to $\psi (a)$, then the
sequence $\seq{a_n}{n\in \natu}$ converges in norm to $a$.\vs

\item  If for any $\omega\in A_+^*$, the sequence
$\seq{\omega(a_n)}{n\in \natu }$ converges to $\norm{\omega}$, then 
the sequence $\seq{a_n}{n\in \natu }$
converges strictly to the unit $1$ of $M(A)$.
\end{enumerate}
\end{Lem}

\begin{pf} Ad 1. We put $\Lambda :=\set{\psi\in 
A^*_+}{\norm{\psi}\leq 1}$. Then
$\Lambda$ is a compact subset of $A^*$ with respect to the
$\text{weak}^*$ topology. By the assumption, $\seq{\psi (a)-\psi 
(a_n)}{n\in \natu }$ is a decreasing
sequence converging to zero. By the Dini's Theorem this convergence 
is uniform on
$\Lambda$.  Therefore we have
\[
\norm{a-a_n}=\sup_{\psi\in\Lambda}|\psi (a)-\psi (a_n)|\rightarrow 0\;.
\]
as $n$ goes to infinity and the assertion follows.

Ad 2. Since $\omega(a_n)\leq\norm{\omega}$ for all $\omega\in A_+^*$,
$a_n$ satisfies $0\leq a_n\leq 1$. Let $b\in A$.  Then the sequence 
$\seq{b^*a_nb}{n\in\natu }$ is
increasing and for any $\psi\in A_+^*$, $\psi(b^*a_nb)=(b\psi 
b^*)(a_n)\rightarrow
\norm{b\psi b^*}=(b\psi b^*)(1)=\psi(b^*b)$ when $n\to\infty$.  Using 
Assertion 1 we see that
$b^*a_nb$ converges in norm to $b^*b$.  Therefore
\[
  \norm{(1-a_n)^{1/2}b}^2=\norm{b^*(1-a_n)b}=\norm{b^*a_nb-b^*b}\rightarrow 0.
\]
and
\[
\norm{(1-a_n)b}\leq
\norm{(1-a_n)^{1/2}} \norm{(1-a_n)^{1/2}b}\rightarrow 0,
\]
for $\norm{(1-a_n)^{1/2}}\leq 1$ for all $n$. It shows that $a_nb$ 
converges to $b$ in norm.\vs

Similarly, we get a proof that the sequence $\seq{ba_n}{n\in\natu}$ 
converges to $b$ in norm. This
proves the strict convergence of the sequence $\seq{a_n}{n\in\natu }$ 
to $1\in M(A)$.
\end{pf}\vs

The following proposition is a consequence of the right invariance of 
the Haar weight.

\begin{Prop}\label{pr3.2} Let $\seq{e_n}{n\in\natu }$ be an 
orthonormal basis of the Hilbert space
$\Hil$.  Then for any $x,\ y$ belonging to the range of the 
GNS-mapping $\eta$, we have
\begin{equation}
\sum_{n=1}^{\infty}W(e_n,x)^*W(e_n,y)=(x|y)1 \;, \label{ww1}
\end{equation}
where the series converges with respect to the strict topology on
$M(A)$ .
\end{Prop}

\begin{pf} Let $\omega\in A_+^*$ and 
$(\Hil_{\omega},\pi_{\omega},\Omega_{\omega})$ be the
GNS-triple associated with $\omega$.  To proceed our computation we 
choose an orthonormal basis
$\seq{\zeta_n}{n\in\natu }$ in $\Hil_{\omega}$.  Then for any 
$x\in\Hil$ we have
\[
\begin{array}{r@{\;}c@{\;}l}
\omega (W(e_n,x)^*W(e_n,x))
    &=&\norm{\pi_{\omega}(W(e_n,x))\Omega_{\omega}}^2
    ={\displaystyle \sum_{k=1}^{\infty}|(\zeta_k|
        \pi_{\omega}(W(e_n,x))\Omega_{\omega})|^2 }\\
    &=&{\displaystyle \sum_{k=1}^{\infty}
        |\varphi_k(W(e_n,x))|^2
    =\sum_{k=1}^{\infty}
        |(e_n|W_{\varphi_k}x)|^2}\;,
\end{array}
\]
where $\varphi_k$ are linear functionals on $A$ introduced by the formula
\[
\varphi_k(a):=(\zeta_k|\pi_{\omega}(a)\Omega_{\omega})
\]
for all $a\in A$ and $k\in\natu $.
Summing over
$n=1,2,\dots,N$ and setting $N\rightarrow \infty$, we obtain
\[
\lim_{N\to\infty}\omega \left(\sum_{n=1}^N W(e_n,x)^*W(e_n,x)\right)
    =\sum_{n=1}^{\infty}\sum_{k=1}^{\infty}|(e_n|W_{\varphi_k}x)|^2
    =\sum_{k=1}^{\infty}\norm{W_{\varphi_k}x}^2.
\]

Assume now that $x=\eta(a)$, where $a\in\Dom(\eta)$.  Then the right 
hand side of the above
relation equals to
\begin{eqnarray*}
\sum_{k=1}^{\infty}\norm{W_{\varphi_k}x}^2
    &=&\sum_{k=1}^{\infty}\norm{\eta (\varphi_k*a)}^2
    =\sum_{k=1}^{\infty}h((\varphi_k*a)^*(\varphi_k*a)) \\
    &=&\lim_{K\to\infty}h\left(\sum_{k=1}^K(\varphi_k*a)^*(\varphi_k*a)\right).
\end{eqnarray*} It turns out below (see the last part of the proof) that
\begin{equation}
\left\norm{\sum_{k=1}^K(\varphi_k*a)^*(\varphi_k*a)-\omega*(a^*a)\right} 
\to 0 \quad \mbox{as} \quad
K\to\infty. \label{3.2}
\end{equation}
Using the lower semicontinuity and the right invariance of $h$, we obtain
\[
\begin{array}{r@{\;}c@{\;}l}
{\displaystyle \lim_{N\to\infty}\omega\left(\sum_{n=1}^N 
W(e_n,x)^*W(e_n,x)\right)} &=& h(\omega*(a^*a))\\
&=&\norm{\omega} h(a^*a)=\norm{\omega}(x|x).
\end{array}
\]
This result holds for any $\omega\in A_+^*$. By Assertion 2 of Lemma 
\reef{lm3.1} we
see that the series\linebreak $\sum_{n=1}^{\infty}W(e_n,x)^*W(e_n,x)$ 
is strictly converging to
$(x|x)1$.  To obtain
$(\reef{ww1})$ in full generality, we apply the polarization argument.\vs

To end the proof we have to show (\reef{3.2}).  Let $\psi\in A_+^*$,
$(\Hil_{\psi},\pi_{\psi},\Omega_{\psi})$ be the GNS-triple associated 
with $\psi$ and
$\seq{\xi_m}{m\in\natu }$ be an orthonormal basis in $\Hil_\psi$.  Then
\begin{eqnarray*}
\psi(\omega *(a^*a)) &=&(\psi\otimes\omega)(\delta(a^*a))
=\norm{(\pi_{\psi}\otimes\pi_{\omega})(\delta(a))(\Omega_{\psi}\otimes\Omega_{\omega})}^2\\
&=&\sum_{m,k=1}^{\infty}|(\xi_m\otimes\zeta_k
|(\pi_{\psi}\otimes\pi_{\omega})(\delta(a))(\Omega_{\psi}\otimes\Omega_{\omega}))|^2\\
&=&\sum_{m,k=1}^{\infty}|(\xi_m 
|\pi_{\psi}((\id\otimes\varphi_k)(\delta(a))\Omega_{\psi})|^2
=\sum_{k=1}^{\infty}\norm{\pi_{\psi}(\varphi_k *a)\Omega_{\psi}}^2\\
&=&\sum_{k=1}^{\infty}\psi((\varphi_k*a)^*(\varphi_k* a))
=\lim_{K\to\infty}\psi\left(\sum_{k=1}^K(\varphi_k*a)^*(\varphi_k* a)\right).
\end{eqnarray*} Now $(\reef{3.2})$ follows from Assertion $1$ of 
Lemma \reef{lm3.1}.
\end{pf}

\begin{Prop}
\label{pr3.3}
Let $\seq{e_n}{n\in\natu }$ be an orthonormal basis of the Hilbert space
$\Hil$.  Then for any $x,\ y$ such that $Jx$ and $Jy$ belong to the 
range of $\eta$, we have
\begin{equation}
\sum_{n=1}^{\infty}W(x,e_n)W(y,e_n)^*=(x|y)1  \label{ww2}
\end{equation}
with respect to the strict topology.
\end{Prop}

\begin{pf} Replacing in the previous Proposition $\seq{e_n}{n\in\natu }$,
$x$ and $y$ by $\seq{Je_n}{n\in\natu }$, $Jy$ and $Jx$ respectively, we obtain
\[
\sum_{n=1}^{\infty}W(Je_n,Jy)^*W(Je_n,Jx)=(Jy|Jx)1=(x|y)1.
\]
The unitary antipode $R$ can be extended to the multiplier algebra
$M(A)$ and the extension denoted by the same symbol is continuous 
with respect to the strict
topology.  Remembering that the above series is strictly converging, we have
\[
\sum_{n=1}^{\infty}R(W(Je_n,Jx))R(W(Je_n,Jy))^*=(x|y)1.
\]
Thus (\reef{ww2}) follows. (cf. Second formula of Proposition \reef{pr2.8})
\end{pf}

\begin{Rem}\label{rmk3.4} In what follows, relations \rf{ww1} and 
\rf{ww2} will be used to
prove the unitarity of the Kac-Takesaki operator. For this purpose, 
it is sufficient to know that the
series
$(\ref{ww1})$ and $(\ref{ww2})$ converge weakly. However to prove 
$(\ref{ww2})$ by using
$(\ref{ww1})$ we have to use strict topology $;$ at the moment the 
continuity of $R$ with respect to
the weak operator topology is not established yet.
\end{Rem}

The Kac-Takesaki operator $W$ is introduced by the following theorem. 
It plays a role of the right
regular representation for the "quantum group" associated with the 
weighted Hopf $C^*$-algebra. We
recall that, for any $x,\ y\in\Hil$, $\omega_{x,y}$ denotes the 
linear functional on
$\skrL(\Hil)$ defined by $\omega_{x,y}(a)=(x|ay)$.

\begin{Thm}
\label{th3.5}
There exists a unique unitary operator $W$ acting on 
$\Hil\otimes\Hil$ such that
\begin{equation}
(x\otimes z|W(y\otimes u))=(z|W(x,y)u),
\label{KT}
\end{equation}
\begin{equation}
(\omega_{x,y}\otimes\id)(W)=W(x,y), \label{matel}
\end{equation}
\begin{equation}
\label{lslice}
(\id\otimes\varphi)W=W_\varphi
\end{equation}
for any $x, y, z, u\in\Hil$ and $\varphi\in{\mathcal L}(\Hil)_*$. On 
the right hand side of the last relation $\varphi\in
A^*$ denotes the restriction of $\varphi\in{\mathcal L}(\Hil)_*$ to 
$A$. Each of the above equations determines $W$ uniquely.
\end{Thm}

\begin{pf} Clearly \rf{KT} and \rf{matel} are equivalent. The computation
\[
\ITS{z}{(\id\otimes \omega_{y,x})(W)}{u}=(z\otimes y|W(u\otimes x))
=(y|W(z,u)x)=(z|W_{\omega_{y,x}}u)
\]
shows the equivalence of \rf{KT} and \rf{lslice}. The uniqueness of 
$W$ is obvious.\vs

  Let $D_0$ be the range of the GNS-mapping $\eta$ and $D_1=JD_0$. 
Then $D_0$ and $D_1$ are dense linear subsets of
$\mathcal H$.  We choose an orthonormal basis $\seq{e_n}{n\in\natu }$ 
in $\mathcal H$.  For any $y\in D_0$, $x\in D_1$,
$z, u\in{A\mathcal H}$ we set
\[
W(y\otimes u):=\sum_{n=1}^{\infty}e_n\otimes W(e_n,y)u
\]
and
\[
W^{\dagger}(x\otimes z):=\sum_{n=1}^{\infty}e_n\otimes W(x,e_n)^*z.
\]
Using $(\ref{ww1})$ and $(\ref{ww2})$ one can easily show that the 
series on the right hand sides
are norm converging in $\Hil\otimes\Hil$. Clearly the above formulae 
introduce linear
mappings $W:D_0\otimes_{alg}{A\mathcal H}\to\Hil\otimes\Hil$ and
$W^{\dagger}:D_1\otimes_{alg}{A\mathcal H}\to\Hil\otimes\Hil$. Using 
$(\ref{ww1})$
and $(\ref{ww2})$ once more, we can show that these mappings are 
isometries.  Extending these
mappings by continuity, we obtain isometric operators $W$ and 
$W^{\dagger}$ defined on the whole
$\Hil\otimes\Hil$.

To show that $W$ is unitary, it is sufficient to prove that the 
adjoint $W^*$ is an isometry.  We
claim that $W^*=W^{\dagger}$. Indeed, for any $y\in D_0$, $x\in D_1$, 
$z, u\in{A\mathcal H}$, we have
\begin{eqnarray*} (x\otimes z|W(y\otimes u)) 
&=&\sum_{n=1}^{\infty}(x\otimes z|e_n\otimes W(e_n,y)u)
=\sum_{n=1}^{\infty}(x|e_n)(z|W(e_n,y)u)\\ 
&=&\IS{z}{W\left(\sum_{n=1}^{\infty}e_n(e_n|x),y\right)u}
=(z|W(x,y)u)
\end{eqnarray*} and
\begin{eqnarray*} (W^{\dagger}(x\otimes z)|y\otimes u) 
&=&\sum_{n=1}^{\infty}(e_n\otimes
W(x,e_n)^*z|y\otimes u)\\ &=&\sum_{n=1}^{\infty}(e_n|y)(z|W(x,e_n)u)\\
&=&\IS{z}{W\left(x,\sum_{n=1}^{\infty}e_n(e_n|y)\right)u}=(z|W(x,y)u).
\end{eqnarray*} By the way we proved $(\ref{KT})$.
\end{pf}\vs

We shall use the leg numbering notation \cite{bs}, for example:
\[
(x_1\otimes x_2\otimes x_3|W_{13}(y_1\otimes y_2\otimes y_3)) 
=(x_1\otimes x_3|W(y_1\otimes
y_3))(x_2|y_2)
\]
for $x_j,\ y_j\in\Hil$ for $j=1,\ 2,\ 3$. We shall also use an exact 
vector presentation (see
Appendix \reef{GNSmap}) of the GNS map $\eta$:
\[
\eta(a)=\sum_{n\in\natu }a\Omega_n
\]
for any $a\in \Dom(\eta)$, where $\Omega_n\in\Hil$ for $n=1,2,\dots$, 
$A\Omega_n\perp A\Omega_m$ whenever
$n\neq m$ and $\Dom(\eta)$ coincides with the set of all $a\in A$ 
such that the above series is
convergent. Moreover $\sum_{n\in\natu}A\Omega_n$ is norm dense in 
$\Hil$, this is because the range
of $\eta$ is dense.\vs

In the following Proposition we collect basic properties of the 
Kac-Takesaki operator.

\begin{Prop}
\label{pr3.6}\Vs{2}

\begin{enumerate}
\item The coproduct $\delta$ is implemented by $W$ $:$
\begin{equation}
\delta(a)=W(a\otimes 1)W^* \label{coproduct}
\end{equation}
for any $a\in A$.
\item The pentagonal equation $W_{12}W_{13}W_{23}=W_{23}W_{12}$ 
holds. In other words
\begin{equation}
\label{character}
(\id\otimes\delta)W=W_{12}W_{13} \;.
\end{equation}
\item The operator $W$ commutes with the operators $Q\otimes Q$ and 
$\Delta\otimes Q^2$:
\begin{equation}
\label{pr3.64a}
W^*(Q\otimes Q)W=Q\otimes Q \;,
\end{equation}
\begin{equation}
\label{pr3.64b}
W^*(\Delta\otimes Q^2)W=\Delta\otimes Q^2 \;.
\end{equation}
\end{enumerate}
\end{Prop}

\begin{pf} \Vs{2}

Ad 1.  For any $x,\ y,\ z\in\Hil$ and $a\in\Dom(\eta)$ we have
\[
\begin{array}{r@{\;}c@{\;}l}
(x\otimes z|W(\eta(a)\otimes y))
&=&(x|(\id\otimes\omega_{z,y})(W)\eta(a))=(x|W_{\omega_{z,y}}\eta(a))\\ \Vs{7}
&=&(x|\eta(\omega_{z,y}*a))={\displaystyle 
\sum_{n\in\natu}(x|(\omega_{z,y}*a)\Omega_n)}.
\end{array}
\]
Let $p_n$ be the projection onto the closure of $A\Omega_n$. Then 
$p_n\in A'$ and $\sum_{n\in\natu
}p_n=1$. Replacing $x$ by $p_nx$, we obtain:
\[
\begin{array}{r@{\;}c@{\;}l}
(x\otimes z|(p_n\otimes 1)W(\eta(a)\otimes y))
&=&(x|(\omega_{z,y}*a)\Omega_n)\\ 
&=&(x|((\id\otimes\omega_{z,y})(\delta(a)))\Omega_n) =(x\otimes
z|\delta(a)(\Omega_n\otimes y))\Vs{5}.
\end{array}
\]
Therefore \[(p_n\otimes 1)W(\eta(a)\otimes 
y)=\delta(a)(\Omega_n\otimes y).\] It shows that
\begin{equation}
W(\eta(a)\otimes y)=\sum_{n\in\natu }\delta(a)(\Omega_n\otimes y),
\label{KT2}
\end{equation}
where the right hand side converges in norm.

Now for any $a\in A$, $b\in\Dom(\eta)$ and $y\in\Hil$ we have
\[
\begin{array}{r@{\;}c@{\;}l}
\delta(a)W(\eta(b)\otimes y)&=& {\displaystyle 
\delta(a)\sum_{n\in\natu}\delta(b)(\Omega_n\otimes
y)}\\ \Vs{6} &=&{\displaystyle 
\sum_{n\in\natu}\delta(ab)(\Omega_n\otimes y)} = W(\eta(ab)\otimes
y)=W(a\otimes 1)(\eta(b)\otimes y).
\end{array}
\]
Therefore $\delta(a)W=W(a\otimes 1)$ and $\delta(a)=W(a\otimes 1)W^*.$\vs

Ad 2.
Let $\varphi,\ \psi\in\skrL(\Hil)_*$.  For any $a\in\skrL(\Hil)$
we set
\[
\rho(a)=(\varphi\otimes\psi)(W(a\otimes 1)W^*).
\]
Then $\rho\in\skrL(\Hil)_*$.  Moreover
\begin{equation}
(\id\otimes\rho)(W)=(\id\otimes\varphi\otimes\psi)(W_{23}W_{12}W_{23}^*). 
\label{KT3}
\end{equation}
If $a\in A$, then by the previous assertion
$\rho(a)=(\varphi\otimes\psi)(\delta(a))=(\varphi*\psi)(a)$.  It 
shows that the restriction of
$\rho$ to $A$ coincides with $\varphi*\psi$.  Using the first assertion we have
\[
(\id\otimes\rho)(W)=W_{\varphi*\psi}=W_{\varphi}W_{\psi}
=(\id\otimes\varphi\otimes\psi)(W_{12}W_{13}).
\]
Comparing this with $(\ref{KT3})$, we see that
$W_{23}W_{12}W_{23}^*=W_{12}W_{13}$ and the pentagonal equation 
follows. Taking into account
\rf{coproduct}, we obtain \rf{character}.\vs

Ad 3. Let $\varphi\in\skrL(\Hil)_*$, $t\in\real$ and
$\varphi_t(a):=\varphi(Q^{2it}aQ^{-2it})$ for any $a\in\skrL(\Hil)$. By
\rf{scaling}, $\varphi_t$ restricted to $A$ coincides with 
$\varphi\comp\tau_t$.  Therefore
\[
(\id\otimes\varphi)((1\otimes Q^{2it})W(1\otimes
Q^{-2it}))=(\id\otimes\varphi_t)(W)=W_{\varphi\comp\tau_t}.
\]
Using now Proposition \reef{pr2.9} and Lemma \reef{lm2.13}, we obtain:
\[
W_{\varphi\comp\tau_t}=\Delta^{-it}W_{\varphi}\Delta^{it}=(\id\otimes\varphi)((\Delta^{-it}\otimes
1)W(\Delta^{it}\otimes 1))
\]
\[
W_{\varphi\comp\tau_t}=Q^{-2it}W_{\varphi}Q^{2it}=(\id\otimes\varphi)((Q^{-2it}\otimes
1)W(Q^{2it}\otimes 1)).
\]
Comparing these with the previous relation, we get
\[
\begin{array}{r@{\;}c@{\;}l} (1\otimes Q^{2it})W(1\otimes Q^{-2it}) 
&=&(\Delta^{-it}\otimes
1)W(\Delta^{it}\otimes 1)\\ &=&(Q^{-2it}\otimes 1)W(Q^{2it}\otimes 1).\Vs{5}
\end{array}
\]
It shows that $\Delta^{it}\otimes Q^{2it}$ and $Q^{2it}\otimes 
Q^{2it}$ commute with $W$ and our
statement follows.
\end{pf}\vs

Combining \rf{wadjoint} with \rf{lslice} we obtain
\begin{equation}
\label{SLWx}
(\id\otimes\varphi\comp\kappa)W=(\id\otimes\varphi)(W^*)
\end{equation}
for any $\varphi\in A^*$ such that $\varphi\comp\kappa\in A^*$.\vs

Let us recall the basic setup of the theory of multiplicative 
unitaries discussed in~\cite{bs,w4}.
Let $\Hil$ be a separable Hilbert space. A unitary operator $W$ 
acting on $\Hil\otimes \Hil$ is called a {\it
multiplicative unitary} if it satisfies the pentagonal equation 
$W_{23}W_{12}=W_{12}W_{13}W_{23}$. Let $J$ be a conjugate
linear bijection $\Hil\ni x \mapsto Jx\in\Hil$ with $(Jx|Jy)=(y|x)$. 
A multiplicative unitary $W$ is said to be {\it
manageable} if there exist a strictly positive self-adjoint operator 
$Q$ acting on $\Hil$ and a unitary operator
$\widetilde W$ acting on $\Hil \otimes \Hil$ such that
\begin{enumerate}
\item $W^*(Q\otimes Q)W=Q\otimes Q$,
\item $(x_1\otimes x_2|W(y_1\otimes y_2))= (Jy_1\otimes 
Qx_2|{\widetilde W}(Jx_1\otimes Q^{-1}y_2))$
for all $x_1,\;y_1\in\Hil$,
$x_2\in\Dom(Q)$ and $y_2\in\Dom(Q^{-1})$.
\end{enumerate} Clearly the definition and the operator $Q$ are 
independent of the choice of $J$.

\begin{Prop}\label{pr3.7} The Kac-Takesaki operator $W$ is manageable.
\end{Prop}

\begin{pf} Let $J$ be the involutive antiunitary introduced in the 
paragraph preceding Lemma
\reef{lm2.4}. We shall prove that the above two conditions are 
satisfied by the operator $Q$
introduced in
$(\ref{Qdef})$ and ${\widetilde W}=W^*$.  We already know that $Q$ is 
strictly positive and that
$Q\otimes Q$ commutes with $W$.

Let $x, y\in\Hil$.  By Proposition \reef{pr2.6}, $W(x,y)\in\Dom(\kappa)$ and
$\kappa(W(x,y))=W(y,x)^*$.  Taking into account the definition of 
$\kappa$ and using the second
formula of Proposition \reef{pr2.8}, we see that 
$W(x,y)\in\Dom(\tau_{i/2})$ and
\[
\tau_{i/2}(W(x,y))=R(W(y,x)^*)=W(Jx,Jy)^*.
\]
Let $z, u\in\Hil$.  Then for any $t\in\real$ we have:
\[
\begin{array}{r@{\;}c@{\;}l} 
(Q^{2it}z|\tau_t(W(x,y))Q^{2it}u)&=&(z|Q^{-2it}\tau_t(W(x,y))Q^{2it}u)\\
&=&(z|W(x,y)u).\Vs{5}
\end{array}
\]
Assume that $z\in\Dom(Q)$ and $u\in\Dom(Q^{-1})$. Making the 
holomorphic continuation up to the
point $t=i/2$, we obtain
\[
(Qz|W(Jx,Jy)^*Q^{-1}u)=(z|W(x,y)u).
\]
Therefore \xxx{Now this formula is numbered; all nonumber are removed}
\begin{equation}
\label{manageable}
\begin{array}{r@{\;}c@{\;}l} (x\otimes z|W(y\otimes u))
&=&(z|W(x,y)u)=\overline{(Q^{-1}u|W(Jx,Jy)Qz)}\\ &=&\overline{(Jx\otimes
Q^{-1}u|W(Jy\otimes Qz))}=(Jy\otimes Qz|W^*(Jx\otimes Q^{-1}u)). \Vs{6}
\end{array}
\end{equation}
\end{pf}\vs

Let us recall the properties of manageable multiplicative unitaries 
obtained in \cite{bs,w4}.

\begin{Thm}
\label{th3.8}
Let $\Hil$ be a Hilbert space and $W$ be a manageable multiplicative unitary
on
$\Hil\otimes\Hil$. We set$:$
\begin{equation}
\label{newA}
A=\set{\Vs{4}(\varphi\otimes{\id})(W)}{\varphi\in\skrL(\Hil)_*}^{\rm 
norm\ closure},
\end{equation}
\begin{equation}
\label{newAhat}
{\hA}=\set{\Vs{4}(\id\otimes{\varphi})(W^*)}{\varphi\in\skrL(\Hil)_*}^{\rm 
norm\ closure}.
\end{equation}
Then
\begin{enumerate}
\item $A$ and $\hA$ are $C^*$-algebras acting non degenerately on $\mathcal H$.
\item The operator $W$ is an element of the multiplier algebra 
$M(\hA\otimes A)$.
\item There exists a unique $\delta\in\Mor(A,A\otimes A)$ such that
$({\id}\otimes\delta)(W)=W_{12}W_{13}$.  The pair $(A,\delta)$ is a 
proper C$^*$-bialgebra with the cancellation
property. The comultiplication is given by the formula
\begin{equation}
\delta(a)=W(a\otimes 1)W^* \label{copro}
\end{equation}
for all $a\in A$.
\item There exists a unique closed linear mapping $\kappa:A\to A$ 
$(A$ is treated as a Banach space$)$ for which the set
$\set{(\varphi\otimes{\id})(W)}{\varphi\in\skrL(\Hil)_*}$ is a core and
\begin{equation}
\kappa((\varphi\otimes{\id})(W))=(\varphi\otimes{\id})(W^*) \label{wantipode}
\end{equation}
for any $\varphi\in\skrL(\Hil)_*$. Furthermore the set $\Dom(\kappa)$ is a
subalgebra of $A$ and the mapping $\kappa:\Dom(\kappa)\to A$ is 
antimultiplicative. Moreover,
$\kappa(\Dom(\kappa))=\set{a}{a^*\in\Dom(\kappa)}$ and
$\kappa(\kappa(a)^*)^*=a$ for $a\in\Dom(\kappa)$.
\item The operator $\kappa$ admits the following polar decomposition $:$
\[
\kappa=R\comp\tau_{i/2},
\]
where $\tau_{i/2}$ is the analytic generator of a one parameter group
$\seq{\tau_t}{t\in\real}$ of automorphisms of the $C^*$-algebra $A$ and
$R$ is an involutive antiautomorphism of $A$ commuting with $\tau_t$ 
for all $t\in\real$.  In
particular, $\Dom(\kappa)=\Dom(\tau_{i/2})$. In this case $R$ and
$\seq{\tau_t}{t\in\real}$ are uniquely determined.
\item We have
\[\delta\comp\tau_t=(\tau_t\otimes\tau_t)\comp\delta \quad for\quad 
t\in{\Bbb R}; \quad \delta\comp
R=\sigma\comp(R\otimes R)\comp\delta.\]
\item $R$ is normal, i.e. continuous with respect to the 
$\sigma$-weak topology on $A$.
\item Let $\widetilde W$ and $Q$ be the operators in the definition 
of manageability for $W$. Then
they satisfy $\tau_t(a)=Q^{2it}aQ^{-2it}$ for $t\in\real$ and 
$W^{\top\otimes R}={\widetilde W}^*$,
where $b^\top=Jb^*J$.
\end{enumerate}
\end{Thm}

We shall apply this theorem to the Kac-Takesaki operator $W$ 
introduced in Theorem \reef{th3.5}. We
know that the set $\set{\omega_{x,y}}{x,y\in\Hil}$ is linearly dense 
in $\skrL({\mathcal
H})_*$. Inserting in (\reef{newA}) $\varphi=\omega_{x,y}$ and using 
(\reef{matel}) we obtain
\begin{equation}
\label{newA1}
A=\set{\Vs{4}W(x,y)}{x,y\in\Hil}^{\rm closed\ linear\ envelope}.
\end{equation}
Now Proposition \reef{pr2.3} shows that the algebra (\reef{newA}) 
coincides with the original
algebra $A$ that we started with. Comparing (\reef{copro}) with 
(\reef{coproduct}), we see that new
$\delta$ coincides with the original one.

Inserting in (\reef{wantipode}) $\varphi=\omega_{x,y}$ for 
$x,y\in\Hil$, we obtain
\[
\begin{array}{r@{\;}c@{\;}l}
\kappa(W(x,y))&=&\kappa((\omega_{x,y}\otimes\id)(W))=(\omega_{x,y}\otimes\id)(W^*)\\
&=&(\omega_{y,x}\otimes\id)(W)^*=W(y,x)^*.\Vs{5}
\end{array}
\]
Comparing now Statement 4 of Theorem \reef{SRIC} and Corollary 
\reef{pr2.6} with Assertion 4 of the above theorem, we see that new
$\kappa$ coincides with the old one (they have the same core on which 
they act in the same way). By
the uniqueness of the polar decomposition, we see that the new $R$ 
and $\tau$ coincide with the old
ones.\vs

\par Now we are ready to construct the dual weighted Hopf 
$C^*$-algebra $(\hA,\hdelta)$
associated with our weighted Hopf $C^*$-algebra $(A,\delta)$. Let 
$\Sigma\in{\mathcal
L}(\Hil\otimes\Hil)$ be the flip map : $\Sigma(x\otimes y)=y\otimes 
x$ for all $x,\ y\in\Hil$.
According to the general theory \cite{bs,w4}, the operator
\begin{equation}
\hW=\Sigma W^*\Sigma \label{sW}
\end{equation}
is a manageable multiplicative unitary.  We shall apply Theorem \reef{th3.8} to
$\hW$.  Clearly, replacing $W$ by $\hW$, we interchange the roles of 
$A$ and $\hA$. The
comultiplication $\hdelta\in\Mor(\hA,\hA\otimes\hA)$ related to $\hW$ 
satisfies the relation
$(\id\otimes\hdelta)\hW=\hW_{12}\hW_{13}$ and is given by the formula 
$\hdelta(b)=\hW(b\otimes 1)\hW^*$ where $b$
runs over $\hA$. Using (\reef{sW}), we obtain
\begin{equation}
\label{deltahat}
(\hdelta\otimes\id)W^*=W^*_{13}W^*_{23}
\end{equation}
and
\[
\hdelta(b)=\Sigma W^*(1\otimes b)W\Sigma.
\]

\Vs{4}By Theorem \reef{th3.8}, $(\hA,\hdelta)$ is a proper 
C$^*$-bialgebra with the cancellation property. It will be
proved in Section \reef{sek7} that $(\hA,\hdelta)$ is a weighted Hopf
$C^*$-algebra.  We say that $(\hA,\hdelta)$ is the dual of $(A,\delta)$.\vs

Assertion 4 of Theorem \reef{th3.8} provide us with the closed linear 
mapping $\hkappa$ acting on
$\hA$ such that the set 
$\set{(\id\otimes\varphi)(W^*)}{\varphi\in\skrL(\Hil)_*}$ is core for
$\hkappa$ and
\[
\hkappa((\id\otimes\varphi)(W^*))=(\id\otimes\varphi)(W).
\]
Clearly, $\hkappa$ is the antipode related to $(\hA,\hdelta)$. By Assertion 5,
\[
\hkappa=\hR\comp\htau_{i/2},
\]
where $\seq{\htau_t}{t\in\real}$ is a one parameter group of 
automorphisms of $\hA$ and $\hR$ is
an involutive antiautomorphism of $\hA$ which commutes with $\htau_t$ 
for $t\in\real$. Clearly,
$\seq{\htau_t}{t\in\real}$ and $\hR$ are the scaling group and the 
unitary antipode related to
$(\hA,\hdelta)$.\vs

According to Statement 3 of Proposition 1.4 of \cite{w4}, $\hW$ is 
manageable. Inspecting the proof
of this statement, we see that the operators $\hQ$ and 
$\widetilde{\hW}$ entering the definition of
the manageability are given by the formulae
\[
\hQ=Q, \quad \widetilde{\hW}=(\Sigma\widetilde{W}^*\Sigma)^{\top\otimes\top}.
\]
We know (cf. the proof of Proposition \reef{pr3.7}) that 
$\widetilde{W}=W^*$. Therefore
\[
\widetilde{\hW}=(\Sigma W\Sigma)^{\top\otimes\top}=(J\otimes J)\Sigma 
W^*\Sigma(J\otimes J).
\] Using now Assertion 8 of Theorem \reef{th3.8}, we obtain
\begin{equation}
\htau_t(b)=Q^{2it}bQ^{-2it} \label{htauQ}
\end{equation}
for any $t\in\real$ and $b\in\hA$ and
\[
\hW^{\top\otimes\hR}=(J\otimes J)\Sigma W\Sigma(J\otimes J).
\]
Applying the flip and *-conjugation to both sides, we obtain
\[ W^{\hR\otimes\top}=(J\otimes J)W^*(J\otimes J)=W^{\top\otimes\top}.
\]
Therefore
\begin{equation}
W^{\hR\comp\top\otimes\id}=W. \label{hRT}
\end{equation}

Let $(A,\delta)$ be a C$^*$-bialgebra. In concrete cases one of the 
difficult task is to show that a given
weight $h$ on $A$ is right invariant and that the strong right 
invariance holds. To this end we shall use the following

\begin{Thm}
\label{newthm}
Let $(A,\delta)$ be a proper C$^*$-bialgebra with the cancellation 
property, $\kappa$ be a closed operator acting on the
Banach space $A$ admitting the polar decomposition
\[
\kappa=R\comp\tau_{i/2},
\]
where $\tau_{i/2}$ is the analytic generator of a one parameter group
$\seq{\tau_t}{t\in\real}$ of automorphisms of the $C^*$-algebra $A$ and
$R$ is an involutive antiautomorphism of $A$ commuting with $\tau_t$ 
for all $t\in\real$ and let $h$ be a strictly
faithful locally finite lower semicontinuous weight on $A$ such that 
$h\comp\tau_t=\lambda^th$ for some fixed $\lambda>0$ and
for all $t\in\real$. We shall use the GNS triple $(\Hil,\pi,\eta)$ 
related to $h$ and identify $A$ with
$\pi(A)\subset\skrL(\Hil)$.\vs

\noindent Assume that a unitary element $W\in M(\skrK(\Hil)\otimes 
A)$ satisfies the following two conditions:\vs

$1.$ For any $a\in\Dom(\eta)$ and any $\varphi\in A^*$ we have:
\begin{equation}
\label{zalozenie}
\left.\begin{array}{c}
\varphi*a\in\Dom(\eta)\mbox{ and}\\
\Vs{5}\eta(\varphi*a)=(\id\otimes\varphi)W\eta(a),
\end{array}\right\}
\end{equation}

$2.$ For any $\varphi\in A^*$ such that $\varphi\comp\kappa\in A^*$ we have:
\begin{equation}
\label{zalozenie1}
(\id\otimes\varphi\comp\kappa)W=(\id\otimes\varphi)(W^*).
\end{equation}
Then $(A,\delta)$ is a weighted Hopf C$^*$-algebra and $h,R,\tau$ and 
$W$ coincide with the right Haar weight, the unitary
antipode, the scaling group and the Kac-Takesaki operator related to 
$(A,\delta)$.

\end{Thm}
\begin{pf}
We have to show that $h$ is right invariant and that the strong right 
invariance holds.
We shall follow the proof of Assertion 1 of Proposition \reef{pr3.6}. Let
\[
\eta(b)=\sum_{n\in\natu}b\Omega_n
\]
be an exact vector presentation of the GNS map $\eta$. Then
\[
h(b)=\sum_{n\in\natu}\IS{\Omega_n}{b\Omega_n}
\]
for any $b\in A_+$. Take $x,y,z\in\Hil$ and $a\in\Dom(\eta)$. Then, 
by \rf{zalozenie}, we have
\begin{equation}
\label{dorazne}
\begin{array}{r@{\;}c@{\;}l}
\ITS{x\otimes z}{W}{\eta(a)\otimes y}
&=&\ITS{x}{(\id\otimes\omega_{z,y})W}{\eta(a)}
=\IS{x}{\eta((\id\otimes\omega_{z,y})\delta(a))}\\&=&\Vs{7}
{\displaystyle 
\sum_{n\in\natu}\ITS{x}{(\id\otimes\omega_{z,y})\delta(a)}{\Omega_n}}=
{\displaystyle \sum_{n\in\natu}\ITS{x\otimes z}{\delta(a)}{\Omega_n\otimes y}}.
\end{array}
\end{equation}
Therefore
\begin{equation}
\label{KT21}
W(\eta(a)\otimes y)=\sum_{n\in\natu}\delta(a)(\Omega_n\otimes y).
\end{equation}
  Computing the norm of the both sides, we obtain:
\[
\begin{array}{r@{\;}c@{\;}l}
\omega_{y,y}(1) h(a^*a)&=&\norm{y}^2 h(a^*a)={\displaystyle 
\sum_{n\in\natu}}\,\norm{\delta(a)(\Omega_n\otimes y)}^2
={\displaystyle \sum_{n\in\natu}}\,\ITS{\Omega_n\otimes 
y}{\delta(a^*a)}{\Omega_n\otimes y}
\\ \Vs{7}
&=&
{\displaystyle \sum_{n\in\natu}}\,\ITS{\Omega_n}{\omega_{y,y}*(a^*a)}{\Omega_n}
=h(\omega_{y,y}*(a^*a)). \Vs{7}
\end{array}
\]
In this way we proved that $\psi(1) h(a^*a)=h(\psi*(a^*a))$ for all 
$\psi$ of the
form
$\omega_{y,y}$ ($y\in\Hil$) i.e. for all $\psi\in A_{*+}$. We have to 
show that this formula holds for all $\psi\in
A_{+}^*$.\vs

Let $p_n$ be the projection onto the closure of $A\Omega_n$. Then 
$p_n\in A'$ and $\sum_{n\in\natu
}p_n=1$. Replacing in \rf{dorazne}, $x$ by $p_nx$, we obtain:
\[
\ITS{x\otimes z}{(p_n\otimes 1)W}{\eta(a)\otimes y}=
\ITS{x\otimes z}{\delta(a)}{\Omega_n\otimes y}.
\]
Therefore
\[
(\omega_{x,\eta(a)}\otimes\id)\left[(p_n\otimes 
1)W\right]=(\omega_{x,\Omega_n}\otimes\id)\delta(a)\in A.
\]
\Vs{4}Let $\pi'$ be a representation of $A$ acting on a Hilbert space 
$\Hil'$. Then
\[
(\omega_{x,\eta(a)}\otimes\pi')\left[(p_n\otimes 
1)W\right]=(\omega_{x,\Omega_n}\otimes\pi')\delta(a)
\]
and for any $z',y'\in\Hil'$ we have
\[
\ITS{x\otimes z'}{(p_n\otimes 1)(\id\otimes\pi')W}{\eta(a)\otimes y'}=
\ITS{x\otimes z'}{(\id\otimes\pi')\delta(a)}{\Omega_n\otimes y'}.
\]
This way we showed that
\[
(p_n\otimes 1)(\id\otimes\pi')W(\eta(a)\otimes y')=
(\id\otimes\pi')\delta(a)(\Omega_n\otimes y')
\]
and summing over $n$ we obtain
\[
(\id\otimes\pi')W(\eta(a)\otimes y')=
\sum_n(\id\otimes\pi')\delta(a)(\Omega_n\otimes y')
\]
with the norm convergent sum on the right hand side. This relation 
holds for any $a\in\Dom(\eta)$. Computing the
square of the norm of both sides, we obtain:
\begin{equation}
\label{777}
\IS{\eta(a)}{\eta(a)}\IS{y'}{y'}=
\sum_n\ITS{\Omega_n\otimes 
y'}{(\id\otimes\pi')\delta(a^*a)}{\Omega_n\otimes y'}.
\end{equation}
For any $b\in A$ we set:
\begin{equation}
\label{stanpsi}
\psi(b)=\IS{y'}{\pi'(b)y'}.
\end{equation}
Then formula \rf{777} takes the form
\[
h(a^*a)\psi(1)=
\sum_n\ITS{\Omega_n}{(\id\otimes\psi)\delta(a^*a)}{\Omega_n}=h(\psi*(a^*a)).
\]
To end the proof of the right invariance we notice that any positive 
functional $\psi$ on $A$ is of the
form \rf{stanpsi} by the GNS construction.\vs

We pass to the strong right invariance. It follows immediately from 
Theorem \reef{SRIC} (Equivalence of Statements 1
and 2). Therefore  $(A,\delta)$ is a weighted Hopf C$^*$-algebra. 
Formula \rf{KT21} is identical with \rf{KT2}. It shows that $W$ is 
the Kac-Takesaki operator related to $(A,\delta)$.
\end{pf}\vs

The reader should notice that Assumption \rf{zalozenie} was used only 
with $\varphi$ of the form $\omega_{z,y}$
where $z,y\in\Hil$. On the other hand the set of pairs 
$(\varphi,a)\in A^*\times\Dom(\eta)$ satisfying \rf{zalozenie} is
closed with respect to the norm topology on $A^*$ and the graph 
topology on $\Dom(\eta)$. This observation leads to the
following

\begin{Rem}
\label{koniec3}
It is sufficient to verify condition \rf{zalozenie} for all $\varphi$ in a norm
dense subset of $A_*$ and for all $a$
in a core of $\eta$.
\end{Rem}

Taking into account the equivalence of Statements 2 and 3 of Theorem 
\reef{SRIC} , we obtain

\begin{Rem}
\label{koniec31}
It is sufficient to verify condition \rf{zalozenie1} for all $\varphi$ from a
weakly$^*$ dense
$\tau$-in\-va\-riant subset of $A^*$.
\end{Rem}

\section{Modular Structure for Haar Weight} \label{sec4}

In what follows $W$ always denotes the Kac-Takesaki operator for $(A,\delta)$,
and $A_*$ denotes the
set  of all $\sigma$-weakly continuous functionals on $A$.
In this section, we investigate the Connes' Radon-Nikodym cocycle of the left
invariant Haar weight $h^L:=h\comp R$ with respect to the right 
invariant Haar weight $h$. For this
purpose, we have to work on the von Neumann algebra $M=A''$ which is 
the weak closure of
the $C^*$-algebra $A$ acting on the GNS Hilbert space $\mathcal H$. 
By \rf{scaling} and
\rf{coproduct}, the scaling automorphisms $\tau_t$ and the 
comultiplication $\delta$ are weakly
continuous with respect to the weak operator topology, so they admit 
weakly continuous extensions
to $M$. Clearly
\[
\delta(a)=W(a\otimes 1)W^*, \quad \tau_t(a)=Q^{2it}aQ^{-2it}
\]
for any element $a\in M$. By Theorem \reef{th3.8} the
unitary antipode $R$ is normal (i.e. weakly continuous), so it admits 
a weakly continuous extension
to an involutive antiautomorphism acting on $M$. If necessary, we 
shall use the notations
$\delta^M$, $R^M$ and $\tau^M$ for them. The involutivity of the 
unitary antipode and the
commutativity of the unitary antipode and the scaling group follow 
immediately.\vs

The scaling group $\tau^M=\seq{\tau^M_t}{t\in\real}$ is no longer 
(pointwise) norm continuous so the theory presented in
Appendix \reef{AnaGene} can not be applied. Nevertheless one may 
speak about the analytic generator of
this group in the weak sense: an element $a\in M$ belongs to the 
domain of $\tau^M_{i/2}$ and
$b=\tau^M_{i/2}(a)$ if and only if for any $x,y\in\Hil$ there exists 
a continuous function $F_{x,y}$
defined on the strip $\set{z\in\compl}{0\leq \Im z\leq 1/2}$ 
holomorphic on the interior of this
strip, such that $F_{x,y}(t)=\IS{x}{\tau^M_t(a)y}$ for all $t\in\real$ and
$F_{x,y}(i/2)=\IS{x}{by}$. The antipode $\kappa^M$ on $M$ is defined by
\[
\kappa^M=\tau_{i/2}^M\comp R^M.
\]
In a sense $\tau^M_{i/2}$ is a closure of $\tau_{i/2}$. If
$a\in\Dom (\tau^M_{i/2})$ and $b=\tau^M_{i/2}(a)$ then there exists a net
$\seq{a_{\varepsilon}}{\varepsilon >0}$ of elements of $\Dom (\tau_{i/2})$
converging to $a$ such that $\tau_{i/2}(a_{\varepsilon})\rightarrow b$ when
$\varepsilon\rightarrow 0$ (both convergence in the sense of $\sigma$-weak
topology). Indeed
$a_{\varepsilon}=\frac 1{\sqrt{\pi\varepsilon}}\int 
e^{-t^2/\varepsilon}\tau_t(a)\;dt$
does the job. In the same sense $\kappa^M$ is the closure of $\kappa$. \vs

Let $\eta$ be the GNS-map associated with the right invariant weight
$h$ on the $C^*$-algebra $A$.  By the theory presented in Appendix 
\reef{GNSmap} the double commutant $\eta''$ is an extension of $\eta$ 
and $\Dom(\eta)$ is a core for $\eta''$. It is a GNS-map defined on 
the von Neumann algebra $M$.\vs

Using \rf{Qdef} one can easily show that $\Dom(\eta'')$ is $\tau^M$ - 
invariant  and
\begin{equation}
\label{4.1}
Q^{2it}\eta''(a)=\lambda^{-t/2}\eta''(\tau^M_t(a))
\end{equation}
for all $a\in\Dom(\eta'')$.\vs

\begin{Prop}
\label{odwr18}
Let $a\in\Dom(\eta'')\cap\Dom(\tau^M_{i/2})$ and 
$\eta''(a)\in\Dom(Q^{-1})$. Then
$\tau^M_{i/2}(a)\in\Dom(\eta'')$ and
\begin{equation}
\label{4.1anal}
\eta''\left(\tau^M_{i/2}(a)\right)=\lambda^{i/4}Q^{-1}\eta''(a).
\end{equation}
\end{Prop}

\begin{pf}
Let $b\in\Dom(\eta')$ and $x\in\Hil$. Applying $b$ to the both sides of
\rf{4.1} and using \rf{komutant}, we obtain:
\[
bQ^{2it}\eta''(a)=\lambda^{-t/2}\tau^M_t(a)\eta'(b)
\]
and
\[
\IS{x}{bQ^{2it}\eta''(a)}=\lambda^{-t/2}\IS{x}{\tau^M_t(a)\eta'(b)}
\]
for any $t\in\real$. Remembering that $a\in\Dom(\tau^M_{i/2})$ and
$\eta''(a)\in\Dom(Q^{-1})$, we see
that both sides of the above formula admit continuations to continuous
functions on the strip
$\set{t\in\compl}{0\leq\Im t\leq 1/2}$ holomorphic on the interior of this
strip. Inserting $t=i/2$ we obtain
\[
\IS{x}{bQ^{-1}\eta''(a)}=\lambda^{-i/4}\IS{x}{\tau^M_{i/2}(a)\eta'(b)}
\]
and $bQ^{-1}\eta''(a)=\lambda^{-i/4}\tau^M_{i/2}(a)\eta'(b)$. 
Therefore (cf. \rf{komutant}) $\tau^M_{i/2}(a)\in\Dom(\eta'')$ and 
\rf{4.1anal} follows.
\end{pf}\vs

Let $\tau_t'(b)=Q^{2it}bQ^{-2it}$ for $b\in A'$. Since the scaling group
$\seq{\tau_t}{t\in\real}$ on $A$ is implemented by $\seq{Q^{2it}}{t\in\real}$,
$\seq{\tau_t'}{t\in\real}$ is a one parameter group of automorphisms 
of $A'$. Let $b\in\Dom(\eta')$.
Then for any $a\in\Dom(\eta)$ we have
\[
\tau_t'(b)\eta(a)=\lambda^{t/2}Q^{2it}b\eta(\tau_{-t}(a))
=\lambda^{t/2}Q^{2it}\tau_{-t}(a)\eta'(b)=\lambda^{t/2}aQ^{2it}\eta'(b).
\]
Hence $\tau_t'(b)\in\Dom(\eta')$ and
\begin{equation}
\eta'(\tau_t'(b))=\lambda^{t/2}Q^{2it}\eta'(b). \label{taucom}
\end{equation}

\begin{Prop}
\label{pr4.1}
The set $\set{c\in\Dom(\eta')}{\eta'(c)\in\Dom(Q^{-1})\cap\Dom(Q)}$ 
is a core for $\eta'$.
\end{Prop}

\begin{pf} For $c\in A'$ we put
\[
{\mathcal R}_{\varepsilon}(c)=\frac
1{\sqrt{\pi\varepsilon}}\int_{\real}e^{-t^2/\varepsilon}\tau_t'(c)dt.
\]

Let $c\in\Dom(\eta')$.
The mapping $t\mapsto \tau_t'(c)$ is strongly continuous and the mapping
$t\mapsto\eta'(\tau_t'(c))=\lambda^{t/2}Q^{2it}\eta'(c)$ is norm 
continuous. Since $\eta'$ is closed,
we find that ${\mathcal R}_{\varepsilon}(c)\in\Dom(\eta')$ and
\[
\begin{array}{r@{\;}c@{\;}l}
\eta'({\mathcal R}_{\varepsilon}(c)) &=&{\displaystyle \frac
1{\sqrt{\pi\varepsilon}}\int_{\real}e^{-t^2/\varepsilon}\eta'(\tau_t'(c))dt}\\ 
&=&{\displaystyle
\frac 1{\sqrt{\pi\varepsilon}} 
\int_{\real}\lambda^{t/2}e^{-t^2/\varepsilon}Q^{2it}dt\;\eta'(c)}
\Vs{8}={\displaystyle e^{-\varepsilon(\log Q-i\log \lambda/4)^2}}\;\eta'(c).
\end{array}
\]

Using this formula one can easily verify that $\eta'\left({\mathcal
R}_{\varepsilon}(c)\right)\in\Dom(Q^{-1})\cap\Dom(Q)$. Clearly 
$\eta'({\mathcal R}_{\varepsilon}(c))$
converges in norm to $\eta'(c)$ as $\varepsilon$ tends to 0. One can 
also easily verify that
${\mathcal R}_{\varepsilon}(c)$ converges strongly to $c$. Therefore the set
$\set{c\in\Dom(\eta')}{\eta'(c)\in\Dom(Q^{-1})\cap\Dom(Q)}$ is a core 
for $\eta'$.
\end{pf}\vs

One can easily verify that the right hand sides of \rf{conpro} make 
sense for $a\in M$ and
$\varphi,\psi\in A_*$. In other words the convolution products 
$\varphi*a$ and $a*\psi$ (with
$\varphi,\psi\in A_*$) defined originally for $a\in A$ admit strongly 
continuous extension to $a\in
M$.\vs

Using \rf{Ftran} and remembering that $\Dom(\eta)$ is a core for $\eta''$, we
see that $\varphi*a\in\Dom(\eta'')$ and
\begin{equation}
W_{\varphi}\eta''(a)=\eta''(\varphi*a).
\label{Ftrans}
\end{equation}
for any $a\in\Dom(\eta'')$ and $\varphi\in A_*$.\vs

In what follows, we denote by the same letter $h$ the natural extension of the
right Haar weight
to $M$. By definition $h$ is the semifinite normal weight on $M$ related to
the GNS map
$\eta''$ via formula \rf{waga}. By the strict faithfulness of the Haar weight,
$h$ is faithful.
Remembering that $R$ is normal we conclude that $h^L:=h\comp R$ is 
also a faithful
semifinite normal weight on $M$. This is the left Haar weight. Let
$\seq{\sigma_t}{t\in\real}$ and
$\seq{\sigma_t^L}{t\in\real}$ be the modular automorphism groups of
$M$ associated with $h$ and $h^L$.

\begin{Lem}
\label{lm4.2}
The following four formulae hold on the von Neumann algebra $M$.
\begin{enumerate}
\item $R\comp\sigma_t\comp R=\sigma_{-t}^L$ for $t\in\real$.
\item $\delta\comp\sigma_t=(\sigma_t\otimes\tau_t)\comp\delta$ for $t\in\real$.
\item $\delta\comp\sigma_s^L=(\tau_{-s}\otimes\sigma_s^L)\comp\delta$ for
$s\in\real$.
\item $\sigma_t\comp\sigma_s^L=\sigma_s^L\comp\sigma_t$ for $t, s
\in\real.$
\end{enumerate}
\end{Lem}

\begin{pf} \Vs{2}

Ad 1. Follows immediately from Proposition \reef{acanon}.

Ad 2. Using \rf{pr3.64b} we compute:
\[
\begin{array}{r@{\;}c@{\;}l}
\delta(\sigma_t(a)) &=&\Ad_W\comp\Ad_{\Delta^{it}\otimes Q^{2it}}(a\otimes 1)\\
\Vs{5}&=&\Ad_{\Delta^{it}\otimes Q^{2it}}\comp\Ad_W(a\otimes 1)\\
\Vs{5} &=&(\sigma_t\otimes\tau_t)\comp\delta(a).
\end{array}
\]

Ad 3. We know (cf. \rf{qdelta}) that $Q$ and $\Delta$ strongly 
commute. Therefore scaling
automorphisms $\tau_t$ commute with modular automorphisms $\sigma_s$.
Furthermore $\tau_t$ commute with $R$. Using Assertion 1 we see that $\tau_t$
commute with $\sigma_s^L$. The rest is a matter of computation. Denoting by
$\sigma$ the flip automorphism acting on $A\otimes A$ and using twice the
first formula of Proposition \reef{pr2.10}, we obtain
\[
\begin{array}{r@{\;}c@{\;}l}
\delta\comp\sigma_t^L&=&\delta\comp R\comp\sigma_{-t}\comp R\\ 
&=&\sigma\comp(R\otimes
R)\comp\delta\comp\sigma_{-t}\comp R\Vs{5}\\
&=&\sigma\comp((R\comp\sigma_{-t})\otimes(R\comp\tau_{-t}))\comp\delta\comp 
R\Vs{5}\\
&=&\sigma\comp((R\comp\sigma_{-t}\comp R)\otimes(R\comp\tau_{-t}\comp R))
\comp\sigma\comp\delta\Vs{5}\\
&=&(\tau_{-t}\otimes\sigma_t^L)\comp\delta.\Vs{5}
\end{array}
\]
Ad 4. We already know that the scaling automorphisms commute with the 
modular automorphisms.
Therefore $(\tau_{-t}\otimes\sigma_t^L)\comp(\sigma_s\otimes\tau_s)=
(\sigma_s\otimes\tau_s)\comp(\tau_{-t}\otimes\sigma_t^L)$ and
\[
\begin{array}{r@{\;}c@{\;}l}
\delta\comp\sigma_t^L\comp\sigma_s 
&=&(\tau_{-t}\otimes\sigma_t^L)\comp\delta\comp\sigma_s
=(\tau_{-t}\otimes\sigma_t^L)\comp(\sigma_s\otimes\tau_s)\comp\delta\\&=&\Vs{5}
(\sigma_s\otimes\tau_s)\comp(\tau_{-t}\otimes\sigma_t^L)\comp\delta=
(\sigma_s\otimes\tau_s)\comp\delta\comp\sigma_t^L=
\delta\comp\sigma_s\comp\sigma_t^L.
\end{array}
\]
Now, using the injectivity of the comultiplication we obtain
$\sigma_t^L\comp\sigma_s=\sigma_s\comp\sigma_t^L$.
\end{pf}\vs

\begin{Cor}
If $h$ is a trace then $\kappa$ coincides with the unitary antipode $R$.
\end{Cor}

\begin{pf}
Let $t\in\real$. If $h$ is a trace then $\sigma_t=\id$ and by
Statement 2, $\id\otimes\tau_t$ coincides with the identity on
$\delta(A)$. It follows that $\id\otimes\tau_t$ coincides with the
identity on $(A\otimes 1)\delta(A)$. By the cancellation property
the latter set is linearly dense in $A\otimes A$. Hence
$\tau_t=\id$ for all $t\in\real$. Formula \rf{poldec} shows now
that $\kappa=R$.
\end{pf}\vs

We now consider the relative modular operator determined by $h$ and 
$h^L=h\comp R$.  Let $h'$ be the
weight on $M'$ associated with the commutant $\eta'$ of the 
GNS-mapping $\eta$ (see Appendices
\reef{B} and \reef{GNSmap}) and let $\Delta_{\rm rel}$ denote the 
spatial derivative
$dh^L/dh'$ in the sense of Connes \cite{Connes} (see Appendix 
\reef{TomTak}). Then $\Delta_{\rm rel}$ is a strictly positive
self-adjoint operator with a core
$\set{\eta''(a)}{a\in\Dom(\eta''),h^L(aa^*)<\infty}$, such that
\[
  h^L(aa^*)=\left\{
\begin{array}{ccc}
  \left \norm{\Delta_{\rm rel}^{1/2}\eta''(a)\right}^2 & \mbox{if} &
\eta''(a)\in\Dom(\Delta_{\rm rel}^{1/2}),\\
  +\infty & \quad & \mbox{otherwise}\Vs{6}
  \end{array}
  \right .
\]
for any $a\in\Dom(\eta'')$. It is known that the modular automorphism group
$\sigma^L$ is implemented by $\Delta_{\rm rel}$:
\[
\sigma_t^L(a)=\Delta_{\rm rel}^{it}a\Delta_{\rm rel}^{-it}
\]
for any $a\in M$ and $t\in\real$. Moreover the Connes' Radon-Nikodym 
cocycle $(Dh^L:Dh)_t$ is expressed
by the formula
\begin{equation}
\label{CRN}
(Dh^L:Dh)_t=\Delta_{\rm rel}^{it}\Delta^{-it}\in M.
\end{equation}

\begin{Prop}
\label{gammarho}
There exist a strictly positive self-adjoint operator $\gamma$ 
affiliated with the center of
the von Neumann algebra $M$ and a strictly positive self-adjoint 
operator $\rho$ affiliated with the
von Neumann algebra $M$ such that
\begin{equation}
\label{defgaro}
(Dh^L:Dh)_t=\gamma^{it^2/2}\rho^{it}
\end{equation}
for all $t\in\real$. Operators $\rho$ and $\gamma$ are uniquely 
determined by the above formula. They
satisfy the following commutation relations$:$
\begin{equation}
\Delta^{it}\rho\Delta^{-it}=\gamma^{t}\rho, \quad
\Delta^{it}\gamma\Delta^{-it}=\gamma,
\label{oprho}
\end{equation}
\begin{equation}
Q^{it}\rho Q^{-it}=\rho,\hspace{7mm} Q^{it}\gamma Q^{-it}=\gamma.
\label{opq}
\end{equation}
for all $t\in\real$.
\end{Prop}

\begin{pf}  For any $t,\ s\in\real$, we set
\begin{equation}
\label{obstruction}
\gamma(t,s):=\Delta^{it}\Delta_{\rm
rel}^{is}\Delta^{-it}\Delta_{\rm rel}^{-is}.
\end{equation}
Clearly 
$\left.\Ad_{\gamma(t,s)}\right|_M=\sigma_t\sigma_s^L\sigma_{-t}\sigma^L_{-s}$. 
By Assertion
4 of Lemma \reef{lm4.2}, the latter automorphism coincides with the identity on
$M$. Therefore $\gamma(t,s)$
commutes with all elements of $M$: $\gamma(t,s)\in M'$. On the other hand one
can easily verify that
\[
\gamma(t,s)=\Ad_{\Delta^{it}}(\Delta_{\rm rel}^{is}\Delta^{-is})
\left(\Delta_{\rm rel}^{is}\Delta^{-is}\right)^*
=\sigma_t\left((Dh^L:Dh)_s\Vs{3}\right)(Dh^L:Dh)_s^*.
\]
Therefore $\gamma(t,s)\in M$. This way we showed that the unitary 
$\gamma(t,s)$ belongs to the center
of $M$. It is known that central elements are invariant under the action of
modular automorphisms.
Therefore $\gamma(t,s)$ commutes with $\Delta^{it}$ and $\Delta_{\rm 
rel}^{is}$.\vs

Rewriting now \rf{obstruction} in the form $\gamma(t,s)\Delta_{\rm 
rel}^{is}=\Delta^{it}\Delta_{\rm
rel}^{is}\Delta^{-it}$ one can easily show that 
$\gamma(t,s+s')=\gamma(t,s)\gamma(t,s')$ for all
$t,s,s'\in\real$. Similarly using the formula
$\Delta^{-it}\gamma(t,s)=\Delta_{\rm rel}^{is}\Delta^{-it}\Delta_{\rm 
rel}^{-is}$ one can easily show
that $\gamma(t+t',s)=\gamma(t,s)\gamma(t',s)$ for all $t,t',s\in\real$.\vs

By the above results $\gamma(t,s)$ is of the form
\[
\gamma(t,s)=\gamma^{ist},
\]
where $\gamma$ is a strictly positive self-adjoint operator 
affiliated with the center of $M$.
Clearly $\gamma$ strongly commutes with $\Delta$ and $\Delta_{\rm 
rel}$. Therefore the second
formula of \rf{oprho} holds.\vs

Let $t,s\in\real$. Then \rf{obstruction} shows that 
$\Delta^{-it}\Delta_{\rm rel}^{is}=
\gamma^{-its}\Delta_{\rm rel}^{is}\Delta^{-it}$. Therefore
\[
\begin{array}{r@{\;}c@{\;}l}
\gamma^{-it^2/2}\Delta_{\rm rel}^{it}\Delta^{-it}
\gamma^{-is^2/2}\Delta_{\rm rel}^{is}\Delta^{-is}&=&
\gamma^{-it^2/2-its-is^2/2}\Delta_{\rm rel}^{it+is}\Delta^{-it-is}\\&=&
\gamma^{-i(t+s)^2/2}\Delta_{\rm rel}^{i(t+s)}\Delta^{-i(t+s)}.\Vs{6}
\end{array}
\]
It shows that $\seq{\gamma^{-it^2/2}\Delta_{\rm 
rel}^{it}\Delta^{-it}}{t\in\real}$ is a one
parameter group of unitaries. By \rf{CRN} these unitaries belong to 
$M$. Therefore there exists a
strictly positive self-adjoint operator $\rho$ affiliated with the 
von Neumann algebra $M$ such that
\[
\rho^{it}=\gamma^{-it^2/2}\Delta_{\rm rel}^{it}\Delta^{-it}
\]
for all $t\in\real$. Clearly \rf{defgaro} is equivalent to the above 
formula.\vs

Inserting $\frac{t}{n}$ instead of $t$ and taking the $n$-th power of 
the both sides we get
\[
\rho^{it}=
\gamma^{it^2/2n}\left(\Delta_{\rm rel}^{it/n}\Delta^{-it/n}\right)^n.
\]
Letting $n\rightarrow\infty$ we obtain
\begin{equation}
\label{Troter}
\rho^{it}=\mbox{s-}\hspace{-2mm}\lim_{n\rightarrow\infty}
\left(\Delta_{\rm rel}^{it/n}\Delta^{-it/n}\right)^n.
\end{equation}
This formula shows the uniqueness of $\rho$ (and $\gamma$).\vs

To prove \rf{opq} we notice that $\tau_t$ scales $h$ and $h^L$ by the 
same factor $\lambda^t$ (this
is because $\tau_t$ and $R$ commute). Therefore the Connes' Radon-Nikodym
cocycle \rf{CRN} is $\tau_t$ invariant. Using the uniqueness of $\rho$ and
$\gamma$ we obtain: $\tau_t(\rho)=\rho$,
$\tau_t(\gamma)=\gamma$ and \rf{opq} follows.\vs

By \rf{defgaro}, $\Delta_{\rm 
rel}^{it}=\gamma^{it^2/2}\rho^{it}\Delta^{it}$. Therefore for any
$t,s\in\real$ we have:
\[
\gamma^{i(t+s)^2/2}\rho^{i(t+s)}\Delta^{i(t+s)}=
\gamma^{it^2/2}\rho^{it}\Delta^{it}\gamma^{is^2/2}\rho^{is}\Delta^{is}.
\]
By simple computation this formula reduces to
\begin{equation}
\label{komrodelta}
\gamma^{its}\rho^{is}\Delta^{it}=\Delta^{it}\rho^{is}.
\end{equation}
Therefore 
$\left(\gamma^{t}\rho\right)^{is}=\Delta^{it}\rho^{is}\Delta^{-it}$ 
and the first formula
of \rf{oprho} follows.\vs
\end{pf}\vs

Formulae \rf{oprho} and \rf{opq} show that
\begin{equation}
\sigma_t(\rho)=\gamma^{t}\rho, \quad
\sigma_t(\gamma)=\gamma,
\label{oprho56}
\end{equation}
\begin{equation}
\tau_t(\rho)=\rho,\hspace{7mm} \tau_t(\gamma)=\gamma.
\label{opq56}
\end{equation}
for all $t\in\real$.\vs

\begin{Prop}\label{pr4.5} The operators $\rho$ and $\gamma$ satisfy 
$R(\rho^{it})=\rho^{-it}$ and
$R(\gamma^{it})=\gamma^{it}$ for $t\in\real$. In other words 
$R(\rho)=\rho^{-1}$ and
$R(\gamma)=\gamma$.
\end{Prop}

\begin{pf} Inserting $h_1=h^L$ in Assertion 2 of Proposition \reef{acanon} we
obtain
\[
(Dh^L:Dh)_t=R\left((Dh^L:Dh)_{-t}\right)=
R\left(\gamma^{it^2/2}\rho^{-it}\right)=R(\gamma)^{it^2/2}R(\rho)^{-it}.
\]
Using the uniqueness of $\rho$ and $\gamma$ we get: $R(\rho)=\rho^{-1}$ and
$R(\gamma)=\gamma$.
\end{pf}\vs

Formula \rf{komrodelta} shows that one parameter groups 
$\seq{\rho^{it}}{t\in\real}$,
$\seq{\Delta^{it}}{t\in\real}$ and $\seq{\gamma^{it}}{t\in\real}$ 
generate a unitary representation
of the Heisenberg group. Using the formula $\Delta_{\rm
rel}^{it}=\gamma^{it^2/2}\rho^{it}\Delta^{it}$ we see that 
$\seq{\Delta_{\rm rel}^{it}}{t\in\real}$
is an another one parameter group appearing in this representation. 
Its infinitesimal generator
$\log\Delta_{\rm rel}=\log\Delta+\log\rho$. Then \rf{Troter} is the 
Trotter product formula
corresponding to the relation $\log\rho=\log\Delta_{\rm 
rel}-\log\Delta$. (cf. \cite{Trot}).\vs

Inserting $s=t$ in \rf{komrodelta} we obtain
$\gamma^{it^2}\rho^{it}\Delta^{it}=\Delta^{it}\rho^{it}$. Combining 
\rf{defgaro} with \rf{CRN} we
get
\begin{equation}
\label{sekunda}
\Delta_{\rm rel}^{it}=\gamma^{it^2/2}\rho^{it}\Delta^{it}
      =\gamma^{-it^2/2}\Delta^{it}\rho^{it}
\end{equation}
for $t\in\real$. Now we put $t:=-i/2$ to obtain
\begin{equation}
\label{prima}
\Delta_{\rm rel}^{1/2}=\gamma^{-i/8}\overline{\rho^{1/2}\comp\Delta^{1/2}}
      =\gamma^{i/8}\overline{\Delta^{1/2}\comp\rho^{1/2}} \;.
\end{equation}
In particular 
$\Dom(\rho^{1/2}\comp\Delta^{1/2})\subset\Dom(\Delta_{\rm 
rel}^{1/2})$. In these formulae $\overline{X}$ denotes
the closure of an operator $X$ and $X\comp Y$ denotes the composition 
of operators $X$ and $Y$. By definition $\Dom(X\comp
Y)=\set{x\in\Dom(Y)}{Yx\in\Dom(X)}$.
\vs

The following result plays an important role in the paper.

\begin{Prop}\Vs{2}
\label{pr6.12}

1. If $a\in M$ and $\delta(a)=1\otimes a$ then $a\in\compl 1$.

2. If $a\in M$ and $\delta(a)=a\otimes 1$ then $a\in\compl 1$.

3. $M\cap\hA'=\compl 1$.
\end{Prop}

\begin{pf}\Vs{2}

Ad 1. By Proposition \reef{pr2.10}, $\delta$ commutes with the 
scaling group.  Therefore the set
\[
\set{a\in M}{\delta(a)=1\otimes a}
\]
  is $\tau$ invariant and we may assume that $a$ is an entire 
analytical element for $\tau$.
Using State\-ment~2 of Lemma \reef{lm4.2} we have:
$\delta(\sigma_t(a))=(\sigma_t\otimes\tau_t)(1\otimes 
a)=1\otimes\tau_t(a)$. It shows that $a$ is
an entire analytical element for the modular automorphism group. 
Hence (see Theorem
\reef{19.08.2002d}), for any $b\in\Dom(\eta)$ we have 
$ba\in\Dom(\eta'')$. Therefore for any
$x\in\Hil$ we have
\[
\begin{array}{r@{\;}c@{\;}l}
W(\eta''(ba)\otimes x)&=&{\displaystyle 
\sum_{n=1}^{\infty}}\delta(ba)(\Omega_n\otimes x)
={\displaystyle \sum_{n=1}^{\infty}}\delta(b)(\Omega_n\otimes ax)\\ 
&=&W(\eta(b)\otimes ax).\Vs{6}
\end{array}
\]
Hence $\eta''(ba)\otimes x=\eta(b)\otimes ax$ and for any 
$c\in\Dom(\eta')$ we have:
\[
ba\eta'(c)\otimes x=c\eta''(ba)\otimes x=c\eta(b)\otimes 
ax=b\eta'(c)\otimes ax.
\]
It shows that $a\otimes 1=1\otimes a$ and $a$ must be a scalar multiples of
the identity.\vs

Ad 2. Using the first Statement of Proposition \reef{pr2.10} we have:
\[
\delta(R^M(a))=\sigma\comp(R^M\otimes 
R^M)\delta(a)=\sigma\comp(R^M\otimes R^M)(a\otimes 1)=1\otimes R^M(a).
\]
Therefore, by Statement 1, $R^M(a)\in\compl 1$ and finally $a\in\compl 1$.\vs

Ad 3. We know that $W\in M(\hA\otimes A)$. Formula \rf{coproduct} 
shows that $\delta(a)=a\otimes 1$ for any $a\in
M\cap\hA'$ and Statement 3 follows immediately from Statement 2.
\end{pf}\vs

Let $LS$ (`LS' stays for left shifts) be the set of all automorphisms 
$\alpha$  of the von Neumann algebra $M$ such that
\begin{equation}
\label{LS}
\delta\comp\alpha=(\alpha\otimes\id)\comp\delta.
\end{equation}
We shall use the $u$-topology (cf. \cite[Section 2.23]{Stra}) on the 
set of automorphisms of $M$. By definition a net
$\alpha_n$ is $u$-convergent to $\alpha_\infty$ if for any $\varphi\in M_*$,
the net of functionals
$\varphi\comp\alpha_n$ converges in norm to $\varphi\comp\alpha_\infty$. One
can easily verify that $LS$ endowed with the $u$-topology is a  topological
group. In particular $\alpha_n^{-1}\rightarrow\alpha_\infty^{-1}$ if
$\alpha_n\rightarrow\alpha_\infty$.\vs

Let $t\in\real$. Combining Statement 2 of Proposition \reef{pr2.10} with
Statement 2 of Lemma \reef{lm4.2}, one can easily verify
that 
$\delta\comp\sigma_t\comp\tau_{-t}=(\sigma_t\comp\tau_{-t}\otimes\id)\comp\delta$. 
It shows that
\begin{equation}
\label{inLS}
\sigma_t\comp\tau_{-t}\in LS
\end{equation}
It is easy to see that the mapping $\real\ni t\mapsto\sigma_t\comp\tau_{-t}$
is $u$-continuous.
\begin{Prop}
\label{leftshifts}
Let $a\in A$. Then $\alpha(a)\in A$ for any $\alpha\in LS$. Moreover 
the mapping
\[
LS\ni\alpha\longmapsto\alpha(a)\in A
\]
is norm continuous $($we use u-topology on $LS)$.
\end{Prop}

\begin{pf} It is sufficient to prove this statement for $a$ running 
over a linearly dense subset of $A$.
Let $\alpha\in LS$. Using \rf{character} and \rf{LS}, one can easily 
verify that
\begin{equation}
\label{pr4.7.1}
(\id\otimes\delta)\left(\Vs{3}(\id\otimes\alpha)(W)W^*\right)=\left(\Vs{3}(\id\otimes\alpha)(W)W^*\right)\otimes 
1.
\end{equation}
Indeed
\[
\begin{array}{r@{\;}c@{\;}l}
(\id\otimes\delta)\left(\Vs{3}(\id\otimes\alpha)(W)W^*\right)&=&
(\id\otimes\delta\comp\alpha)W(\id\otimes\delta)W^*=
(\id\otimes\alpha\otimes\id)(\id\otimes\delta)W(\id\otimes\delta)W^*\\
\Vs{5}&=&\left(\Vs{4}(\id\otimes\alpha\otimes\id)W_{12}W_{13}\right)\left(W_{12}W_{13}\Vs{4}\right)^*
=\left(\Vs{3}(\id\otimes\alpha)(W)W^*\right)_{12}
\end{array}
\]
and \rf{pr4.7.1} follows. Using now Statement 2 of Proposition 
\reef{pr6.12} we conclude that $(\id\otimes\alpha)(W)W^*$ is of
the form $v\otimes 1$, where $v\in\skrL(\Hil)$. Hence
\[
(\id\otimes\alpha)(W)=(v\otimes 1)W.
\]
Clearly $v$ is unitary. Passing from $\Hil\otimes\Hil$ to 
$\Hil\otimes\Hil\otimes\Hil$ we obtain:
\[
\begin{array}{r@{\;}c@{\;}l}
(\id\otimes\alpha\otimes\id)(W_{12})&=&(v\otimes 1\otimes 1)W_{12},\\
(\id\otimes\id\otimes\alpha)(W_{13})&=&(v\otimes 1\otimes 1)W_{13}.\Vs{5}
\end{array}
\]
Eliminating $(v\otimes 1\otimes 1)$ we get:
\[
(\id\otimes\alpha\otimes\alpha)\left(W_{12}^*W_{13}\right)=W_{12}^*W_{13}
\]
and finally
\begin{equation}
\label{pr4.7.2}
(\id\otimes\id\otimes\alpha)\left(W_{12}^*W_{13}\right)=(\id\otimes\alpha^{-1}\otimes\id)\left(W_{12}^*W_{13}\right).
\end{equation}

Let $x,y\in\Hil$ and $\mu\in M_*$. Applying
$\omega_{x,y}\otimes\mu^*\otimes\id$ to the both sides of \rf{pr4.7.2} we
obtain:
\[
\alpha\left(W(W_\mu x,y)\Vs{3}\right)=W\left(W_{\mu\comp\alpha^{-1}}x,y\right).
\]
The reader should notice that the right hand side of the above formula belongs
to $A$ and depends continuously on $\alpha$.
This way our statement is proved for all $a$ of the form $a=W(W_\mu 
x,y)$. To end the proof we recall that $\hA$ act on $\Hil$
in a nondegenerate way. Therefore the set of vectors $\set{W_\mu 
x}{\mu\in\skrL(\Hil)_*,x\in\Hil}$ is dense in $\Hil$.
Proposition \reef{pr2.3} shows now that the set \[\set{W(W_\mu 
x,y)}{\mu\in\skrL(\Hil)_*;x,y\in\Hil}\] is linearly dense
in $A$.
\end{pf}\vs

We shall use Proposition \reef{leftshifts} with $\alpha$ replaced by 
\rf{inLS}. It shows that $\sigma_t\comp\tau_{-t}$ acts
within $A$ and that for any $a\in A$, $\sigma_t\comp\tau_{-t}(a)$ 
depends on $t\in\real$ in a norm continuous way. Combining
this result with the known properties of the scaling group, we obtain

\begin{Prop}
The modular automorphism group $\seq{\sigma_t}{t\in\real}$ of $M$ associated
with the Haar weight $h$ may be restricted to $A$. It defines the pointwise
norm continuous one parameter group of automorphisms of $A$.
\end{Prop}

In this way  we have shown Assertion 3 of Theorem \reef{th1.5}.  Restricting
the second formula of Lemma \reef{lm4.2} to $A$, we obtain Assertion 4 of
Theorem \reef{th1.5}.\vs

The end of this Section is devoted to the proof of the uniqueness of the
antipode and the Haar weight.

\begin{Thm}
\label{th7.3}
Let $(A,\delta)$ be a weighted Hopf
$C^*$-algebra with the Haar weight $h$ and antipode $\kappa$.  Assume that we
have another antipode $\kappa_1$ with corresponding Haar weight $h_1$
satisfying all the requirements of Definition \reef{def1.2}. Then
$\kappa_1=\kappa$ and $h_1=\mu h$ for some positive scalar $\mu$.
\end{Thm}

\begin{pf}
Let $(\Hil_1,\pi_1,\eta_1)$ be the GNS triple related to $h_1$.  Then 
we have another
Kac-Takesaki operator $W_1\in\skrL(\Hil_1\otimes\Hil_1)$ that gives 
rise to $(A,\delta)$. By a
Theorem in \cite{w5} we may assume that $\Hil_1=\Hil$ and 
$\pi_1=\pi$. Combining the formula
$(\id\otimes\delta)(W_1)=(W_1)_{12}(W_1)_{13}$ with
\rf{coproduct} we obtain: 
$W_{23}(W_1)_{12}W_{23}^*=(W_1)_{12}(W_1)_{13}$. It means that $W_1$ 
is
adapted to $W$ in the sense of
\cite[Defini\-tion~1.3]{w4}. Statement 4 of Theorem 1.6 of the same 
reference shows now that
$\kappa_1\subset\kappa$. Interchanging the roles of $W$ and
$W_1$ we obtain the converse inclusion. Therefore
\begin{equation}
\label{uniquekappa}
\kappa_1=\kappa.
\end{equation}

We shall use the relative conjugate linear Tomita-Takesaki operators
$\widetilde{S}=\widetilde{F}^*$. By definition the set
$\set{\eta(a)}{a\in\Dom(\eta)\cap\Dom(\eta_1)^*}$ is a core for
$\widetilde{S}$ and
\begin{equation}
\label{Stilde}
\widetilde{S}\eta(a)=\eta_1(a^*)
\end{equation}
for any  $a\in\Dom(\eta)\cap \Dom(\eta_1)^*$.
For any $a\in\Dom(\eta)\cap \Dom(\eta_1)^*$ and $\varphi\in A^*$ we have
$\varphi^* *a\in\Dom(\eta)$ and $\varphi*a^*\in\Dom(\eta_1)$. Therefore for any
$z\in\Dom(\widetilde{F})$ we have
\[
\begin{array}{r@{\;}c@{\;}l}
\varphi(W_1(z,\eta_1(a^*)))&=&(z|(W_1)_{\varphi}\eta_1(a^*)) 
=(z|\eta_1(\varphi*a^*))\\ \Vs{6}
&=&(z|\widetilde{S}\eta(\varphi^* 
*a))=\overline{(\widetilde{F}z|W_{\varphi^*}\eta(a))}\\ \Vs{6}
  &=&\overline{\varphi^*(W(\widetilde{F}z,\eta(a)))} = 
\varphi(W(\widetilde{F}z,\eta(a))^*).
\end{array}
\]
Hence we have $W_1(z,\eta_1(a^*))=W(\widetilde{F}z,\eta(a))^*$. Remembering
that $\set{\eta(a)}{a\in\Dom(\eta)\cap
\Dom(\eta_1)^*}$ is a core for $\widetilde{S}$, we conclude that
$W_1(z,\widetilde{S}x)=W(\widetilde{F}z,x)^*$ for any 
$x\in\Dom(\widetilde{S})$.
Setting $y=\widetilde{F}z$ we get
\begin{equation}
\label{lm7.4}
W_1(\widetilde{F}^{-1}y,\widetilde{S}x)=W(y,x)^*.
\end{equation}
This formula holds for any $x\in\Dom(\widetilde{S})$ and
$y\in\Dom(\widetilde{F}^{-1})$.\vs

Let $\widetilde{S}=\widetilde{J}\widetilde{\Delta}^{1/2}$ be the 
polar decomposition of $\widetilde{S}$.  Then
$\widetilde{F}=\widetilde{\Delta}^{1/2}\widetilde{J}^{-1}$. For any 
$x\in\Dom(\widetilde{\Delta})$ and 
$y\in\Dom(\widetilde{\Delta}^{-1})$, we have
$\widetilde{S}x\in\Dom(\widetilde{F})$ and 
$\widetilde{F}^{-1}y\in\Dom(\widetilde{S}^{-1})$. Using \rf{lm7.4} 
and \rf{uniquekappa} and Proposition \reef{pr2.6}
we get
\[
\begin{array}{r@{\;}c@{\;}l}
\kappa^2(W(x,y))&=&\kappa(W(y,x)^*) 
=\kappa(W_1(\widetilde{F}^{-1}y,\widetilde{S}x))=
\kappa_1(W_1(\widetilde{F}^{-1}y,\widetilde{S}x))\\ \Vs{6}
&=&W_1(\widetilde{S}x,\widetilde{F}^{-1}y)^*
=W(\widetilde{F}\widetilde{S}x,\widetilde{S}^{-1}\widetilde{F}^{-1}y)
=W(\widetilde{\Delta}x,\widetilde{\Delta}^{-1}y).
\end{array}
\]
Remembering that $\kappa^2=\tau_i$ and using Theorem \reef{drugie} , we obtain
\[
\tau_t(W(x,y))=W(\widetilde{\Delta}^{it}x,\widetilde{\Delta}^{it}y).
\]
for all $t\in\real$. On the other hand, 
$\tau_t(W(x,y))=W(\Delta^{it}x,\Delta^{it}y)$ by Proposition 
\reef{pr2.8}.
Therefore 
$W(\widetilde{\Delta}^{it}x,\widetilde{\Delta}^{it}y)=W(\Delta^{it}x,\Delta^{it}y)$ 
and
  $\widetilde{\Delta}^{-it}\Delta^{it}$ commutes with $W_{\varphi}$ 
for all $\varphi\in
A^*$.  Hence it belongs to the commutant of $\hA$.  On the other hand 
the Radon-Nikodym cocycle
$(Dh_1:Dh)_{-t}=\widetilde{\Delta}^{-it}\Delta^{it}\in M$. Therefore 
by Proposition \reef{pr6.12},
it must be a scalar multiple of the identity and there exists a 
positive scalar $\mu$ such that
$\widetilde{\Delta}^{it}=\mu^{it}\Delta^{it}$ for all $t\in\real$.  Hence
$\widetilde{\Delta}=\mu\Delta$. In particular
$\Dom(\widetilde{S})=\Dom(\widetilde{\Delta}^{1/2})=\Dom(\Delta^{1/2})$.\vs

The reader should notice that $\Dom(\eta)^*$ is a right ideal and
$\Dom(\eta_1)$ is a left ideal in $A$. Therefore
$\Dom(\eta)^*\Dom(\eta_1)\subset\Dom(\eta)^*\cap\Dom(\eta_1)$. Using the norm
density of $\Dom(\eta)^*$ one can easily show that $\Dom(\eta)^*\Dom(\eta_1)$
is a core for $\eta_1$. So is $\Dom(\eta)^*\cap\Dom(\eta_1)$. Let
$a\in\Dom(\eta)^*\cap\Dom(\eta_1)$. Then (cf. \rf{Stilde}) 
$\eta(a^*)\in\Dom(\widetilde{S})=\Dom(\Delta^{1/2})$ and using Theorem
\reef{19.08.2002c} we see that $a\in\Dom(\eta)$. Moreover
\begin{equation}
\label{dorazne1}
\eta_1(a)=\widetilde{J}\widetilde{\Delta}^{1/2}\eta(a^*)
=\sqrt{\mu}\,\widetilde{J}\Delta^{1/2}\eta(a^*) 
=\sqrt{\mu}\,{\widetilde J}J\eta(a).
\end{equation}
Remembering that $\Dom(\eta)^*\cap\Dom(\eta_1)$ is a core for $\eta_1$ we see
that $\Dom(\eta_1)\subset\Dom(\eta)$ and that
the above formula holds for all $a\in\Dom(\eta_1)$. Interchanging the roles of
$\eta$ and $\eta_1$ we obtain the converse inclusion. Hence
$\Dom(\eta_1)=\Dom(\eta)$. Formula \rf{dorazne1} shows now that
$h_1(a^*a)=\mu h(a^*a)$ for any $a\in A$.
\end{pf}

\section{Modular Structure for the Dual} \label{sek5}

We now go into the discussion on the density theorem. Combining 
(\reef{Q}) and (\reef{opq}), we know
that $J\rho^{-1}J$ and $Q$ are strongly commuting strictly positive 
self-adjoint operators.
Therefore the operator
\begin{equation}
\hDelta:=\overline{J\rho^{-1}J\comp Q^2} \label{hDelta}
\end{equation}
is also strictly positive and self-adjoint. This section is devoted 
to the proof of the following
theorem. It contains the inverse of the square root of (\reef{hDelta}).

\begin{Thm}\label{th5.1} Let
\begin{equation}
\Dom_0=\set{b\in\Dom(\eta'')\cap\Dom(\kappa^M)}{ \kappa^M(b)^*\in\Dom(\eta'')}.
\label{d0}
\end{equation}
Then the set $\set{\eta''(b)}{b\in\Dom_0}$ is a core for $\hDelta^{-1/2}$ and
\begin{equation}
\norm{\hDelta^{-1/2}\eta''(b)}=\norm{\eta''(\kappa^M(b)^*)} \label{dmconj0}
\end{equation}
for any element $b\in \Dom_0$.
\end{Thm}

We shall use the notation introduced in the previous section.

\begin{Lem}
\label{lm5.2}
Let $a$ be an element in $\Dom(\eta'')$. If $\eta''(a)\in {\mathcal
D}(\Delta_{\rm rel}^{1/2})$, then $R(a)\in \Dom(\eta'' )$ and
\begin{equation}
\norm{\eta'' (R(a))}=\norm{\Delta_{\rm rel}^{1/2}\eta'' (a)}. \label{spade0}
\end{equation}
\end{Lem}

\begin{pf} If $a$ is an element in $\Dom(\eta'')$ satisfying 
$\eta''(a)\in {\mathcal
D}(\Delta_{\rm rel}^{1/2})$, then $h(R(a)^*R(a))=h^L(aa^*)=\norm{\Delta_{\rm
rel}^{1/2}\eta''(a)}^2<\infty$. Hence we obtain $R(a)\in 
\Dom(\eta'')$ and $\norm{\eta''
(R(a))}=\norm{\Delta_{\rm rel}^{1/2}\eta'' (a)}$.
\end{pf}

\begin{Lem}\label{lm5.3} Let $a$ be an element in $\Dom(\eta'')$ with 
$\eta''(a)\in {\mathcal
D}(\Delta^{1/2})$.  If
$\eta''(a)\in \Dom(J\rho^{1/2}J)$, then $R(a)^*\in \Dom(\eta'')$ and
\begin{equation}
\norm{\eta''(R(a)^*)}=\norm{J\rho^{1/2}J\eta''(a)}. \label{spade1}
\end{equation}
\end{Lem}

\begin{pf} Since $\eta''(a)\in \Dom(\Delta^{1/2})$, we have $a^*\in 
\Dom(\eta'')$
(cf. Theorem \reef{19.08.2002c}) and
\[
\eta''(a^*)=J\Delta^{1/2}\eta''(a)=\Delta^{-1/2}J\eta''(a);
\]
hence, $\eta''(a^*)\in \Dom(\Delta^{1/2})$ and 
$\Delta^{1/2}\eta''(a^*)=J\eta''(a)$. Assume
that $\eta''(a)\in \Dom(J\rho^{1/2}J)$.  Then $J\eta''(a)\in 
\Dom(\rho^{1/2})$ and
$\eta''(a^*)\in \Dom(\rho^{1/2}\Delta^{1/2})\subset\Dom(\Delta_{\rm 
rel}^{1/2})$. Now
we apply Lemma \reef{lm5.2} with $a$ replaced by $a^*$ to obtain 
$R(a)^*=R(a^*)\in {\mathcal
D}(\eta'')$ and
\begin{eqnarray*}
\norm{\eta''(R(a)^*)}&=&\norm{\eta'' (R(a^*))}=\norm{\Delta_{\rm 
rel}^{1/2}\eta''(a^*)}\\
&=&\norm{\rho^{1/2}\Delta^{1/2}\eta''(a^*)}=\norm{\rho^{1/2}J\eta''(a)} 
=\norm{J\rho^{1/2}J\eta''(a)}.
\end{eqnarray*}
\end{pf}\vs

Next we have an improved version of the above statement, which plays 
the fundamental role in the
subsequent argument.

\begin{Prop}
\label{pr5.4}
Let $a$ be an element in $\Dom(\eta'')$.  If $\eta''(a)\in {\mathcal
D}(J\rho^{1/2}J)$, then $R(a)^*\in \Dom(\eta'')$ and
\begin{equation}
\norm{\eta''(R(a)^*)}=\norm{J\rho^{1/2}J\eta''(a)}. \label{spade2}
\end{equation}
\end{Prop}

\begin{pf} For any $b\in \Dom(\eta'')$, $b^*a\in \Dom(\eta'')$ and 
$\eta''(b^*a)\in
\Dom(\Delta^{1/2})$. Moreover due to the assumption $\eta''(a)\in 
\Dom(J\rho^{1/2}J)$
and the operator inclusion $J\rho^{1/2}Jb^*\supset b^*J\rho^{1/2}J$, we have
$\eta''(b^*a)=b^*\eta''(a)\in \Dom(J\rho^{1/2}J)$. Therefore the 
element $b^*a$ satisfies the
conditions of the previous Lemma and hence we have 
$R(b)R(a)^*=R(b^*a)^*\in \Dom(\eta'')$ and
$\norm{\eta''(R(b^*a)^*)}=\norm{J\rho^{1/2}Jb^*\eta''(a)}$. Therefore we have
\[
\norm{\eta''(R(b)R(a)^*)}=\norm{b^*J\rho^{1/2}J\eta''(a)}\leq
\norm{b}\; \norm{J\rho^{1/2}J\eta''(a)}.
\]
Now we shall use an exact vector presentation (see Definition 
\reef{defexact} and Theorem \reef{Bszesc}):
\[
\eta(a)=\sum_{n\in\natu} a\Omega_n
\]
of the GNS mapping $\eta$. Then $\eta''(R(b^*a)^*)=\sum_{n\in\natu 
}R(b)R(a^*)\Omega_n$, so for any natural $k$ and any 
$b\in\Dom(\eta'')$ with $\norm{b}\leq 1$ we have
\[
\sum_{n=1}^k\norm{R(b)R(a^*)\Omega_n}^2\leq\norm{J\rho^{1/2}J\eta''(a)}^2.
\]
The set $\set{b\in A''}{b\in\Dom(\eta''), \norm{b}\leq 1}$ is weakly 
dense in the unit ball of the von Neumann algebra $A''$. It follows 
that for any natural $k$,
\[
\sum_{n=1}^k\norm{R(a^*)\Omega_n}^2\leq\norm{J\rho^{1/2}J\eta''(a)}^2.
\]
and the series
$
\sum_{n\in\natu}\norm{R(a^*)\Omega_n}^2
$
is convergent. Remembering that our vector presentation is exact we 
see that $R(a^*)\in\Dom(\eta'')$. Now by taking elements $b$ in the 
unit ball of
$\Dom(\eta'')$, we obtain
\begin{eqnarray*}
\norm{\eta''(R(a^*))} &=&\sup_{\norm{b}\leq 1}
\norm{R(b)\eta''(R(a^*))}\\ &=&\sup_{\norm{b}\leq 1}
\norm{b^*J\rho^{1/2}J\eta''(a)} =\norm{J\rho^{1/2}J\eta''(a)}
\end{eqnarray*}
\end{pf}\vs

Let
\begin{equation}
\Dom_1=\set{a\in\Dom(\eta'')\cap\Dom(\kappa^M)}{
\begin{array}{c}
\eta''(a)\in\Dom(Q^{-1})\cap\Dom(J\rho^{1/2}J) \\
\mbox{and }Q^{-1}\eta''(a)\in\Dom(J\rho^{1/2}J)\Vs{4}
\end{array}
}.
\label{d1}
\end{equation}

\begin{Prop}\label{pr5.5}
$\Dom_1\subset\Dom_0$ and
\begin{equation}
\norm{\eta''(\kappa^M(a)^*)}=\norm{J\rho^{1/2}JQ^{-1}\eta''(a)} \quad 
\mbox{for} \quad a\in\Dom_1.
\label{dmconj}
\end{equation}
\end{Prop}

\begin{pf} Suppose that $a\in\Dom_1$. We know (cf. Proposition 
\reef{odwr18}) that
$\tau_{i/2}^M(a)\in\Dom(\eta'')$ and \linebreak
$\eta''(\tau_{i/2}^M(a))=\lambda^{i/4}Q^{-1}\eta''(a)$. By our assumption,
$Q^{-1}\eta''(a)\in\Dom(J\rho^{1/2}J)$. Applying Proposition 
\reef{pr5.4} to the element
$\tau_{i/2}^M(a)$ instead of $a$, we see that
$\kappa^M(a)^*=R(\tau_{i/2}^M(a))^*\in\Dom(\eta'')$. In
this case equality (\reef{spade2}) converts into (\reef{dmconj}).
\end{pf}\vs

Up to this point, it is not clear whether the sets $\Dom_0$ and 
$\Dom_1$ are not trivial.  To show that they
are large enough, we have to tame unbounded operators that appear in 
our theory. In order to tame
$\rho$ and $\gamma$, for any $\varepsilon>0$ we introduce two 
elements in $A''$ as follows:
\begin{equation}
\label{opGamma'}
\begin{array}{r@{\;}c@{\;}l}
\Gamma_{\varepsilon}&:=&\exp\left\{-\varepsilon\left((\log\rho )^2
     +\frac 14 (\log\gamma )^2\right)\right\},\\
\Gamma_{\varepsilon}'&:=&\exp\left\{-\varepsilon\left((\log\rho )^2
     -i(\log\rho )(\log\gamma )\right)\right\}.\Vs{5}
\end{array}
\end{equation}
Notice that $0\leq\Gamma_{\varepsilon}\leq 1$ and
$\Gamma_{\varepsilon}'$ is bounded.  Moreover, $\Gamma_{\varepsilon}$ 
and $\Gamma_{\varepsilon}'$
converges strongly to 1 as $\varepsilon$ tends to 0.

\begin{Lem}\label{lm5.6}
Let $a$ be an element in $\Dom(\eta'')$. Then
$a\Gamma_{\varepsilon}\in \Dom(\eta'')$ and
$\eta''(a\Gamma_{\varepsilon})=J\Gamma_{\varepsilon}'J\eta''(a)$.
\end{Lem}

\begin{pf}
Using (\reef{oprho}),  we can easily check that
\[
\sigma^h_t(\Gamma_{\varepsilon})=\Delta^{it}\Gamma_{\varepsilon}\Delta^{-it}
=\exp\left\{-\varepsilon\left((\log\gamma^t\rho)^2
+\frac 14(\log\gamma)^2\right)\right\}.
\]
The right hand side admits the holomorphic continuation for $t$ 
running over all $\compl$. It means that 
$\Gamma_{\varepsilon}\in\Dom(\sigma^h_{i/2})$ and setting in the 
above formula $t=i/2$ we obtain
\[
\sigma_{i/2}(\Gamma_{\varepsilon})=
\exp\left\{-\varepsilon\left(\left(\log\rho+\frac i2\log\gamma\right)^2 +\frac
14(\log\gamma)^2\right)\right\} = \left(\Vs{4}\Gamma_{\varepsilon}'\right)^*.
\]

Now our statement follows immediately from Theorem \reef{19.08.2002d}.
\end{pf}\vs

In order to tame the effects coming from the unboundedness of $\log 
Q$, for any  $\varepsilon>0$ and
any $a\in M$ we set
\begin{equation} {\mathcal R}_{\varepsilon}(a) 
=\frac{1}{\sqrt{\varepsilon\pi}}\int_\real
\tau_t^M(a\Gamma_{\varepsilon})\lambda^{-t/2}e^{-t^2/\varepsilon}dt.
\label{regular1}
\end{equation}

\begin{Rem}\label{rmk5.7} By \rf{opq}, $\rho$ and $\gamma$ are
$\seq{\tau_t}{t\in\real}$ invariant.  So is $\Gamma_{\varepsilon}$. Therefore
$\Gamma_{\varepsilon}\in\Dom(\tau_{i/2})$ and
$\tau_{i/2}(\Gamma_{\varepsilon})=\Gamma_{\varepsilon}$.  Applying 
$R$ to the both sides, we conclude
that $\Gamma_{\varepsilon}\in\Dom(\kappa^M)$ and
$\kappa^M(\Gamma_{\varepsilon})=R(\Gamma_{\varepsilon})=\Gamma_{\varepsilon}$. 
The latter equality
follows from Proposition $\reef{pr4.5}$. Since $\Dom(\kappa^M)$ is a 
subalgebra and that $\kappa^M$ is
antimultiplicative, we see that $a\Gamma_{\varepsilon}\in\Dom(\kappa^M)$ and
$\kappa^M(a\Gamma_{\varepsilon})^* =\kappa^M(a)^*\Gamma_{\varepsilon}$ for any
$a\in\Dom(\kappa^M)$. Using this result together with the closedness 
of $\kappa^M$, one can easily
show that ${\mathcal R}_{\varepsilon}(a)\in\Dom(\kappa^M)$ and
\begin{equation}
{\mathcal R}_{\varepsilon}(\kappa^M(a)^*)=\kappa^M({\mathcal 
R}_{\varepsilon}(a))^* \label{regular2}
\end{equation}
for any $a\in\Dom(\kappa^M)$.
\end{Rem}

\begin{Lem}
\label{lm5.8}\Vs{4}

$1.$ For any $a\in M$, ${\mathcal R}_{\varepsilon}(a)$ converges 
strongly to $a$ when $\varepsilon
\longrightarrow 0$.\vs

$2.$ Let $a\in\Dom(\eta'')$ and $\varepsilon>0$.  Then ${\mathcal 
R}_{\varepsilon}(a)\in\Dom_1$
and
\begin{equation}
\eta''({\mathcal R}_{\varepsilon}(a)) =\exp\left\{-\varepsilon(\log
Q)^2\right\}J\Gamma_{\varepsilon}'J\eta''(a). \label{regular3}
\end{equation}
In particular $\eta''({\mathcal R}_{\varepsilon}(a))$  converges in 
norm  to $\eta''(a)$ when
$\varepsilon\longrightarrow 0$.\vs

$3.$ The set $\Dom_1$ is a core for $\eta''$.
\end{Lem}

\begin{pf}
Ad 1. One can easily verify that
\[
\delta_\varepsilon(t)=\frac{1}{\sqrt{\varepsilon\pi}} 
\lambda^{-t/2}e^{-t^2/\varepsilon}
\]
is a $\delta$-like sequence when $\varepsilon\rightarrow 0$. More
precisely: $\delta_\varepsilon(t)>0$ for all $t\in\real$ and
$\varepsilon>0$, $\int_\real \delta_\varepsilon(t)dt=1$ and for
any neighborhood ${\mathcal O}$ of $0$ in $\real$,
$\int_{\real\backslash{\mathcal
O}}\delta_\varepsilon(t)dt\longrightarrow 0$ when
$\varepsilon\rightarrow 0$. Therefore for any continuous bounded
function $\varphi$ on $\real$ we have
\begin{equation}
\label{deltalike}
\lim_{\varepsilon\rightarrow 
0}\int_\real\varphi(t)\delta_\varepsilon(t)dt=\varphi(0).
\end{equation}

Let $x\in\Hil$.  Then
\[
\begin{array}{r@{\;}c@{\;}l}
\norm{{\mathcal R}_{\varepsilon}(a)x-ax} &=&\left\norm{ 
{\displaystyle \int_\real
\left(\tau_t^M(a)\Gamma_{\varepsilon}x-ax\Vs{4}\right) 
\delta_\varepsilon(t)dt} \right}
\\  \Vs{8}&\leq&{\displaystyle \int_\real
\left\norm{\tau_t^M(a)\left(\Gamma_{\varepsilon}x-x\right)\Vs{3}\right} 
\delta_\varepsilon(t)dt +
  \int_\real
\left\norm{\tau_t^M(a)x-ax\Vs{3}\right} \delta_\varepsilon(t)dt}
\\  \Vs{7}&\leq&{\displaystyle
\norm{a}\left\norm{\Gamma_{\varepsilon}x-x\right} +
  \int_\real
\left\norm{\tau_t^M(a)x-ax\Vs{3}\right} \delta_\varepsilon(t)dt}
\end{array}
\]
The first term on the right hand side converges to $0$ as 
$\varepsilon$ goes to $0$. Inserting in
\rf{deltalike}, $\varphi(t)=\left\norm{\tau_t^M(a)x-ax\Vs{3}\right}$ 
we see that also the second term
converges to $0$. It shows that ${\mathcal R}_{\varepsilon}(a)x$ 
converges in norm to $ax$.

\vs

Ad 2. Assume that $a\in\Dom(\eta'')$. We already know by Lemma 
\reef{lm5.6} that
$a\Gamma_{\varepsilon}\in\Dom(\eta'')$ and
$\eta''(a\Gamma_{\varepsilon})=J\Gamma_{\varepsilon}'J\eta''(a)$. By 
(\reef{4.1}),
$\tau_t^M(a\Gamma_{\varepsilon})\in\Dom(\eta'')$ and
$\eta''(\tau_t^M(a\Gamma_{\varepsilon})) 
=\lambda^{t/2}Q^{2it}J\Gamma_{\varepsilon}'J\eta''(a)$.
Therefore
\begin{eqnarray*}
\frac 1{\sqrt{\varepsilon\pi}}\int_\real
\eta''(\tau_t^M(a\Gamma_{\varepsilon}))\lambda^{-t/2}e^{-t^2/\varepsilon}dt
&=&\frac{1}{\sqrt{\varepsilon\pi}}\int_\real
Q^{2it}e^{-t^2/\varepsilon}J\Gamma_{\varepsilon}'J\eta''(a)dt \\ 
&=&e^{-\varepsilon(\log
Q)^2}J\Gamma_{\varepsilon}'J\eta''(a).
\end{eqnarray*} Remembering that $\eta''$ is closed, we conclude that
\begin{equation}
{\mathcal R}_{\varepsilon}(a) =\frac 1{\sqrt{\varepsilon\pi}}\int_\real
\tau_t^M(a\Gamma_{\varepsilon})\lambda^{-t/2} 
e^{-t^2/\varepsilon}dt\in\Dom(\eta'') \label{regular}
\end{equation}
and $\eta''({\mathcal R}_{\varepsilon}(a)) =e^{-\varepsilon(\log
Q)^2}J\Gamma_{\varepsilon}'J\eta''(a)$. Due to the dumping factor 
$e^{-\varepsilon(\log
Q)^2}J\Gamma_{\varepsilon}'J$, the vector $\eta''({\mathcal 
R}_{\varepsilon}(a))$
belongs to $\Dom(Q^{-1})$ and the vectors $\eta''({\mathcal 
R}_{\varepsilon}(a))$ and
$Q^{-1}\eta''({\mathcal R}_{\varepsilon}(a))$ belong to $\Dom(J\rho^{1/2}J)$.
\vs

One can easily verify that ${\mathcal R}_{\varepsilon}(a)$ is an entire
analytic element for the scaling group. Therefore
${\mathcal R}_{\varepsilon}(a)\in\Dom(\tau_{i/2}^M)=\Dom(\kappa^M)$. It
means that ${\mathcal R}_{\varepsilon}(a)\in\Dom_1$.\vs

Ad 3. Follows immediately from the first two assertions.
\end{pf}\vs

Now we recall a fact concerning a core for positive self-adjoint operators

\begin{Lem}\label{lm5.9} Let $H$ be a strictly positive self-adjoint 
operator on $\Hil$. Let $K$ be a
bounded operator such that $K\Hil$ is dense in $\Hil$, 
$K\Hil\subset\Dom(H)$, $HK\in{\mathcal
L}(\Hil)$ and $HKx=KHx$ for
$x\in\Dom(H)$. Then, for any dense linear subset $\Dom$ of $\Hil$, 
$K\Dom$ is a core for $H$.
\end{Lem}

\begin{pf} It suffices to show that the set $(H+1)K\Dom$ is dense in 
$\Hil$. Suppose that $y$ is
orthogonal to the above set. Then $(y|(H+1)Kx)=0$ holds for any 
vector $x$ in $\Dom$. Since $HK\in{\mathcal L}(\Hil)$ , it holds for 
$x$ in
$\Hil$. Therefore $(y|K(H+1)x)=(y|(H+1)Kx)=0$ holds for any 
$x\in\Dom(H)$. Since $(H+1)\Dom(H)=\Hil$
and $K$ has a dense range, it follows that $y=0$.
\end{pf}\vs

Now we shall use the operator $\hDelta$ introduced by (\reef{hDelta}). Clearly,
$\hDelta^{-1/2}=\overline{J\rho^{1/2}J\comp Q^{-1}}$.

\begin{Prop}
\label{pr5.10}
The set $\set{\eta''(b)}{b\in\Dom_1}$ is a core for $\hDelta^{-1/2}$.
\end{Prop}

\begin{pf} At first we notice that $\set{\eta''(b)}{b\in\Dom_1}
\subset\Dom(\hDelta^{-1/2})$ by (\reef{d1}). We use the previous lemma with
$H=\hDelta^{-1/2}=\overline{J\rho^{1/2}J\comp Q^{-1}}$,
$K=J\Gamma_1'Je^{-(\log Q)^2}$ and 
$\Dom=\set{\eta''(a)}{a\in\Dom(\eta'')}$. It is clear that $K$ is
bounded and the range is dense in $\Hil$. Let $x\in\Dom(H)$.  There 
exists a sequence
$\seq{x_n}{n\in\natu }$ in
$\Dom(J\rho^{1/2}J\comp Q^{-1})$ such that $x_n\to x$ and $Hx_n\to 
Hx$ in norm. For $x_n$ we see that
\begin{eqnarray*} KHx_n&=&J\Gamma'_1Je^{-(\log 
Q)^2}J\rho^{1/2}JQ^{-1}x_n =J\rho^{1/2}\Gamma'_1Je^{-(\log
Q)^2}Q^{-1}x_n\\ &=&J\rho^{1/2}\Gamma'_1JQ^{-1}e^{-(\log Q)^2}x_n 
=J\rho^{1/2}JQ^{-1}J\Gamma'_1Je^{-(\log
Q)^2}x_n =HKx_n.
\end{eqnarray*} Since $Kx_n\to Kx$ and $HKx_n=KHx_n\to KHx$, the 
closedness of $H$ implies that
$Kx\in\Dom(H)$ and $HKx=KHx$. Hence $K\Hil$ is contained in the 
domain of $H$ and the operator $HK$
is bounded on $\Hil$. Thus above lemma tells us that $K\Dom$ is a 
core for $\hDelta^{-1/2}$. On the
other hand, according to Lemma \reef{lm5.8} (with $\varepsilon=1$),
$K\Dom$ is contained in the set $\set{\eta''(b)}{b\in\Dom_1}$ and 
assertion follows.
\end{pf}\vs

We are now ready to prove the main theorem of this section.

\vspace{3mm}

\noindent {\it Proof of Theorem} \reef{th5.1}.\ \ First we prove that
\begin{equation}
\set{\eta''(b)}{b\in\Dom_0}\subset \Dom(\hDelta^{-1/2}) \label{dense}
\end{equation}
Let $b\in\Dom_0$. Then, Lemma \reef{lm5.8} shows that ${\mathcal 
R}_{\varepsilon}(b)\in\Dom_1$ for
any $\varepsilon >0$. Inserting $a={\mathcal 
R}_{\varepsilon}(b)-{\mathcal R}_{\varepsilon'}(b)$ in
(\reef{dmconj}) and using (\reef{regular2}), we get
\begin{equation}
\norm{\hDelta^{-1/2}\eta''({\mathcal R}_{\varepsilon}(b)) 
-\hDelta^{-1/2}\eta''({\mathcal
R}_{\varepsilon'}(b))}\Vs{4}  =\norm{\eta''({\mathcal 
R}_{\varepsilon}(\kappa^M(b))^*) -\eta''({\mathcal
R}_{\varepsilon'}(\kappa^M(b))^*)},
\label{Cauchy}
\end{equation}

\noindent for any $\varepsilon,\ \varepsilon'>0.$ According to 
Assertion 2 of Lemma \reef{lm5.8}, $\eta''({\mathcal
R}_{\varepsilon}(b))$ and $\eta''({\mathcal 
R}_{\varepsilon}(\kappa^M(b)^*))$ converge in norm to
$\eta''(b)$ and $\eta''(\kappa^M(b)^*)$ respectively as $\varepsilon$ 
tends to 0. By the Cauchy
criterion, $\hDelta^{-1/2}\eta''({\mathcal R}_{\varepsilon}(b))$ has 
a norm limit when
$\varepsilon$ tends to 0. Since $\hDelta^{-1/2}$ is closed, we see that
$\eta''(b)\in\Dom(\hDelta^{-1/2})$ and (\reef{dense}) follows.

The rest is easy. Remembering that $\Dom_1\subset\Dom_0$ (cf. 
Proposition \reef{pr5.5}) and using
Proposition \reef{pr5.10}, we see that the set 
$\set{\eta''(b)}{b\in\Dom_0}$ is a core for
$\hDelta^{-1/2}$.  Let $b\in\Dom_0$. Replacing $a$ in (\reef{dmconj}) 
by ${\mathcal
R}_{\varepsilon}(b)$ and setting $\varepsilon\to 0$, we obtain the 
norm equality (\reef{dmconj0}).

\qed

\vs

\par

Now, we shall consider the conjugate linear operator
\begin{equation}
{\hJ}: \hDelta^{-1/2}\eta''(b) \mapsto \eta''(\kappa^M(b)^*), \label{hJ}
\end{equation}
where $b$ runs over $\Dom_0$. By Theorem \reef{th5.1}, $\hJ$ is a 
densely defined isometry.
By the same letter we denote its continuous extension on the whole $\Hil$. Let
\begin{equation}
{\hF}={\hJ}\hDelta^{-1/2}.  \label{hF}
\end{equation}
Then
\begin{equation}
{\hF}\eta''(b)=\eta''(\kappa^M(b)^*) \label{fhat2}
\end{equation}
for all $b\in\Dom_0$. Moreover the set $\set{\eta''(b)}{b\in\Dom_0}$ 
is a core for $\hF$, because
it is a core for $\hDelta^{-1/2}$. In other words, $\hF$ is the 
closure of the mapping
\[
\eta''(b)\mapsto \eta''(\kappa^M(b)^*)
\]
where $b$ runs over $\Dom_0$.

The mapping $b\mapsto\kappa^M(b)^*$ is an involution acting on 
$\Dom_0$.  Therefore $\hF$ is also an
involution. Clearly (\reef{hF}) is the polar decomposition of $\hF$. 
Now the relation
$\hF^{-1}=\hF$ implies that
\begin{equation}
{\hJ}^2=1 \quad \mbox{and} \quad \hJ\hDelta\hJ=\hDelta^{-1}. \label{jhDelta}
\end{equation}
It means that the operator $\hJ$ is an antiunitary involution.

\begin{Thm}
\label{th5.11}
We have$:$
\begin{equation}
\label{th5.11a}
R(a)={\hJ}a^*{\hJ},
\end{equation}
\begin{equation}
\label{th5.11b}
(J\otimes \hJ)W(J\otimes \hJ)=W^*,
\end{equation}
\begin{equation}
\label{th5.11c}
J\hA J=\hA.
\end{equation}
  The first formula holds for all $a\in A$.
\end{Thm}

\begin{pf}
Ad \rf{th5.11a}. Let $a\in\Dom(\kappa)$ and $b\in\Dom_0$. Using 
(\reef{d0}), one can easily check that
$ab\in\Dom_0$. Therefore 
$a\eta''(b)=\eta''(ab)\in\Dom(\hDelta^{-1/2})$ and replacing $b$ in
(\reef{hJ}) by $ab$, we get
\[
\begin{array}{r@{\;}c@{\;}l}
\hJ\hDelta^{-1/2}a\eta''(b) 
&=&\hJ\hDelta^{-1/2}\eta''(ab)=\eta''(\kappa^M(ab)^*)\\
&=&\eta''(\kappa(a)^*\kappa^M(b)^*) =\kappa(a)^*\hJ\hDelta^{-1/2}\eta''(b).
\end{array}
\]
Since $\set{\eta''(b)}{b\in\Dom_0}$ is a core for $\hDelta^{-1/2}$, we see that
\begin{equation}
\hJ\kappa(a)^*\hJ\hDelta^{-1/2}\subset\hDelta^{-1/2}a. \label{jkappa}
\end{equation}
On the other hand, using (\reef{scaling}), one can show that
$\tau_{i/2}(a)Q^{-1}\subset Q^{-1}a$. For 
$x\in\Dom(J\rho^{1/2}JQ^{-1})$ we see that
\[
\tau_{i/2}(a)J\rho^{1/2}JQ^{-1}x=J\rho^{1/2}J\tau_{i/2}(a)Q^{-1}x
=J\rho^{1/2}JQ^{-1}ax=\hDelta^{-1/2}ax.
\]
Hence
$\tau_{i/2}(a)\hDelta^{-1/2}\subset\hDelta^{-1/2}a$. Comparing this 
inclusion with (\reef{jkappa}), we
obtain
\[
\begin{array}{rcl}
\hJ\kappa(a)^*\hJ=\tau_{i/2}(a).
\end{array}
\]
Replacing $a$ by $\kappa^{-1}(a)$, we get
\[
\begin{array}{rcl}
\hJ a^*\hJ=R(a)
\end{array}
\]
for all $a$ in the range of $\kappa$. To end the proof we notice that 
the range of
$\kappa$ equal to $\Dom(\kappa)^*$ is dense in $A$.\vs

Ad \rf{th5.11b}. Remembering that $R(W(x,y))=W(Jy,Jx)$ (cf. Assertion 
2 of Proposition \reef{pr2.2}) and
$R(a)=\hJ a^*\hJ$, we compute
\[
\begin{array}{r@{\;}c@{\;}l}
\IS{Jx\otimes\hJ z}{W(Jy\otimes\hJ u)}&=&\IS{\hJ z}{W(Jx,Jy)\hJ u}\\ \Vs{5}
&=&\IS{\hJ z}{R(W(y,x))\hJ u} =\IS{\hJ z}{\hJ 
W(y,x)^*u}=\IS{u}{W(y,x)z}\\ \Vs{5}
&=&\IS{y\otimes u}{W(x\otimes z)} =\IS{Jx\otimes\hJ z}{(J\otimes\hJ 
)W^*(y\otimes u)}.
\end{array}
\]
Thus we have $W(J\otimes\hJ)=(J\otimes\hJ)W^*$.\vs

Ad \rf{th5.11c}. By \rf{newAhat} it follows immediately from \rf{th5.11b}.
\end{pf}

\begin{Rem}\label{rmk5.13} By definitions $(\reef{hDelta})$ and $(\reef{hJ})$,
we see that $Q^{it}$ commute with $\hDelta$ and $\hJ$ for all $t\in\real$. It
implies that $Q$ and $\hDelta$ strongly
commute and
\begin{equation}
\label{hJQkom}
\hJ Q\hJ=Q^{-1}.
\end{equation}
\end{Rem}

In Section 3 we used the transpose mapping acting on $\skrL(\Hil)$. 
It was defined in the
following way :
\[
m^{\top}=Jm^*J
\]
for any $m\in\skrL(\Hil)$.  Formula \rf{th5.11c} shows that this transpose
maps $\hA$ onto itself.  Remembering that $\top$ is an involution, we have
$W^{\top\top\otimes\id}=W$.  Comparing this formula with (\reef{hRT}), we see
that $\hR=\top$. More explicitly
\begin{equation}
\label{hRT1}
\hR(b)=Jb^*J
\end{equation}
for all $b\in\hA$.

\begin{Prop}\label{pr5.12} Let $\eta_L$ be a linear mapping from $A$ 
to $\Hil$ with the domain
$\Dom(\eta_L)=R(\Dom(\eta)^*)=\set{a\in A}{R(a^*)\in\Dom(\eta)}$ 
defined by the formula
\begin{equation}
\eta_L(a)=\hJ\eta(R(a^*)).
\label{alpha}
\end{equation}
  Then $\eta_L$ is a closed GNS-mapping and the weight corresponding to
$\eta_L$ coincides with $h^L=h\comp R$.
\end{Prop}

\begin{pf} Let $a\in\Dom(\eta_L)$ and $b\in A$.  Then $R(a^*)\in\Dom(\eta)$ and
$R((ba)^*)=R(b^*)R(a^*)\in\Dom(\eta)$.  Therefore $ba\in\Dom(\eta_L)$ 
and using \rf{th5.11a} in the last
step, we obtain
\[
\eta_L(ba)=\hJ\eta(R(b^*)R(a^*))=\hJ R(b^*)\hJ\hJ\eta(R(a^*))=b\eta_L(a).
\]
It
shows that $\eta_L$ is a GNS-mapping.  The closedness of $\eta_L$ follows from
that of $\eta$.  To show the last statement we notice that
\[
(h\comp R)(a^*a)<\infty\Leftrightarrow
h(R(a)R(a^*))<\infty\Leftrightarrow R(a^*)\in\Dom(\eta) 
\Leftrightarrow a\in\Dom(\eta_L).
\]
  Moreover, in this case
\[
(h\comp 
R)(a^*a)=\norm{\eta(R(a^*))}^2=\norm{\hJ\eta(R(a^*))}^2=\norm{\eta_L(a)}^2.
\]
\end{pf}\vs

Let
\begin{equation}
\label{vecpre}
\eta(a)={\sum_k}^\oplus a\Omega_k,
\end{equation}
where $a\in\Dom(\eta)$, be a vector presentation of the GNS map 
$\eta$ (see Definition \reef{Bsz} and Theorem
\reef{Bszesc}). Then for any
$a\in\Dom(\eta_L)$ we have:
\[
\eta_L(a)=\hJ\eta(R(a)^*)=\hJ\;{\sum_k}^\oplus R(a)^*\Omega_k=
\hJ\;{\sum_k}^\oplus \hJ a\hJ\Omega_k={\sum_k}^\oplus a\hJ\Omega_k.
\]
Clearly
\begin{equation}
\label{vecpreR}
\eta_L(a)={\sum_k}^\oplus a\hJ\Omega_k,
\end{equation}
is a vector presentation of the GNS map $\eta_L$. If \rf{vecpre} is 
exact, so is \rf{vecpreR}.

\begin{Prop}
\label{removeJ}
For any $a\in\Dom(\eta'')$ we have
\begin{equation}
\left(a\in\Dom(\eta_L'')\Vs{4}\right)\Longleftrightarrow
\left(\eta''(a)\in\Dom(J\rho^{1/2}J)\Vs{4}\right).
\label{alphabeta}
\end{equation}
If this is the case, then
\begin{equation}
\eta_L''(a)=\lambda^{-i/4}J\rho^{1/2}J\eta''(a).
\label{betabeta}
\end{equation}
Furthermore $\Dom(\eta'')\cap\Dom(\eta_L'')$ is a core for $\eta''$ 
and $\eta_L''$.
\end{Prop}

\begin{pf} We know that $R$ is normal.  Passing to the weak closures 
in (\reef{alpha}), we obtain
$\Dom(\eta_L'')=\set{a\in M}{R(a)^*\in\Dom(\eta'')}$ and
\[
\eta_L''(a)=\hJ\eta''(R(a)^*).
\]
  Using
Proposition \reef{pr5.4}, we see that
\[
\left(\eta''(a)\in\Dom(J\rho^{1/2}J)\Vs{4}\right)\quad \Longrightarrow\quad
\left(\begin{array}{c} a\in\Dom(\eta_L'')\quad {\rm and}\\
\norm{\eta_L''(a)}=\norm{J\rho^{1/2}J\eta''(a)}.\Vs{5}
\end{array}\right)
\]
Therefore there exists a linear isometry $V$ such that
\begin{equation}
\eta_L''(a)=VJ\rho^{1/2}J\eta''(a)
\label{delta}
\end{equation}
for any $a\in\Dom(\eta'')$ such that 
$\eta''(a)\in\Dom(J\rho^{1/2}J)$.  We have to show that
$V=\lambda^{-i/4}1$.\vs

Let $b\in\Dom_1$ (cf. (\reef{d1})).  Then $\tau_{i/2}(b)\in\Dom(\eta'')$ and
$\eta''(\tau_{i/2}(b))=\lambda^{i/4}Q^{-1}\eta''(b)\in\Dom(J\rho^{1/2}J)$. 
Replacing $a$ in
(\reef{delta}) by $\tau_{i/2}(b)$ and using (\reef{hDelta}), we obtain
\begin{equation}
\eta_L''(\tau_{i/2}(b))=\lambda^{i/4}VJ\rho^{1/2}JQ^{-1}\eta''(b)
=\lambda^{i/4}V\hDelta^{-1/2}\eta''(b).
\label{epsilon}
\end{equation}
On the other hand
\[
\eta_L''(\tau_{i/2}(b))=\hJ\eta''((R\comp\tau_{i/2}(b))^*)=\hJ\eta''(\kappa^M(b)^*).
\]
We know that
$\hJ$ is an involution.  Combining the above formula with (\reef{hJ}), we get
$\eta_L''(\tau_{i/2}(b))=\hDelta^{-1/2}\eta''(b)$.  Therefore from 
(\reef{epsilon}) we have
$V\hDelta^{-1/2}\eta''(b)=\lambda^{-i/4}\hDelta^{-1/2}\eta''(b)$.

By Proposition \reef{pr5.10} the set $\set{\eta''(b)}{b\in\Dom_1}$ is 
a core for $\hDelta^{-1/2}$. Hence
$\set{\hDelta^{-1/2}\eta''(b)}{b\in\Dom_1}$ is dense in $\Hil$, 
$V=\lambda^{-i/4}1$ and (\reef{betabeta})
follows.\vs

Assume now that $a\in\Dom(\eta'')$.  By Lemma \reef{lm5.6}, we see that
$a\Gamma_{\varepsilon}\in\Dom(\eta'')$ and
\[\eta''(a\Gamma_{\varepsilon})=J\Gamma_{\varepsilon}'J\eta''(a)\in\Dom(J\rho^{1/2}J).\] 
Using the
part of the proposition proved so far, we obtain 
$a\Gamma_{\varepsilon}\in\Dom(\eta_L'')$ and
\[\eta_L''(a\Gamma_{\varepsilon})=\lambda^{-i/4}J\rho^{1/2}JJ\Gamma_{\varepsilon}'J\eta''(a).\] 
If
moreover $a\in\Dom(\eta_L'')$, then $R(a)^*\in\Dom(\eta'')$ and
\begin{eqnarray*}
\eta_L''(a\Gamma_{\varepsilon})
&=&\hJ\eta''(R(a\Gamma_{\varepsilon})^*)=\hJ\eta''(R(a)^*\Gamma_{\varepsilon})\\ 
&=&\hJ
J\Gamma_{\varepsilon}'J\eta''(R(a)^*)=\hJ 
J\Gamma_{\varepsilon}'J\hJ\eta_L''(a).
\end{eqnarray*} It shows that \[\hJ 
J\Gamma_{\varepsilon}'J\hJ\eta_L''(a)=\lambda^{-i/4}J\rho^{1/2}J
J\Gamma_{\varepsilon}'J\eta''(a).\] When $\varepsilon\to 0$, then
$J\Gamma_{\varepsilon}'J\eta''(a)\to\eta''(a)$ and $\hJ
J\Gamma_{\varepsilon}'J\hJ\eta_L''(a)\to\eta_L''(a)$.  Since 
$J\rho^{1/2}J$ is closed, we see that
$\eta''(a)\in\Dom(J\rho^{1/2}J)$.  The proof of (\reef{alphabeta}) is 
complete.\vs

By Lemma \reef{lm5.8} the set $\Dom_1$ is a core for $\eta''$. 
According to definition (\reef{d1}),
$\eta''(a)\in\Dom(J\rho^{1/2}J)$ for any $a\in\Dom_1$.  Using 
(\reef{alphabeta}), we see that
$\Dom_1\subset\Dom(\eta_L'')$.  Therefore 
$\Dom_1\subset\Dom(\eta'')\cap\Dom(\eta_L'')$ and hence
$\Dom(\eta'')\cap\Dom(\eta_L'')$ is a core for $\eta''$.  By the 
symmetry it is also a core for
$\eta_L''$.
\end{pf}\vs

With the vector presentations \rf{vecpre} and \rf{vecpreR}, formula 
\rf{betabeta} takes the form:
\[
{\sum_k}^\oplus 
a\hJ\Omega_k=\lambda^{-i/4}J\rho^{1/2}J\;{\sum_k}^\oplus a\Omega_k.
\]
We shall need the following extended version of this formula

\begin{Prop}
\label{removeJx}
Let $\Hil_1$ be a Hilbert space, $c\in 
M\;\overline{\otimes}\;\skrL(\Hil_1)$ and
$y\in\Hil_1$. Assume that the vector presentation \rf{vecpre} is 
exact, the series $\sum_k^\oplus
c(\Omega_k\otimes y)$ is norm convergent and that
\begin{equation}
\label{gwzdx}
{\sum_k}^\oplus c(\Omega_k\otimes y)\in\Dom(J\rho^{1/2}J\otimes 1).
\end{equation}
Then
\begin{equation}
\label{gwzdx2}
{\sum_k}^\oplus c(\hJ\Omega_k\otimes y)=
\lambda^{-i/4}\left(J\rho^{1/2}J\otimes 1\right){\sum_k}^\oplus 
c(\Omega_k\otimes y),
\end{equation}
where the series on the left hand side is norm convergent.
\end{Prop}

\begin{pf} We choose an orthonormal basis $\seq{e_s}{s=1,2,\dots}$ in 
$\Hil_1$. Then there exists a
sequence $\seq{c_s}{s=1,2,\dots}$ of elements of $M$ such that
\begin{equation}
\label{gwzwx2}
c(x\otimes y)={\sum_s}^\oplus c_sx\otimes e_s
\end{equation}
for all $x\in\Hil$. Inserting $x=\Omega_k$ and summing over $k$ we 
get a norm converging double
series
\[
{\sum_k}^\oplus c(\Omega_k\otimes y)=
{\sum_s}^\oplus\left({\sum_k}^\oplus c_s\Omega_k\otimes e_s\right).
\]
In particular for all $s$ the series $\sum_k^\oplus c_s\Omega_k$ are 
norm convergent. Therefore (cf. Definition \reef{defexact})
$c_s\in\Dom(\eta'')$ and
\[
{\sum_k}^\oplus c(\Omega_k\otimes y)=
{\sum_s}^\oplus\eta''(c_s)\otimes e_s.
\]
Assumption (\reef{gwzdx}) shows now that $\eta''(c_s)\in\Dom(J\rho^{1/2}J)$ and
$c_s\in\Dom(\eta''_L)$ for all $s$. Now we have:
\[
\begin{array}{c@{\;}l} {\displaystyle
\left(\lambda^{-i/4}J\rho^{1/2}J\otimes 1\right){\sum_k}^\oplus
c(\Omega_k\otimes y)}&= {\displaystyle
{\sum_s}^\oplus\lambda^{-i/4}J\rho^{1/2}J\eta''(c_s)\otimes e_s=
{\sum_s}^\oplus\eta''_L(c_s)\otimes e_s}\\
  &={\displaystyle {\sum_s}^\oplus{\sum_k}^\oplus c_s\hJ\Omega_k\otimes e_s=
{\sum_k}^\oplus
c(\hJ\Omega_k\otimes y)},\Vs{6}
\end{array}
\]
where in the last step we used (\reef{gwzwx2}).
\end{pf}\vs

We end this section with the following interesting
\begin{Thm}
\label{newAx}
The correspondence
\[
(a\otimes 1)\delta(b)\,\longmapsto\,\delta(a)(b\otimes 1),
\]
  \oldversion{where $a,b\in A$, extends to a closed linear mapping $\Psi$
acting on
$A\otimes A$. By definition the linear span of $\set{(a\otimes
1)\delta(b)}{a,b\in A}$ is a core for
$\Psi$. The mapping $\Psi$ has trivial kernel$\,:$ $\ker\Psi=\{0\}$.}
\newversion{where $a,b\in A$, extends to a closed
linear mapping acting on $A\otimes A$.}
\end{Thm}
\begin{pf}
The statement follows immediately from the formula
\begin{equation}
\label{lastf}
\ITS{x\otimes z}{(a\otimes 1)\delta(b)}{W(y\otimes u)}
=\ITS{W(x\otimes\hJ Q^{-1}u)}{\delta(a)(b\otimes 1)}{y\otimes\hJ Qz},
\end{equation}
that holds for all $a,b\in A$, $x,y\in\Hil$, $z\in\Dom(Q)$ and 
$u\in\Dom(Q^{-1})$ (the reader should notice
that each of the vectors $x\otimes z$, $W(y\otimes u)$, 
$W(x\otimes\hJ Q^{-1}u)$ and
$y\otimes\hJ Qz$ run  over a linearly dense subset of $\Hil\otimes\Hil$).\vs

To prove \rf{lastf} we start with the formula \rf{manageable} proving 
the manageability
of the Kac-Takesaki operator. The right hand side may be rewritten in 
the following way:
\[
\begin{array}{r@{\;}c@{\;}l}
\IS{x\otimes z}{W(y\otimes u)}&=&\IS{Jy\otimes Qz}{W^*(Jx\otimes Q^{-1}u)}\\
&=&\IS{(J\otimes\hJ)(y\otimes\hJ Qz)}{W^*(Jx\otimes Q^{-1}u)}\Vs{6}\\
&=&\IS{(J\otimes\hJ)W^*(Jx\otimes Q^{-1}u)}{y\otimes\hJ Qz}\Vs{6}\\
&=&\IS{W(J\otimes\hJ)(Jx\otimes Q^{-1}u)}{y\otimes\hJ Qz}\Vs{6}\\
&=&\IS{x\otimes\hJ Q^{-1}u}{W^*(y\otimes\hJ Qz)}\Vs{6},
\end{array}
\]
where in the fourth step we used \rf{th5.11b}. Replacing in the above 
formula $x$ by $a^*x$ and $y$ by $by$ we obtain
\[
\ITS{x\otimes z}{(a\otimes 1)W(b\otimes 1)}{y\otimes u}
=\ITS{x\otimes\hJ Q^{-1}u}{(a\otimes 1)W^*(b\otimes 1)}{y\otimes\hJ Qz}
\]
and \rf{lastf} follows (cf. \rf{coproduct}).

\end{pf}\vs

\section{A GNS-Mapping on $\hA$} \label{sek6}

One of the aim of this paper is to show that $(\hA,\hdelta)$ is a weighted Hopf
$C^*$-algebra.  To this end we have to introduce a right invariant 
weight $\hh$ on $\hA$. It will be related to a
GNS-mapping $\heta:\hA\to\Hil$.  However the definition of $\heta$ is 
not simple. We
start with the following technical result, in which the commutant 
$\eta'$ of the GNS mapping $\eta$ is used.

\begin{Prop}\label{pr6.1} Let $z\in\Dom(\hF)$.  Then for any $a,\ 
b\in\Dom(\eta')$ and
$x,\ y\in\Hil$ we have $:$
\begin{equation}
\IS{z\otimes y}{W(a^*x\otimes\eta'(b))} 
=\IS{b^*y\otimes\eta'(a)}{W^*(\hF z\otimes x)}. \label{heta}
\end{equation}
\end{Prop}

\begin{pf} We may assume that $z=\eta''(c)$ for some $c\in\Dom_0$.  Then
$\hF z=\eta''(\kappa^M(c)^*)$ by Proposition \reef{lm4.2}. 
Furthermore, we may assume that $x,\
y\in\Dom(Q)$ and $\eta'(a),\ \eta'(b)\in\Dom(Q^{-1})$. For any $d\in A$ we set
\[
\varphi_1(d)=(x|d\eta'(a)), \quad
\varphi_2(d)=(Q^{-1}\eta'(a)|R(d)Qx)
\]
\[
\psi_1(d)=(y|d\eta'(b)), \quad \psi_2(d)=(Q^{-1}\eta'(b)|R(d)Qy).
\]
Then $\varphi_1,\ \varphi_2,\ \psi_1,\ \psi_2\in A_*$,
$\varphi_2=\varphi_1^*\comp\kappa$, $\psi_2=\psi_1^*\comp\kappa$ and 
hence by (\reef{wadjoint})
$W_{\varphi_1}^*=W_{\varphi_2}$ and
$W_{\psi_1}^*=W_{\psi_2}$.  With this notation, the left hand side of 
(\reef{heta})
\begin{eqnarray*}
LHS&=&(\eta''(c)|W_{\psi_1}a^*x)=(W_{\psi_2}\eta''(c)|a^*x)=(\eta''(\psi_2*c)|a^*x)\\
&=&(a\eta''(\psi_2*c)|x)=((\psi_2*c)\eta'(a)|x)=\overline{\varphi_1(\psi_2*c)}\\&=&\overline{(\varphi_1*\psi_2)(c)},
\end{eqnarray*} where $\varphi_1*\psi_2\in A_*$. Replacing in the 
above computation $x,\ y,\ a,\ b,\
c$ by
$y,\ x,\ b,\ a,\ \kappa^M(c)^*$ respectively, we see that
$\psi_1*\varphi_2\in A_*$ and the complex conjugate of the right hand 
side of (\reef{heta}) equals
$\overline{(\psi_1*\varphi_2)(\kappa^M(c)^*)}$. For the proof of 
(\reef{heta}) it remains to show that
\[
(\psi_1*\varphi_2)(\kappa^M(c)^*)=\overline{(\varphi_1*\psi_2)(c)}.
\]
Since
\[
(W_{\psi_1*\varphi_2})^* =(W_{\psi_1}W_{\varphi_2})^* 
=W_{\varphi_2}^*W_{\psi_1}^*
=W_{\varphi_1}W_{\psi_2} =W_{\varphi_1*\psi_2},
\]
we have
$\varphi_1*\psi_2=(\psi_1*\varphi_2)^*\comp\kappa$ by Proposition 
\reef{pr2.14}. Hence
$(\psi_1*\varphi_2)^*\comp\kappa^M|_A=(\psi_1*\varphi_2)^*\comp\kappa$ 
is $\sigma$-weakly continuous.
It shows that
\[
\begin{array}{r@{\;}c@{\;}l}
(\psi_1*\varphi_2)(\kappa^M(c)^*) 
&=&\overline{(\psi_1*\varphi_2)^*(\kappa^M(c))}
=\overline{((\psi_1*\varphi_2)^*\comp\kappa^M)(c)}\\ 
&=&\overline{(\psi_1*\varphi_2)^*\comp\kappa(c)}
=\overline{(\varphi_1*\psi_2)(c)}\Vs{5}
\end{array}
\]
and (\reef{heta}) follows.
\end{pf}\vs

Relation (\reef{heta}) may be rewritten in the following way :
\[
(y|W(z,a^*x)\eta'(b))
=\IS{b^*y}{W_{\omega_{x,\eta'(a)}}^*\hF z}
\]
It shows that
\[
W(z,a^*x)\eta'(b)=bW_{\omega_{x,\eta'(a)}}^*\hF z.
\]
This relation holds for all $b\in\Dom(\eta')$. Therefore
\[ W(z,a^*x)\in\Dom(\eta'') \quad \mbox{and} \quad
\eta''(W(z,a^*x))=W_{\omega_{x,\eta'(a)}}^*\hF z.
\] Remembering that $W(z,a^*x)\in A$ and $\eta''|_A=\eta$, we obtain

\begin{Prop}\label{pr6.2} Let $z\in\Dom(\hF)$, $a\in\Dom(\eta')$ and 
$x\in\Hil$.  Then
$W(z,a^*x)\in\Dom(\eta)$ and
\[
\eta(W(z,a^*x))=W_{\omega_{x,\eta'(a)}}^*\hF z.
\]
\end{Prop}

It should be remarked that the vectors on the right hand side form a 
dense subset of $\Hil$.  Indeed,
the range of $\hF$ is dense in $\Hil$, the set of operators 
$W_{\omega_{x,\eta'(a)}}^*$ is dense in
$\hA$ and
$\hA$ is a non degenerate $C^*$-algebra acting on $\Hil$.\vs

Let $\Phi$ be the set :
\[
\Phi=\set{(x,y)\in\Hil\times\Hil}{W(x,y)\in\Dom(\eta)}.
\]
By Proposition \reef{pr6.2} and the above
remark, $\Phi$ is dense in $\Hil\times\Hil$  and
$\set{\eta(W(x,y))}{(x,y)\in\Phi}$ is dense in $\Hil$.  Therefore the relation
\begin{equation}
(y|ax)=(\eta(W(x,y))|z) \quad \mbox{for all} \quad (x,y)\in\Phi. 
\label{dgnsmap}
\end{equation}
gives a one to one correspondence between elements $a\in\hA$ and 
vectors $z\in\Hil$.
This correspondence will be denoted by
$\heta$ : we write
\[
z=\heta(a),
\]
if the relation (\reef{dgnsmap}) holds.  Clearly $\heta$ is a linear 
mapping from $\hA$ to
$\Hil$ and the domain $\Dom(\heta)$ consists of all elements 
$a\in\hA$ such that there
exists $z\in\Hil$ satisfying (\reef{dgnsmap}).\vs

Let $\varphi\in A^*$.  We say that $\varphi$ is $L^2$-{\rm bounded} 
if there exists
$z_{\varphi}\in\Hil$ such that
\[
\varphi(a)=(z_{\varphi}|\eta(a))
\]
for all $a\in\Dom(\eta)$.

\begin{Prop}\label{pr6.3} Let $\varphi\in A^*$ be an $L^2$-bounded 
functional.  Then $\varphi$ is
$\sigma$-weakly continuous.
\end{Prop}

\begin{pf} At first we notice that for any $b\in\Dom(\eta)$ the 
functional $b\varphi$ is
$\sigma$-weakly continuous. Indeed for any $a\in A$, $ab\in\Dom(\eta)$ and
\[
(b\varphi)(a)=\varphi(ab)=(z_{\varphi}|\eta(ab))=(z_{\varphi}|a\eta(b)).
\]
The set of $\sigma$-weakly continuous functionals is norm closed. 
Inserting $b=e_\alpha$ (where $\seq{e_\alpha}{}$ is an approximate 
identity for $A$) and using Appendix \reef{A}, we see that $\varphi$ 
is $\sigma$-weakly continuous.
\end{pf}\vs

Remark that the set of all $L^2$-bounded functionals is norm dense in 
$A_*$ (cf. the proof of
Assertion 1 of Proposition \reef{pr6.5}).

\begin{Thm}\label{th6.4}
\begin{enumerate}
  \item The domain $\Dom(\heta)$ is dense in $\hA$
  and the range of $\heta$ is dense in $\Hil$.
  \item The set $\Dom(\heta)$ is a left ideal in $\hA$ and $\heta(ba)=b\heta(a)$
  for all $b\in\hA$ and $a\in\Dom(\heta)$.
  \item The mapping $\heta$ of $\hA$ to $\Hil$ is injective and closed 
with respect to the strong topology on $\hA$
  and the norm topology on $\Hil$.
  \item Let $\varphi$ be an $L^2$-bounded functional on $A$.
  Then $W_{\varphi}^*\in\Dom(\heta)$ and $\heta(W_{\varphi}^*)=z_\varphi$.
\end{enumerate}
\end{Thm}

\begin{pf} \Vs{2}

Ad 4. Using Proposition \reef{pr6.3} we see that $\varphi\in A_*$ and 
hence $W_{\varphi}\in \hA$. For any $(x,y)\in\Phi$ we have
\begin{eqnarray*} 
(y|W_{\varphi}^*x)&=&\overline{(x|W_{\varphi}y)}=\overline{\varphi(W(x,y))}\\
&=&\overline{(z_\varphi|\eta(W(x,y)))}=(\eta(W(x,y))|z_\varphi).
\end{eqnarray*}
It shows that $W_{\varphi}^*\in\Dom(\heta)$ and 
$\heta(W_{\varphi}^*)=z_\varphi$.\vs

Ad 1. Now, let $y\in\Hil$, $c\in\Dom(\eta')$ and
\begin{equation}
\varphi(b)=(y|b\eta'(c)) \quad \mbox{for} \quad b\in A. \label{Lbdd}
\end{equation}
Then $\varphi\in A_*$ and for any $b\in\Dom(\eta)$, 
$\varphi(b)=(c^*y|\eta(b))$. By Assertion 4,
$W_{\varphi}^*\in\Dom(\heta)$ and $\heta(W_{\varphi}^*)=c^*y$. The 
set of functionals $\varphi$ of
the form (\reef{Lbdd}) is linearly dense in $A_*$. Therefore 
$\Dom(\heta)$ is dense
in $\hA$. Moreover, the set of vectors of the form $c^*y$ (where 
$c\in\Dom(\eta')$ and
$y\in\Hil$) is dense in $\Hil$. So is the range of $\heta$.\vs

Ad 2. Let $(x,y)\in\Phi$ and $\varphi\in A_*$. Then for any $\psi\in 
A^*$ we have
\begin{eqnarray*}
\psi(\varphi*W(x,y))&=&(\psi*\varphi)(W(x,y))=(x|W_{\psi*\varphi}y)\\
&=&(x|W_{\psi}W_{\varphi}y)=\psi(W(x,W_{\varphi}y)).
\end{eqnarray*}
It shows that $\varphi*W(x,y)=W(x,W_{\varphi}y)$. Since 
$W(x,y)\in\Dom(\eta)$, we see
that $W(x,W_{\varphi}y)\in\Dom(\eta)$. Hence 
$(x,W_{\varphi}y)\in\Phi$.  We also have
\[
\eta(W(x,W_{\varphi}y))=W_{\varphi}\eta(W(x,y)).
\]
Let $a\in\Dom(\heta)$.  Then for any
$(x,y)\in\Phi$ we have $(\eta(W(x,y))|\heta(a))=(y|ax)$. Replacing 
$(x,y)$ by $(x,W_{\varphi}y)$ we
obtain
\[
(W_{\varphi}\eta(W(x,y))|\heta(a))=(W_{\varphi}y|ax).
\]
Therefore
$(\eta(W(x,y))|W_{\varphi}^*\heta(a))=(y|W_{\varphi}^*ax)$. It shows that
$W_{\varphi}^*a\in\Dom(\heta)$ and 
$\heta(W_{\varphi}^*a)=W_{\varphi}^*\heta(a)$. Remembering that
$\set{W_{\varphi}^*}{\varphi\in A_*} 
=\set{(\id\otimes\varphi)(W^*)}{\varphi\in\skrL(\Hil)_*}$ is
dense in $\hA$ (cf. (\reef{newAhat})), we obtain Assertion 2 in full 
generality.\vs

Ad 3. It follows immediately from $(\reef{dgnsmap})$.
\end{pf}\vs

In this way we have shown that $\heta:\hA\to\Hil$ is a closed, 
densely defined GNS-mapping.
In what follows $\heta':\hA'\to\Hil$ will denote the commutant of $\heta$.\vs

If $\varphi\in A^*$ is $L^2$-bounded, then according to Assertion 4,
$W_{\varphi}^*\in\Dom(\heta)$ and $\varphi(b) 
=(\heta(W_{\varphi}^*)|\eta(b))$ for any
$b\in\Dom(\eta)$. Remembering that $\varphi$ is $\sigma$-weakly 
continuous (the same letter will denote its $\sigma$-weakly continuous
extension to $A''$) and that $\eta''$ is the extension of $\eta$ with 
$\Dom(\eta)$ being a core for $\eta''$ (cf. Theorem 
\reef{podkomutant}), we obtain
\begin{equation}
\varphi(b)=(\heta(W_{\varphi}^*)|\eta''(b)) \label{Fourier}
\end{equation}
for any $b\in\Dom(\eta'')$. Let
\begin{equation}
\hA_0=
\set{W_{\varphi}^*}{\varphi,\ \varphi^*\comp\kappa\in A^*\
\mbox{and}\ \varphi,\ \varphi^*\comp\kappa\ \mbox{are}\ L^2\mbox{-bounded}}
\label{dhalg}
\end{equation}
and
\begin{equation}
\hA_1=
\set{W_{\varphi}^*}{\varphi,\ \varphi^*\comp\kappa\in A^*\
\mbox{and}\ \varphi\ \mbox{is}\ L^2\mbox{-bounded}}.\hspace*{12mm}
\label{dhalg1}
\end{equation}

\begin{Prop}\label{pr6.5}
\begin{enumerate}
\item $\hA_0$ is a norm dense linear subset of $\hA$.
\item $\hA_0$ is a $*$-subalgebra of $\skrL(\Hil)$.
\item $\hA_0\subset\Dom(\heta)$ and $\heta(\hA_0)$ is dense in $\Hil$.
\item If $a\in \hA_0$, then $\heta(a)\in\Dom(\hF^*)$ and
$\hF^*\heta(a)=\heta(a^*)$.
\end{enumerate}
\end{Prop}

\begin{pf} Clearly $\hA_0\subset\hA_1$. Assume that $\varphi,\ \psi$ 
are $L^2$-bounded functionals in $A^*$ such that
$\varphi^*\comp\kappa$,\ $\psi^*\comp\kappa\in A^*$.  Then using 
(\reef{wproduct}) and (\reef{wadjoint}) we
obtain
\[
W_{\varphi}W_{\psi}^*=(W_{\psi}W_{\varphi^*\comp\kappa})^* 
=W_{\psi*(\varphi^*\comp\kappa)}^* \;,
\]
\[
W_{\varphi}W_{\psi}^*=W_{\varphi}W_{\psi^*\comp\kappa} 
=W_{\varphi*(\psi^*\comp\kappa)}.
\]
Let $\rho=\psi*(\varphi^*\comp\kappa)$.  Using Proposition 
\reef{pr2.14}, we see that
$\rho^*\comp\kappa=\varphi*(\psi^*\comp\kappa)\in A^*$. We shall 
prove that $\rho$
and $\rho^*\comp\kappa$ are $L^2$-bounded. Let $a\in\Dom(\eta)$. Then 
(cf. \rf{pr2.0})
$(\varphi^*\comp\kappa)*a\in\Dom(\eta)$ and
\begin{eqnarray*}
\rho(a) &=&(\psi*(\varphi^*\comp\kappa))(a)=\psi((\varphi^*\comp\kappa)*a)
=(z_{\psi}|\eta((\varphi^*\comp\kappa)*a))\\
&=&(z_{\psi}|W_{\varphi^*\comp\kappa}\eta(a))=(W_{\varphi}z_{\psi}|\eta(a)).
\end{eqnarray*}
Hence $\rho$ is $L^2$-bounded. Exchanging the
role of $\varphi$ and $\psi$ in the above argument, we see that 
$\rho^*\comp\kappa$ is $L^2$-bounded.
Hence $W_{\varphi}W_{\psi}^*=W_\rho^*\in \hA_0$. In this way we showed that
\begin{equation}
\label{fund}
\hA_1^*\hA_1\subset\hA_0\subset\hA_1.
\end{equation}

Ad 1. Let $y\in\Dom(Q)$ and $c\in\Dom(\eta')$ such that 
$\eta'(c)\in\Dom(Q^{-1})$
and
$\varphi=\omega_{y,\eta'(c)}$. Then for any $a\in\Dom(\eta)$ we have
\begin{equation}
\label{naka}
\varphi(a)=(y|a\eta'(c))=(y|c\eta(a))=(c^*y|\eta(a)).
\end{equation}
It shows that $\varphi$ is $L^2$-bounded.
On the other hand, $\varphi^*=\omega_{\eta'(c),y}$ and
\[
\varphi^*\comp\kappa=\omega_{\eta'(c),y}\comp\tau_{i/2}\comp R 
=\omega_{Q^{-1}\eta'(c),Qy}\comp R\in
A^*.\
\]
Therefore $W_{\omega_{y,\eta'(c)}}^*\in\hA_1$. One can easily check 
that the set of functionals
$\omega_{\eta'(b),\eta'(c)}$ (where $b$ and $c$ satisfy the above 
condition) is linearly dense in
$\skrL(\Hil)_*$ by Proposition \reef{pr4.1}. Hence $\hA_1$ is dense 
in $\hA$ (cf. (\reef{newAhat})). The first inclusion
in \rf{fund} shows now that $\hA_0$ is dense in $\hA$.\vs

Ad 2. If $W_{\varphi}^*\in\hA_0$, then $W_{\varphi^*\comp\kappa}^*\in\hA_0$ and
$W_{\varphi}=W_{\varphi^*\comp\kappa}^*$ by (\reef{wadjoint}).  Hence 
$\hA_0$ is invariant under the
involution. Formula \rf{fund} shows now that $\hA_0$ is a subalgebra.\vs

Ad 3. Statement 4 of Theorem \reef{th6.4} shows that 
$\hA_1\subset\Dom(\heta)$. Let $y,\ c$ and $\varphi$ be as in the
proof of Assertion 1. Using Assertion 4 of Theorem \reef{th6.4} and 
formula \rf{naka} we see that
\begin{equation}
\label{naka1}
\heta(W^*_{\omega_{y,\eta'(c)}})=c^*y.
\end{equation}
Using Proposition \reef{pr4.1}, we find that $\heta(\hA_1)$ is dense in
$\Hil$. So is $\heta(\hA_0)$; by \rf{fund},
$\heta(\hA_0)$ contains $\heta(\hA^*_1\hA_1)=\hA^*_1\heta(\hA_1)$.\vs

Ad 4. If $a\in\hA_0$, then $a=W_{\varphi}^*$ and $a^*=W_{\psi}^*$, where
$\varphi,\ \psi$ are bounded and $L^2$-bounded functionals on $A$ such
that
$\psi=\varphi^*\comp\kappa$. Remembering that $\psi$ and $\varphi^*$ are
$\sigma$-weakly continuous and that $\kappa^M$ is a closure of $\kappa$,
we see that $\psi(b)=\varphi^*\comp\kappa^M(b)$ for any
$b\in\Dom(\kappa^M)$.

Let $b\in\Dom_0$ (cf. (\reef{d0})).  Then $b,
\kappa^M(b)^*\in\Dom(\eta'')$ and using (\reef{Fourier}), we have
\begin{eqnarray*} (\eta''(\kappa^M(b)^*)|\heta(a))
&=&\overline{\varphi(\kappa^M(b)^*)}=\varphi^*(\kappa^M(b))\\
&=&\psi(b)=(\heta(a^*)|\eta''(b)).
\end{eqnarray*} Since $\set{\eta''(b)}{b\in \Dom_0}$ is a core for $\hF$
and
$\hF\eta''(b)=\eta''(\kappa^M(b)^*)$, we obtain $\heta(a)\in\Dom(\hF^*)$
and
$\hF^*\heta(a)=\heta(a^*)$.
\end{pf}

\begin{Rem}\label{rmk6.6} \oldversion{To show that the density of
$\hA_0\Hil$ in $\Hil$}
\newversion{To show that
$\hA_0\Hil$ is dense in $\Hil$} we suppose that $x\in\Hil$ is orthogonal
to $W_{\omega_{\eta'(b),\eta'(c)}}^*\Hil$ for all $b,c\in\Dom(\eta')$ with
$\eta'(b)\in\Dom(Q)$ and $\eta'(c)\in\Dom(Q^{-1})$.  Then
\[(\eta'(b)|W(y,x)\eta'(c))=(W_{\omega_{\eta'(b),\eta'(c)}}^*y|x)=0\] and
hence $W(y,x)=0$ for all $y\in\Hil$.
\oldversion{Since $\norm{x}^21=\sum_{n\in{\Bbb
N}}W(e_n,x)^*W(e_n,x)=0$ by Proposition $\reef{pr3.2}$, it follows that
$x=0$.}
\newversion{Using Proposition \reef{pr3.2} we see that
\[\norm{x}^21=\sum_{n\in{\Bbb N}}W(e_n,x)^*W(e_n,x)=0\]
and $x=0$.} Thus $\hA_0\Hil$ is dense in $\Hil$ and hence
$\hA_0\heta(\hA_0)$ is dense in $\heta(\hA_0)$.
\end{Rem}

\begin{Prop}\label{pr6.7} Let $\ha\in\skrL(\Hil)$ and $z\in\Hil$. 
Then the following three
conditions are equivalent
$:$
\begin{enumerate}
\item $\ha\in\Dom(\heta')$ and $\heta'(\ha)=z$.
\item For any $\hb\in \hA_0$,
\begin{equation}
\hb z=\ha\heta(\hb).
\label{rbdd}
\end{equation}
\item $(z,x)\in\Phi$ and $\eta(W(z,x))=\ha^*x$ for all $x\in\Hil$.
\end{enumerate}
\end{Prop}

\begin{pf} Ad $1\Rightarrow 2$. It is clear.

Ad $2\Rightarrow 3$. Let $y,\ c$ and $\varphi$ be as in the proof of 
Assertion 1 of Proposition
\reef{pr6.5}. Then $W_{\varphi}^*\in\hA$. Moreover, 
$\varphi(a)=(c^*y|\eta(a))$ for all
$a\in\Dom(\eta)$ and $\heta(W_{\varphi}^*)=c^*y$ by \rf{naka1}. Replacing $\hb$
by $W_{\varphi}^*$ in (\reef{rbdd}), we obtain
\[
W_{\varphi}^*z=\ha c^*y.
\]
Let $x\in\Hil$.  Then
\[
(y|W(z,x)\eta'(c))=\varphi(W(z,x)) =(W_{\varphi}^*z|x)=(\ha c^*y|x)=(y|c\ha^*x)
\]
for all $y\in\Dom (Q)$. Therefore $W(z,x)\eta'(c)=c\ha^*x$ for any 
$c\in\Dom(\eta')$ such that
$\eta'(c)\in\Dom(Q^{-1})$. Using Proposition \reef{pr4.1}, we see 
that $W(z,x)\in\Dom(\eta'')$ (hence $(z,x)\in\Phi$) and
$\eta''(W(z,x))=\ha^*x$. Since $W(z,x)\in A$ and $\eta''|A=\eta$, we see that
$W(z,x)\in\Dom(\eta)$ and $\eta(W(z,x))=\ha^*x$ for all $x\in\Hil$.

Ad $3\Rightarrow 1$. Assume that $(z,x)\in\Phi$ and 
$\eta(W(z,x))=\ha^*x$ for all $x\in\Hil$. Let
$\hb\in\Dom(\heta)$.  Then
\[
(x|\ha\heta(\hb))=(\ha^*x|\heta(\hb))=(\eta(W(z,x))|\heta(\hb))=(x|\hb z),
\]
and hence we have $\hb
z=\ha\heta(\hb)$ for all $\hb\in\Dom(\heta)$. Thus $\ha\in\hA'$, 
$\ha\in\Dom(\heta')$ and
$\heta'(\ha)=z$.
\end{pf}

\begin{Cor}\label{cor6.8} The set $\hA_0$ is a core for $\heta$ 
$($resp. $\heta''$ $)$ with respect to
the strong topology in $\hA$ $($resp. $\hA''$ $)$ and the norm 
topology in $\Hil$.
\end{Cor}

\begin{pf} Let $\heta_0=\heta|_{\hA_0}$.  By Proposition \reef{pr6.7} 
(equivalence of the first two
conditions), the commutant of $\heta_0$ coincides with $\heta'$. 
Therefore $\heta_0''=\heta''$. By Theorem \reef{podkomutant}, 
$\hA_0=\Dom(\heta_0)$ is a core for $\heta''$. Consequently it is a 
core for $\heta$.
\end{pf}\vs

Now we shall show that the operator $\hDelta$ introduced in 
(\reef{hDelta}) is the dual modular
operator. Let $a,\ b\in\Dom(\heta')$. By Proposition \reef{pr6.7},
$W(\heta'(a),\heta'(b))\in\Dom(\eta)$ and
\[
\eta(W(\heta'(a),\heta'(b)))=a^*\heta'(b).
\]
Exchanging $a$ and $b$, we see that $W(\heta'(b),\heta'(a))\in\Dom(\eta)$ and
$\eta(W(\heta'(b),\heta'(a)))=b^*\heta'(a)$. Moreover, by Proposition 
\reef{pr2.6},
$W(\heta'(a),\heta'(b))\in\Dom(\kappa)$ and 
$\kappa(W(\heta'(a),\heta'(b)))^*=W(\heta'(b),\heta'(a))$.
Let $\Dom_0$ be the subset of $A''$ introduced in (\reef{d0}). Then 
$W(\heta'(a),\heta'(b))\in\Dom_0$,
and using (\reef{fhat2}), we see that $a^*\heta'(b)\in\Dom(\hF)$ and
\begin{equation}
\hF a^*\heta'(b)=b^*\heta'(a). \label{hF2}
\end{equation}

We shall use the Tomita-Takesaki theory for the von Neumann algebra 
$\hA''$ and the GNS-mapping
$\heta''$. In particular, $\hS_T$ and $\hF_T$ (the subscript $T$ refers to
Tomita-Takesaki) will denote the left and the right involution 
operators, respectively. By
definition, the linear span of $\set{\hb^*\heta(\ha)}{\ha,\ 
\hb\in\Dom(\heta)}$ is a core for $\hS_T$ and
\begin{equation}
\hS_T\hb^*\heta(\ha)=\ha^*\heta(\hb). \label{linv}
\end{equation}
Similarly the linear span of
$\set{b^*\heta'(a)}{a,\ b\in\Dom(\heta')}$ is a core for $\hF_T$ and
\begin{equation}
\hF_Tb^*\heta'(a)=a^*\heta'(b). \label{rinv}
\end{equation}
It is known that $\hF_T=\hS_T^*$. Formula (\reef{hF2}) shows that 
$\hF_T\subset \hF$.

\begin{Thm}
\label{th6.9}
The operator $\hF$ introduced by $(\reef{fhat2})$ coincides with $\hF_T$.
\end{Thm}

\begin{pf} We have to show that $\hF\subset\hF_T$. Let 
$x\in\Dom(\hF)$. Let us consider linear
operators
\begin{equation}
\label{Tdef}
\begin{array}{r@{\;}c@{\;}l}
T: \heta(a)&\longmapsto&ax, \\
T^{\flat}: \heta(b)&\longmapsto&b\hF x,\Vs{5}
\end{array}
\end{equation}
where $a$ and $b$ run over $\hA_0$.  By Statement 3 of Proposition 
\reef{pr6.5}, the domains
$\Dom(T)=\Dom(T^{\flat})$ are dense in $\Hil$.

According to Assertion 2 of Proposition \reef{pr6.5} $a^*b,\ 
b^*a\in\hA_0$ for any $a,\ b\in\hA_0$.
Using Assertion 4 of the same proposition, we have
\[
\begin{array}{r@{\;}c@{\;}l}
(\heta(b)|T\heta(a)) &=&(\heta(b)|ax)=(\heta(a^*b)|x)=(\hF^*\heta(b^*a)|x)\\
&=&(\hF x|\heta(b^*a))=(b\hF x|\heta(a))=(T^{\flat}\heta(b)|\heta(a)).\Vs{5}
\end{array}
\]
It shows that $T^{\flat}\subset T^*$.  Therefore $T^*$ is densely 
defined and $T$ is
closable.  In what follows $T$ will denote the closure of the first 
operator in (\reef{Tdef}). Replacing $a$ by $b^*a$ in
(\reef{Tdef}), we obtain
\[
Tb^*\heta(a)=T\heta(b^*a)=b^*ax=b^*T\heta(a)
\]
for any $a,\ b\in\hA_0$. It shows that $b^*T\subset
Tb^*$.  Using this formula, one can easily show that $T$ is 
affiliated with the von Neumann algebra
$\hA'$.

Let $f$ be a continuous function on $\real_+$ such that 
$\set{tf(t^2)}{t\in\real_+}$ is bounded, and
$T_f=Tf(T^*T)$. Then $T_f\in\hA'$, $T_fy=f(TT^*)Ty$ for any 
$y\in\Dom(T)$, and $T_f^*z=f(T^*T)T^*z$
for any $z\in\Dom(T^*)$.  Replacing $y$ by $\heta(a)$, we obtain
\[
T_f\heta(a)=f(TT^*)ax=af(TT^*)x.
\]
This formula holds for all $a\in\hA_0$.  Then Proposition
\reef{pr6.7} shows now that
$T_f\in\Dom(\heta')$ and $\heta'(T_f)=f(TT^*)x$. Similarly, replacing 
$z$ by $\heta(b)$, we obtain
\[
T_f^*\heta(b)=f(T^*T)b\hF x=bf(T^*T)\hF x.
\]
This formula holds for all $b\in\hA_0$. Then
Proposition \reef{pr6.7} shows now that $T_f^*\in\Dom(\heta')$ and 
$\heta'(T_f^*)=f(T^*T)\hF x$.
Combining the two results, we see that $f(TT^*)x\in\Dom(\hF_T)$ and
\[
\hF_Tf(TT^*)x=f(T^*T)\hF x.
\]
Making $f$ converge to the constant function 1, we see that
$x\in\Dom(\hF_T)$ and
$\hF_Tx=\hF x$.   Thus we have $\hF\subset \hF_T$.  The converse inclusion
$\hF_T\subset \hF$ has already been known.
\end{pf}

\begin{Thm}
\label{th6.10}
The operators $\hDelta$ and $\hJ$ $($introduced by $(\reef{hDelta})$ and
$(\reef{hJ})$ respectively$)$ coincide with the modular operator and 
the modular conjugation operator
related to the GNS-mapping $\heta$.  In particular,
\[
\hJ\hA''\hJ=\hA',
\]
\[
\hDelta^{it}\hA''\hDelta^{-it}=\hA''
\]
hold for all $t\in\real$.
\end{Thm}

\begin{pf} See formula (\reef{hF}).
\end{pf}

\oldversion{
\begin{Rem}\label{rmk6.11} Combining two formulae $(\reef{dgnsmap})$ and
$(\reef{Fourier})$, we can deduce the following two new formulae
\[
(\eta(a)|\heta(\ha))=(\heta(\ha^*)|\eta(\kappa(a)^*))
\]
for $\ha\in\Dom(\heta)\cap\Dom(\heta)^*$ and for
$a\in\Dom(\eta)$ with
$\kappa(a)^*\in\Dom(\eta)$ $;$ and
\[
(\eta(a)|\heta(W_{\varphi}^*))=(\heta(W_{\varphi^*}^*)|\eta(a^*))
\]
for $a\in\Dom(\eta)\cap\Dom(\eta)^*$ and for $W_{\varphi}^*\in\Dom(\heta)$ with
$W_{\varphi^*}^*\in\Dom(\heta)$.
\end{Rem}
}

\begin{Prop}\label{pr6.121}
$A''\cap\hA''=\compl 1$.
\end{Prop}
\begin{pf}\ We know that $\hJ\hA''\hJ=\hA'$ and
$\hJ A''\hJ=A''$ (cf. \rf{th5.11a}). Therefore 
$\hJ(A''\cap\hA'')\hJ=A''\cap\hA'$ and our assertion reduces to
Statement 3 of Proposition \reef{pr6.12}.
\end{pf}

\begin{Prop}
\label{pr6.13}\Vs{2}

\begin{enumerate}
  \item $\rho^{it}\in M(A)$ for $t\in\real$ and
  the mapping $t\in\real\mapsto\rho^{it}\in M(A)$ is strictly continuous.
  \item $\delta(\rho^{it})=\rho^{it}\otimes\rho^{it}$ for $t\in\real$.
  \item The operator $\gamma$ in Proposition \reef{gammarho} is a scalar
  multiple of the identity.
\end{enumerate}
\end{Prop}

\begin{pf} Ad 1 and 2. At first we shall prove that
\begin{equation}
\delta^M(\rho^{it})=\rho^{it}\otimes\rho^{it}. \label{mrho0}
\end{equation}
Since $\hDelta^{it}=J\rho^{it}JQ^{2it}$ for $t\in\real$ and $J\hA J=\hA$
by \rf{th5.11c}, it follows that
$\rho^{it}=J\hDelta^{it}Q^{-2it}J$ implements an action of
$\real$ on $\hA''$ : $\rho^{-it}\hA''\rho^{it}\subset \hA''$. 
Therefore the right hand side
of the equation
\begin{equation}
(\rho^{-it}\otimes 1)\delta^M(\rho^{it}) =(\rho^{-it}\otimes 
1)W(\rho^{it}\otimes 1)W^*,
\label{mrho}
\end{equation}
commutes with $\hA'\otimes 1$. The left hand side clearly commutes with
$A'\otimes 1$. Since $A''\cap\hA''=\compl 1$ by Proposition 
\reef{pr6.121}, the element
in (\reef{mrho}) is of the form $1\otimes v_t$ for some unitary $v_t$. Hence
$\delta^M(\rho^{it})=\rho^{it}\otimes v_t$. It is easy to see that 
$\seq{v_t}{t\in\real}$ is a strongly
continuous one parameter unitary group. Now we use a formula
$\delta^M\comp R^M=\sigma\comp(R^M\otimes R^M)\comp\delta^M$. Since 
$R^M(\rho^{it})=\rho^{-it}\;$ by
Proposition \reef{pr4.5}, we find that 
$\delta^M(\rho^{it})=R^M(v_{-t})\otimes\rho^{it}$. Thus
$\rho^{it}\otimes v_t=R^M(v_{-t})\otimes\rho^{it}$. Therefore $v_t$ 
is a scalar multiple of
$\rho^{it}$. Hence there exists a positive scalar $\mu$ with 
$v_t=\mu^{it}\rho^{it}$. Inserting this
$v_t$ to the last formula, we get $\mu^{-it}=\mu^{it}$ for all 
$t\in\real$. Hence we see that $\mu=1$
and $v_t=\rho^{it}$.

Now formula $(\reef{mrho})$ takes the form 
$1\otimes\rho^{it}=(\rho^{-it}\otimes 1)W(\rho^{it}\otimes
1)W^*.$ All the factors on the right hand side belong to $M({\mathcal 
K}(\Hil)\otimes A)$ and depend
continuously on $t\in\real$ with respect to the strict topology. 
Therefore $\rho^{it}\in M(A)$ and
the mapping $t\mapsto\rho^{it}$ is continuous with respect to the 
strict topology.  At this moment one
can erase the superscript $M$ in (\reef{mrho0}).\vs

Ad 3. Due to \rf{oprho56} and \rf{opq56} we have: 
$\sigma_s(\rho^{it})=\gamma^{ist}\rho^{it}$ and
$\tau_s(\rho^{it})=\rho^{it}$ for all $t,s\in\real$.  Taking into 
account Assertion 2 of Lemma \reef{lm4.2} we
compute:
\[
\begin{array}{r@{\;}c@{\;}l}
\delta(\gamma^{ist})\delta(\rho^{it}) 
&=&\delta(\gamma^{ist}\rho^{it})=\delta(\sigma_s(\rho^{it}))\\
&=&(\sigma_s\otimes\tau_s)(\delta(\rho^{it})) 
=\sigma_s(\rho^{it})\otimes\tau_s(\rho^{it})\Vs{5}\\
&=&\gamma^{ist}\rho^{it}\otimes\rho^{it} =(\gamma^{ist}\otimes 
1)\delta(\rho^{it}).\Vs{5}
\end{array}
\]
  Thus we get $\delta(\gamma^{ist})=\gamma^{ist}\otimes 1$. Using 
Proposition \reef{pr6.12} we get
$\gamma^{ist}\in\compl 1$. Hence $\gamma$ is a positive scalar 
multiple of the identity.
\end{pf}\vs

In short, Assertion 1 of Proposition \reef{pr6.13} means that $\rho$ 
and $\rho^{-1}$ are affiliated
to the C$^*$-algebra $A$. Assertion 2 says that 
$\delta(\rho)=\rho\otimes\rho$. In other words
\begin{equation}
\label{deltarho}
W(\rho\otimes 1)W^*=\rho\otimes\rho.
\end{equation}

\begin{Prop}
We have
\begin{equation}
\label{gammalambda}
\gamma=\lambda^{-1}\,1,
\end{equation}
where $\lambda$ is the constant appearing in Definition $\reef{def1.2}$.
\end{Prop}

\begin{pf}
We shall use the first formula of \rf{oprho56}. It shows that 
$\sigma_s(\rho^{it})=\gamma^{ist}\rho^{it}$ for any $s,t\in\real$. 
Therefore
$\rho^{it}$ is entire analytic for $\seq{\sigma_s}{s\in\real}$ and 
$\sigma_{i/2}(\rho^{it})=\gamma^{-t/2}\rho^{it}$.  By Theorem 
\reef{19.08.2002d}, $b\rho^{it}\in\Dom(\eta)$ and
\[
\gamma^{t/2}\eta(b\rho^{it})=J\left(\rho^{it}\right)^*J\eta(b)=J\rho^{-it}J\eta(b)
\]
for any $b\in\Dom(\eta)$. Now we have:
\[
\hDelta^{-it}\eta(b)=Q^{-2it}J\rho^{-it}J\eta(b)=\gamma^{t/2}Q^{-2it}\eta(b\rho^{it})=(\lambda\gamma)^{t/2}\eta(\tau_{-t}(b\rho^{it})),
\]
where in the last step we used \rf{Qdef}.\vs

Let $\varphi\in A^*$ be $L^2$-bounded. Taking into account \rf{Fourier} we get:
\[
\begin{array}{r@{\;}c@{\;}l}
\IS{\hDelta^{it}\heta(W_{\varphi}^*)}{\eta(b)}
&=&\IS{\heta(W_{\varphi}^*)}{\hDelta^{-it}\eta(b)}
=(\lambda\gamma)^{t/2}\IS{\heta(W_{\varphi}^*)}{\eta(\tau_{-t}(b\rho^{it}))}\\
&=&(\lambda\gamma)^{t/2}\varphi(\tau_{-t}(b\rho^{it}))=(\lambda\gamma)^{t/2}(\rho^{it}\varphi\comp\tau_{-t})(b)\Vs{5}.
\end{array}
\]
This formula holds for all $b\in\Dom(\eta)$. It shows that
$\rho^{it}\varphi\comp\tau_{-t}$ is $L^2$-bounded and that
\[
\hDelta^{it}\heta(W_{\varphi}^*)
=(\lambda\gamma)^{t/2}\heta(W_{\rho^{it}\varphi\comp\tau_{-t}}^*).
\]
On the other hand, according to the Tomita-Takesaki theory (Statement 
1 of Theorem \reef{18.08.2002a}), $\hDelta^{it}\heta(W_{\varphi}^*) 
=\heta(\hsigma_t(W_{\varphi}^*))$, where
$\hsigma$ is the modular automorphism group related to the weight 
$\hh$ associated with $\heta$ (see \rf{hh} in the next section). 
Since $\heta$ is injective, we get
\[
\hsigma_t(W_{\varphi}^*) 
=(\lambda\gamma)^{t/2}W_{\rho^{it}\varphi\comp\tau_{-t}}^* .
\]
Remembering that $\hsigma_t$ and $\tau_t$ are implemented by 
$\hDelta^{it}$ and $Q^{2it}$ we obtain

\[
\begin{array}{r@{\;}c@{\;}l}
(\id\otimes\varphi)((\hDelta^{it}\otimes 1)W(\hDelta^{-it}\otimes 1))
&=&(\lambda\gamma)^{t/2}(\id\otimes(\rho^{it}\varphi\comp\tau_{-t}))(W)\\
&=&(\lambda\gamma)^{t/2}(\id\otimes\varphi) ((1\otimes 
Q^{-2it})W(1\otimes Q^{2it}\rho^{it})).\Vs{5}
\end{array}
\]
Hence we have
\[
(\hDelta^{it}\otimes(\lambda\gamma)^{-t/2}Q^{2it})W 
=W(\hDelta^{it}\otimes Q^{2it}\rho^{it}).
\]
In this formula, all entries except
$(\lambda\gamma)^{-t/2}$ are unitary operators. Therefore
$|(\lambda\gamma)^{-t/2}|=1$ and
$\gamma=\lambda^{-1}\,1$.
\end{pf}\vs

\begin{Rem}
Notice that the last formula of the above proof reduces now to
\begin{equation}
\label{lech}
(\hDelta^{it}\otimes Q^{2it})W =W(\hDelta^{it}\otimes Q^{2it}\rho^{it}).
\end{equation}
\end{Rem}
\vs

At the end of this section we investigate the behavior of the right 
invariant Haar weight with respect to left translations. In the 
following proposition we use a vector presentation of the GNS map 
$\eta$:
\[
\eta(a)={\sum_{k}}^\oplus\,a\Omega_k
\]
for any $a\in\Dom(\eta)$. Then (cf. \rf{KT2}), $W(\eta''(a)\otimes
y)={\sum_{k}}^\oplus\,\delta(a)(\Omega_k\otimes y)$ for any 
$a\in\Dom(\eta'')$ and $y\in\Hil$.

\begin{Prop}
Let $a\in\Dom(\eta'')$ and $x\in\Dom\left(\rho^{1/2}\right)$. Then
\begin{equation}
\label{ahoha}
{\sum_{k}}^\oplus\,\delta(a)(x\otimes\Omega_k)=
(\hJ\otimes\hJ)\hW^*(\hJ\otimes\hJ)\left(\rho^{1/2}x\otimes\eta''(a)\right).
\end{equation}
\end{Prop}
\begin{pf} We may assume that $a\in\Dom(\eta'')\cap\Dom(\eta''_L)$ 
(by Proposition \reef{removeJ},
the latter set is a core for $\eta''$). Using \rf{betabeta} we obtain:
\[
\begin{array}{r@{\;}c@{\;}l}
W\left(\eta''_L(a)\otimes\hJ\rho^{1/2}x\right)&=&
W\left(\lambda^{-i/4}J\rho^{1/2}J\eta''(a)\otimes\hJ\rho^{1/2}x\right)\\ \Vs{6}
&=&\lambda^{-i/4}(J\otimes\hJ)W^*\left(\rho^{1/2}\otimes\rho^{1/2}\right)(J\eta''(a)\otimes 
x)\\
\Vs{6} &=&\lambda^{-i/4}(J\otimes\hJ)\left(\rho^{1/2}\otimes 
1\right)W^*(J\eta''(a)\otimes x)\\
\Vs{6}&=&\left(\lambda^{-i/4}J\rho^{1/2}J\otimes 
1\right)W(\eta''(a)\otimes\hJ x)\\
\Vs{6}&=&\left(\lambda^{-i/4}J\rho^{1/2}J\otimes 1\right)
{\displaystyle {\sum_{k}}^\oplus\,\delta(a)(\Omega_k\otimes\hJ x)}
={\displaystyle {\sum_{k}}^\oplus\,\delta(a)(\hJ\Omega_k\otimes\hJ x)}.\\
\end{array}
\]
In this computation we used in the second and fourth steps formula 
\rf{th5.11b} and in
the third step formula \rf{deltarho}. At the end we used Proposition 
\reef{removeJx}. It shows
that the last series is convergent.\vs

Replacing in the above formula $a$ by $R(a)^*$ we obtain
\[
W\left(\hJ\eta''(a)\otimes\hJ\rho^{1/2}x\right)=
{\sum_{k}}^\oplus\,\delta(R(a)^*)(\hJ\Omega_k\otimes\hJ x).
\]
and
\[
\hW^*\left(\hJ\rho^{1/2}x\otimes\hJ\eta''(a)\right)=
{\sum_{k}}^\oplus\,\sigma\comp\delta(R(a)^*)(\hJ x\otimes\hJ \Omega_k),
\]
where $\sigma$ denotes the flip acting on $A\otimes A$.
Using Assertion 1 of Proposition \reef{pr2.10} and formula \rf{th5.11a} we
obtain:
\[
\sigma\comp\delta(R(a)^*)=(R\otimes R)\delta(a^*)=
(\hJ\otimes\hJ)\delta(a)(\hJ\otimes\hJ).
\]
Therefore
\[
\hW^*\left(\hJ\rho^{1/2}x\otimes\hJ\eta''(a)\right)=
{\sum_{k}}^\oplus\,(\hJ\otimes\hJ)\delta(a)(x\otimes\Omega_k)
\]
and \rf{ahoha} follows.

\end{pf}

\begin{Prop}
Let $\varphi\in A^*_+$ and $b\in A_+$. Assume that 
$\varphi(\rho)<\infty$ and $h(b)<\infty$. Then
\begin{equation}
\label{leftsh}
h(b*\varphi)=\varphi(\rho)h(b).
\end{equation}
\end{Prop}
It shows that $\rho$ is an analogue of the modular function on the 
group. It says, how the right
invariant Haar measure transforms under left shifts.

\begin{pf}
Let $a=b^{1/2}$. Computing the norm of the both sides of \rf{ahoha} we obtain:
\[
\sum_k\,\norm{\delta(a)(x\otimes\Omega_k)}^2=\left\norm{\rho^{1/2}x\right}^2 
h(a^*a)=\omega_{x,x}(\rho)h(b).
\]
On the other hand
\[
\sum_k\,\norm{\delta(a)(x\otimes\Omega_k)}^2=
\sum_k\,\ITS{x\otimes\Omega_k}{\delta(a^*a)}{x\otimes\Omega_k}=
\sum_k\,\ITS{\Omega_k}{b*\omega_{x,x}}{\Omega_k}=h(b*\omega_{x,x}).
\]
So we proved \rf{leftsh} for all $\varphi$ of the form $\omega_{x,x}$ 
($x\in\Dom(\rho^{1/2})$) i.e. for all
$\varphi\in A_{*+}$ with $\varphi(\rho)<\infty$. \vs

Let $\varphi\in A^*_+$. By the GNS construction one may assume that
$\varphi(b)=\IS{\Omega_\varphi}{\pi(b)\Omega_\varphi}$ for any $b\in 
A$. In this formula $\pi$ is a
representation of $A$ acting on a Hilbert space and $\Omega_\varphi$ is 
an element of this space.
Applying the method developed in the proof of Theorem \reef{newthm} 
to formula \rf{ahoha}, one
may show that
\[
{\sum_{k}}^\oplus\,(\pi\otimes\id)\delta(a)(\Omega_\varphi\otimes\Omega_k)=
(\pi\otimes\id)(\hJ\otimes\hJ)\hW^*(\hJ\otimes\hJ)\left(\pi(\rho^{1/2})
\Omega_\varphi\otimes\eta''(a)\right).
\]
Computing the norm of the both sides we obtain \rf{leftsh} in full 
generality. The details are left
to the reader.
\end{pf}\vs

\section{Dual weighted Hopf $C^*$-algebra and Duality} \label{sek7}

In Section \reef{sek3}, applying Theorem \reef{th3.8} to the 
multiplicative unitary $\hW$, we introduced the proper 
$C^*$-bialgebra $(\hA,\hdelta)$ with cancellation property together 
with the unitary antipode $\hR$ and the scaling group 
$\seq{\htau_t}{t\in\real}$ acting on $\hA$. In Section \reef{sek6} we 
constructed the closed densely defined GNS map $\heta$ from $\hA$ 
into $\Hil$. Let $\hh$ be the locally finite lower semicontinuous 
weight on $\hA$ related to the GNS map $\heta$ by Theorem 
\reef{wagaaa}: for any $b\in\hA$ we have
\begin{equation}
\label{hh}
\hh(b^*b)=\left\{
\begin{array}{c@{\hspace{9mm}}c}
\IS{\heta(b)}{\heta(b)}&\mbox{if }b\in\Dom(\heta)\\
+\infty&\mbox{otherwise.}\Vs{5}
\end{array}\right.
\end{equation}

  Let $a\in\Dom(\heta'')$ and $\heta''(a)=0$. Then there exists a 
sequence $\seq{a_n}{n\in\natu}$ of elements of $\Dom(\heta)$
converging strongly to $a$ such that the sequence 
$\seq{\heta(a_n)}{n\in\natu}$ converges in norm $0$. Taking into 
account
$(\reef{dgnsmap})$ we obtain: $\IS{y}{ax}=\IS{\eta(W(x,y))}{z}=0$ for 
all $(x,y)\in\Phi$. Remembering that $\Phi$ is dense in
$\Hil\times\Hil$ we see that $a=0$. It shows that $\heta''$ is 
injective. Taking into account Lemma
\reef{lem:strictfaithful} we conclude that $\hh$ is strictly faithful.\vs

We shall prove that $(\hA,\hdelta)$ is a weighted Hopf $C^*$-algebra 
with $\hh$ playing the role of the right Haar
weight. To this end we shall use Theorem \reef{newthm} with 
$A,\delta,R,\tau,h$ and $W$ replaced by
$\hA,\hdelta,\hR,\htau,\hh$ and $\hW=\Sigma W^*\Sigma$ respectively. 
With this replacement formulae \rf{zalozenie} and
\rf{zalozenie1} take the form:
\begin{equation}
\label{zalozeniex}
\heta\left((\id\otimes\varphi)\hdelta(b)\right)=\left((\id\otimes\varphi)\hW\right)\heta(b)\vspace{-3mm}
\end{equation}
and
\begin{equation}
\label{zalozenie1x}
(\id\otimes\varphi\comp\hkappa)\hW=(\id\otimes\varphi)(\hW^*).
\end{equation}
We have also to verify that for all $t\in\real$,
\begin{equation}
\label{zalozenie2x}
\hh\comp\htau_t=\hlambda^t\hh,
\end{equation}
where $\hlambda>0$  is a fixed number independent of $t$.

\begin{Lem}
\label{lm7.1}
Let  $\varphi\in\hA^*$ and $b\in\Dom(\heta)$. Then 
$(\id\otimes\varphi)\hdelta(b)\in\Dom(\heta)$ and formula
\rf{zalozeniex} holds.
\end{Lem}

\begin{pf} According to Corollary \reef{cor6.8}, the set $\hA_0$ (cf. 
\rf{dhalg}) is a core for $\heta$. Therefore we may assume that
$b=W_{\psi}^*$ where $\psi$ is a $L^2$-bounded functional on $A^*$. 
Applying $\id\otimes\id\otimes\psi^*$ to the
both sides of \rf{deltahat} we get 
$\hdelta(b)=(\id\otimes\id\otimes\psi^*)\left(W^*_{13}W^*_{23}\right)$. 
Setting
$d=(\varphi\otimes\id)(W^*)\in M(A)$ we obtain
\[
\begin{array}{r@{\;}c@{\;}l}
(\id\otimes\varphi)\hdelta(b)
&=&(\id\otimes\varphi\otimes\psi^*)(W_{13}^*W_{23}^*)\\
&=&(\id\otimes\psi^*)(W^*(1\otimes d))\Vs{5}\\&=&(\id\otimes 
d\psi^*)(W^*)=W\Vs{3}^*_{\psi d^*}\Vs{5}.
\end{array}
\]
We have to show that
\begin{equation}
\label{teza71}
W\Vs{3}^*_{\psi d^*}\in\Dom(\heta)\hspace{5mm}{\rm 
and}\hspace{5mm}\heta\left(W\Vs{3}^*_{\psi
d^*}\right)=d\heta(b).
\end{equation}
We shall use formula \rf{Fourier}. For any $a\in\Dom(\eta)$, we have
\[
(\psi 
d^*)(a)=\psi(d^*a)=(\heta(W_{\psi}^*)|\eta(d^*a))=(d\heta(W_{\psi}^*)|\eta(a)).
\]
Now, \rf{teza71} follows immediately from Assertion 4 of Theorem \reef{th6.4}.
\end{pf}

\begin{Lem}
\label{lm7.2}
Formula \rf{zalozenie1x} holds for any $\varphi\in\hA^*$ such that 
$\varphi\comp\hkappa\in\hA^*$.
\end{Lem}
\begin{pf}
Let $\psi\in\skrL(\Hil)_*$.  Formula \rf{wantipode} applied to the manageable
multiplicative unitary $\hW$ shows that 
$(\psi\otimes\id)(\hW)\in\Dom(\hkappa)$ and
\[
\hkappa((\psi\otimes\id)(\hW))=(\psi\otimes\id)(\hW^*).
\]
Therefore we have
\[
\begin{array}{r@{\;}c@{\;}l}
\psi((\id\otimes\varphi\comp\hkappa)(\hW))
&=&(\varphi\comp\hkappa)((\psi\otimes\id)(\hW))\\ 
&=&\varphi((\psi\otimes\id)(\hW^*))
=\psi((\id\otimes\varphi)(\hW^*))\Vs{5}
\end{array}
\]
and \rf{zalozenie1x} follows.
\end{pf}\vs

Formula \rf{zalozenie2x} with $\hlambda=\lambda^{-1}$ follows 
immediately from the following

\begin{Prop}
\label{pr6.16}
Let $b\in\Dom(\heta)$. Then $\htau_t(b)\in\Dom(\heta)$ and
\[
\heta(\htau_t(b))=\lambda^{-t/2}Q^{2it}\heta(b)
\]
for all $t\in\real$.
\end{Prop}

\begin{pf} Assume that $b\in\Dom(\heta)$.  Then
\begin{equation}
(y|bx)=(\eta(W(x,y))|\heta(b)) \label{dgnsmap3}
\end{equation}
for all $(x,y)\in\Phi$.  Let $t\in\real$. According to the second 
assertion of Lemma \reef{lm2.13},
$\tau_{-t}(W(x,y))=W(Q^{-2it}x,Q^{-2it}y)$. Formula (\reef{scaling}) 
shows now that
$(Q^{-2it}x,Q^{-2it}y)\in\Phi$ and
\[
\eta(W(Q^{-2it}x,Q^{-2it}y))=\lambda^{-t/2}Q^{-2it}\eta((W(x,y)).
\]
Replacing $(x,y)$ in (\reef{dgnsmap3}) by $(Q^{-2it}x,Q^{-2it}y)$, we obtain
\[
(y|\htau_t(b)x)=(\eta(W(x,y))|\lambda^{-t/2}Q^{2it}\heta(b))
\]
for all $(x,y)\in\Phi$. It shows that
$\htau_t(b)\in\Dom(\heta)$ and
$\heta(\htau_t(b))=\lambda^{-t/2}Q^{2it}\heta(b)$.
\end{pf}\vs

Using now Theorem \reef{newthm} we get

\begin{Thm}
\label{th7.2}
The $C^*$-bialgebra $(\hA,\hdelta)$ is a weighted Hopf $C^*$-algebra with
antipode $\hkappa$, Haar weight $\hh$ and scaling constant $\widehat{\lambda}$
equal to the inverse of the corresponding constant $\lambda$ for $(A,\delta)$.
\end{Thm}

By definition $(\hA,\hdelta)$ is the dual of the weighted Hopf 
$C^*$-algebra $(A,\delta)$. It is easy to see that the
passage $W\longmapsto\hW$ is involutive:
\[
\widehat{\hW}=\Sigma\hW^*\Sigma=\Sigma\left(\Sigma W^*\Sigma\right)^*\Sigma=W.
\]
Therefore we have the following

\begin{Thm}[Duality]
\label{th7.d}
The dual of the dual of a weighted Hopf $C^*$-algebra is isomorphic 
to the original weighted Hopf $C^*$-algebra.
\end{Thm}

Developing the theory of weighted Hopf $C^*$-algebras we constructed 
a number of objects related to
a given weighted Hopf $C^*$-algebra $(A,\delta)$. We list some of 
them in order of appearance:
\[
\kappa ,\; h , \;\tau , \;R , \;\lambda , \;\Hil , \;\eta , \;Q, \;\Delta ,
\;J, \;W, \;\hA , \;\hdelta , \;\hR , \;\htau , \;\hDelta ,
\;\hJ \;\mbox{ and } \;\heta .
\]
Clearly the corresponding objects related to the dual
weighted Hopf $C^*$-algebra $(\hA,\hdelta)$ equals:
\[
\hkappa , \;\hh , \;\htau ,
\;\hR , \;\hlambda=\lambda^{-1} , \;\Hil , \;\heta , \;Q, \;\hDelta ,
\;\hJ , \;\hW=\Sigma W^*\Sigma , \;A, \;\delta , \;R, \;\tau , 
\;\Delta , \;J \;\mbox{ and } \;\eta .
\]
The objects dual to $\rho$ and $\gamma$ will be denoted by $\hrho$ and
$\hgamma$, respectively.
Many of the formulae appearing in the
paper form dual pairs. The examples are: (\rf{character},\rf{deltahat}),
(\rf{hJQkom},\rf{Q}), (\rf{hRT1},\rf{th5.11a}), etc. The formulae 
\rf{th5.11b}, \rf{pr3.64a}, \rf{april0},
  etc. are selfdual.\vs

\begin{Prop}
Let $\hDelta_{\rm rel}$ be the relative modular operator determined 
by the weight $\hh$ and the weight $\hh\circ\hR$ $(\hDelta_{\rm rel}$ 
is the object
dual to $\Delta_{\rm rel})$. Then
\begin{equation}
\label{april0}
\hDelta_{\mathrm{rel}}\otimes\Delta=W(\hDelta\otimes\Delta_{\mathrm{rel}})W^*.
\end{equation}
\end{Prop}
\begin{pf}
We start with formula \rf{deltarho}. Applying to the both sides 
${\Ad}_{J\otimes\hJ}$ and using
\rf{th5.11b} we obtain
\[
W^*(J\rho J\otimes 1)W=J\rho J\otimes\hJ\rho\hJ.
\]
Combining Proposition \reef{pr4.5} with \rf{th5.11a} we see that 
$\hJ\rho\hJ=\rho^{-1}$. Therefore
\begin{equation}
\label{april1}
W^*(J\rho J\otimes 1)W=J\rho J\otimes\rho^{-1}.
\end{equation}
We shall also use the dual formula
\[
\hW^*(\hJ\hrho\hJ\otimes 1)\hW=\hJ\hrho\hJ\otimes\hrho^{-1}.
\]
By easy computation it reduces to
\begin{equation}
\label{april2}
W^*(\hrho^{-1}\otimes\hJ\hrho\hJ)W=1\otimes\hJ\hrho\hJ.
\end{equation}
Let $t\in\real$. Rising both sides of \rf{april1} and \rf{april2} to 
the power $-it$ and multiplying side by side the two equations we 
obtain
\[
W^*(J\rho^{it}J\hrho^{it}\otimes\hJ\hrho^{it}\hJ)W=J\rho^{it}J\otimes\rho^{it}\hJ\hrho^{it}\hJ.
\]
Remembering that $W$ commutes with $Q\otimes Q$ we get
\begin{equation}
\label{april3}
W^*(J\rho^{it}J\hrho^{it}Q^{2it}\otimes\hJ\hrho^{it}\hJ 
Q^{2it})W=J\rho^{it}JQ^{2it}\otimes\rho^{it}\hJ\hrho^{it}\hJ Q^{2it}.
\end{equation}

Formula \rf{hDelta} shows that $J\rho^{it}JQ^{2it}=\hDelta^{it}$.
By duality $\hJ\hrho^{it}\hJ Q^{2it}=\Delta^{it}$. Consequently,
by \rf{sekunda} and \rf{gammalambda} we have
$\rho^{it}\hJ\hrho^{it}\hJ
Q^{2it}=\rho^{it}\Delta^{it}=\gamma^{-it^2/2}\Delta_{\rm
rel}^{it}=\lambda^{it^2/2}\Delta_{\rm rel}^{it}$. By duality
$\hrho^{it}J\rho^{it}J Q^{2it}=\lambda^{-it^2/2}\hDelta_{\rm
rel}^{it}$. Replacing $t$ by $-t$ and taking hermitian conjugation
we obtain $J\rho^{it}J
Q^{2it}\hrho^{it}=\lambda^{it^2/2}\hDelta_{\rm rel}^{it}$. We know
(cf. \rf{opq56}) that $\rho$ is $\tau_t$-invariant. Therefore it
commutes with $Q$. So does $\hrho$ (duality!). Therefore
$J\rho^{it}J\hrho^{it} Q^{2it}=\lambda^{it^2/2}\hDelta_{\rm
rel}^{it}$. Inserting these data into \rf{april3} we obtain:
\[
W^*(\hDelta_{\rm 
rel}^{it}\otimes\Delta^{it})W=\hDelta^{it}\otimes\Delta_{\rm rel}^{it}
\]
and \rf{april0} follows.
\end{pf}

\vs

We end this Section with an explicit formula for the relative 
Tomita-Takesaki operator $S_{\rm rel}$. By definition
$S_{\rm rel}$ is a closed conjugate linear operator acting on $\Hil$ such that
\[
S_{\rm rel}b^*\eta(a)=a^*\eta_L(b)
\]
for any $a\in\Dom(\eta'')$ and $b\in\Dom(\eta_L'')$. The self-adjoint 
part of the polar decomposition of $S_{\rm rel}$, denoted by
$\Delta_{\rm rel}^{1/2}$, has been already investigated in Section 
\reef{sec4}. Let $J_{\rm rel}$ be the antiunitary part of
this decomposition:
\[
S_{\rm rel}=J_{\rm rel}\Delta_{\rm rel}^{1/2}.
\]
We shall prove that $J_{\rm rel}=\lambda^{-i/8}J$. Let $a\in\Dom (\eta'')$ and
$b\in\Dom (\eta_L'')\cap\Dom(\eta'')$. Then we have
\[
\begin{array}{r@{\;}c@{\;}l}
JJ^*_{\rm rel}a^*\eta_L''(b) &=& JJ^*_{\rm rel}S_{\rm rel}b^*\eta''(a)
  = J\Delta^{1/2}_{\rm rel}b^*\eta''(a) \\
  &=& \gamma^{i/8}J\rho^{1/2}\Delta^{1/2}b^*\eta''(a)
  = \gamma^{i/8}J\rho^{1/2}Ja^*\eta''(b)  \Vs{5} \\
  &=& \gamma^{i/8}\lambda^{i/4}\eta_L''(a^*b)
  = \lambda^{i/8}a^*\eta_L''(b), \Vs{5}
\end{array}
\]
where we used \rf{prima} in the third step, \rf{betabeta} in the 
fifth step and \rf{gammalambda} in the last step.

Now we notice that due to the last statement of Proposition \reef{removeJ},
$a^*\eta_L''(b)$ runs over a dense subset of $\Hil$. Therefore
$JJ^*_{\rm rel}=\lambda^{i/8}1$, $J_{\rm rel}=\lambda^{-i/8}J$. In this way we
showed that
\begin{equation}
\label{klopoty}
S_{\rm rel}=\lambda^{-i/8}J\Delta_{\rm rel}^{1/2}.
\end{equation}

\section{Quantum Codouble}
\label{sek8}

In \cite{pw} P. Podle\'s and the third author of the present paper 
have introduced the so called `double group
construction'. It was the construction dual to the quantum double of 
Drinfeld \cite{Drin}. In this section we describe
this construction within the framework of weighted Hopf 
$C^*$-algebras. The weighted Hopf C$^*$-algebra
$(B,\delta^B)$ arising with this construction is built of a weighted 
Hopf $C^*$-algebra $(A,\delta)$ and its
dual $(\hA,\hdelta)$. We shall call $(B,\delta^B)$ the quantum 
codouble of $(A,\delta)$. The theory presented in \cite{pw}
deals with the case of compact $(\hA,\hdelta)$ (and discrete $(A,\delta)$).\vs

A weighted Hopf $C^*$-algebra is said to be {\it unimodular}, if the 
right Haar weight is left
invariant. Then $h\comp R=h$. \vs

Let $(A,\delta)$ be a weighted Hopf $C^*$-algebra,
$({\hA},{\hdelta})$ be its dual and $W\in M(\hA\otimes A)$ be the 
Kac-Takesaki operator related to $(A,\delta)$.  We
shall consider the twisted flip $\sigma_W$ acting from $\hA\otimes A$ 
to $A\otimes\hA$ defined by
\[
\sigma_W=\sigma\comp\Ad_{W^*}=\Ad_{\hat W}\comp\sigma,
\]
where $\sigma$ is the flip map: $\sigma(b\otimes a)=a\otimes b$ for 
any $a\in A$ and $b\in\hA$. Then
$\sigma_W\in\Mor(\hA\otimes A,A\otimes\hA)$.\vs

The quantum codouble of $(A,\delta)$ is by definition the pair 
$(B,\delta^B)$, where
\[
\begin{array}{r@{\;}c@{\;}l}
B&=&\hA\otimes A\\ \Vs{5}
\delta^B&=&(\id\otimes\sigma_W\otimes\id)\comp(\hdelta\otimes\delta). \Vs{5}
\end{array}
\]

In the second formula
$\hdelta\otimes\delta\in\Mor(\hA\otimes A,\hA\otimes\hA\otimes A\otimes A)=
\Mor(B,\hA\otimes\hA\otimes A\otimes A)$ and 
$\id\otimes\sigma_W\otimes\id\in\Mor(\hA\otimes\hA\otimes A\otimes A,
\hA\otimes A\otimes\hA\otimes A)=\Mor(\hA\otimes\hA\otimes A\otimes 
A,B\otimes B)$. Therefore $\delta^B\in\Mor(B,
B\otimes B)$. We start with the following
\begin{Prop}
\label{lm8.3}
The pair $(B,\delta^B)$ is a proper C$^*$-bialgebra with the 
cancellation property.
\end{Prop}
\begin{pf}
Clearly $B$ is a separable C$^*$-algebra. Using formulae 
\rf{character} and \rf{deltahat} one can
easily verify that
\begin{equation}
\label{warkocz}
\begin{array}{r@{\;}c@{\;}l}
(\delta\otimes\id)\sigma_W&=&(\id\otimes\sigma_W)(\sigma_W\otimes\id)(\id\otimes\delta),\\
(\id\otimes\hdelta)\sigma_W&=&(\sigma_W\otimes\id)(\id\otimes\sigma_W)(\hdelta\otimes\id).\Vs{5}
\end{array}
\end{equation}
Now the coassociativity of $\delta^B$ follows from a simple 
computation. For details see
\cite[the proof of Theorem 4.1]{pw}.\vs

For any subset $E\subset M(B\otimes B)$ we denote by $\CLS{E}$ the 
norm-closed linear span of $E$.
To prove the properness and the cancellation property we have to show 
that $\CLS{\delta^B(B)(1\otimes B)}=
\CLS{(B\otimes 1)\delta^B(B)}=B\otimes B$.

We know that the set $\delta(A)(1\otimes A)$ is a linearly dense 
subset of $A\otimes A$. Therefore
\[
\begin{array}{r@{\;}c@{\;}l}
\CLS{\delta^B(B)(1\otimes 
B)}&=&\Sigma_{23}\CLS{W_{23}^*\left(\hdelta(\hA)\otimes\delta(A)\right)W_{23}
(1\otimes\hA\otimes 1\otimes A)}\Sigma_{23}\\ \Vs{5}
&=&\Sigma_{23}\CLS{W_{23}^*\left(\hdelta(\hA)\otimes\delta(A)(1\otimes 
A)\right)W_{23}
(1\otimes\hA\otimes 1\otimes 1)}\Sigma_{23}\\ \Vs{5}
&=&\Sigma_{23}\CLS{W_{23}^*\left(\hdelta(\hA)\otimes A\otimes A\right)W_{23}
(1\otimes\hA\otimes 1\otimes 1)}\Sigma_{23}\\ \Vs{5}
&=&\Sigma_{23}\CLS{W_{23}^*\hW_{12}\left(\hA\otimes 1\otimes A\otimes 
A\right)\hW_{12}^*W_{23}
(1\otimes\hA\otimes 1\otimes 1)}\Sigma_{23},
\end{array}
\]
where in the last step we used the formula $\hdelta(b)=\hW(b\otimes 
1)\hW^*$. The pentagonal
equation implies that $\hW_{12}^*W_{23}=W_{13}\hW_{12}^*W_{13}^*$. 
Therefore using in the second step the formula
$(\hA\otimes A)W=(\hA\otimes A)$ we obtain
\[
\begin{array}{r@{\;}c@{\;}l}
\CLS{\delta^B(B)(1\otimes 
B)}&=&\Sigma_{23}\CLS{W_{23}^*\hW_{12}\left(\hA\otimes 1\otimes 
A\otimes
A\right)W_{13}\hW_{12}^*W_{13}^* (1\otimes\hA\otimes 1\otimes 
1)}\Sigma_{23}\\ \Vs{5}
&=&\Sigma_{23}\CLS{W_{23}^*\hW_{12}\left(\hA\otimes 1\otimes A\otimes
A\right)\hW_{12}^*W_{13}^* (1\otimes\hA\otimes 1\otimes 1)}\Sigma_{23}\\ \Vs{5}
&=&\Sigma_{23}W_{23}^*\CLS{(\hdelta(\hA)\otimes A\otimes
A) (1\otimes\hA\otimes 1\otimes 1)}W_{13}^*\Sigma_{23}.
\end{array}
\]
We know that the set $\hdelta(\hA)(1\otimes\hA)$ is a linearly dense 
subset of $\hA\otimes\hA$. Therefore
\[
\begin{array}{r@{\;}c@{\;}l}
\CLS{\delta^B(B)(1\otimes B)}
&=&\Sigma_{23}W_{23}^*(\hA\otimes\hA\otimes
A\otimes A)W_{13}^*\Sigma_{23}\\ \Vs{5}
&=&\Sigma_{23}(\hA\otimes\hA\otimes
A\otimes A)\Sigma_{23}\\ \Vs{5}
&=&\hA\otimes A\otimes
\hA\otimes A=B\otimes B\Vs{5}.
\end{array}
\]
In the similar way one can show that $\CLS{(B\otimes 
1)\delta^B(B)}=B\otimes B$.
\end{pf}\vs

Let $\eta$, $\heta$ and $\eta_L$ be the GNS maps considered in 
previous sections. We shall also
use the GNS map $\heta_L$ defined on $\hA$, introduced by the formula
\begin{equation}
\label{dualpha}
\heta_L(b)=J\heta(\hR(b^*))
\end{equation}
dual to \rf{alpha}. By definition 
$\Dom(\heta_L)=\set{b\in\hA}{\hR(b^*)\in\Dom(\heta)}$. As we know, 
$\eta$ and $\eta_L$
are GNS maps related to the right and left Haar weights $h$ and 
$h\comp R$ on $(A,\delta)$.  In same way, $\heta$ and
$\heta_L$ are GNS maps related to the right and left Haar weights 
$\hh$ and $\hh\comp\hR$ on $(\hA,\hdelta)$. We shall use the tensor 
product of GNS mappings
(see Appendix \reef{Tensor}). Let
\begin{equation}
\label{GNSB}
\begin{array}{r@{\;}c@{\;}l}
\eta^B&=&\heta_L\otimes\eta,\\
\eta^B_d&=&\heta\otimes\eta_L.\Vs{4}
\end{array}
\end{equation}
Then $\eta^B$ and $\eta^B_d$ are closed densely defined GNS maps on 
$B=\hA\otimes A$ with values in $\Hil\otimes\Hil$.
  Let $U\in\skrL(\Hil\otimes\Hil)$ be  the unitary operator introduced by
\begin{equation}
\label{U}
U=W(\hJ\otimes J)W(\hJ\otimes J).
\end{equation}
We recall that $W\in M(\hA\otimes A)=M(B)$.

\begin{Prop}
\label{pro8.3}
We have$\,:$ $\Dom(\eta^B_d)=\Ad_{W^*}(\Dom(\eta^B))$ and
\begin{equation}
\label{for8.3}
U^*\eta^B(c)=\lambda^{i/4}\eta^B_d(\Ad_{W^*}(c))
\end{equation}
for any $c\in\Dom(\eta^B)$.
\end{Prop}
\begin{pf}
Recall that $S_{\mathrm{rel}}\eta(a)=\eta_L(a^*)$ for
$a\in\Dom(\eta)\cap\Dom(\eta_L)^*$ and (cf. \rf{klopoty})
$S_{\mathrm{rel}}=\lambda^{-i/8}J\Delta_{\mathrm{rel}}^{1/2}$. By
duality $\hS_{\mathrm{rel}}\heta(a)=\heta_L(a^*)$ for
$a\in\Dom(\heta)\cap\Dom(\heta_L)^*$ and
$\hS_{\mathrm{rel}}=\lambda^{i/8}\hJ\hDelta_{\mathrm{rel}}^{1/2}$. To
make the notation shorter we set $J_B=\hJ\otimes J$. For any
$a\in\Dom(\eta)\cap\Dom(\eta_L)^*$ and
$b\in\Dom(\heta)\cap\Dom(\heta)^*$ we have
\[
\begin{array}{r@{\;}c@{\;}l}
U\eta^B_d(b^*\otimes a^*) &=&
WJ_BWJ_B(\heta(b^*)\otimes\eta_L(a^*))
  = WJ_BWJ_B(\hS\otimes S_{\rm rel})(\heta(b)\otimes\eta(a))\Vs{5}\\
  &=& 
\lambda^{-i/8}WJ_BW(\hDelta\otimes\Delta_{\mathrm{rel}})^{1/2}(\heta(b)\otimes\eta(a))\Vs{5}\\
&=&\lambda^{-i/8}WJ_B(\hDelta_{\mathrm{rel}}\otimes\Delta)^{1/2}W(\heta(b)\otimes\eta(a)) 
\Vs{5}\\
&=&\lambda^{-i/8}WJ_B(\hDelta_{\mathrm{rel}}\otimes\Delta)^{1/2}(\heta\otimes\eta)(W(b\otimes
a))\Vs{5}\\
&=&\lambda^{-i/4}W(\widehat{S}_{\rm rel}\otimes 
S)(\heta\otimes\eta)(W(b\otimes a))\Vs{5}
,
\end{array}
\]
where in the fourth step we used relation \rf{april0}. Since 
$W(b\otimes a)\in\Dom(\heta\otimes\eta)$ and
$(\heta\otimes\eta)(W(b\otimes
a))\in\Dom(\hS_{\mathrm{rel}}\otimes S)$, it follows (cf. Theorem 
\reef{19.08.2002c}) that $(b\otimes 
a)^*W^*\in\Dom(\heta_L\otimes\eta)$.  Hence we see that
\[
\begin{array}{r@{\;}c@{\;}l}
U\eta^B_d(b^*\otimes a^*) &=&
\lambda^{-i/4}W(\widehat{S}_{\rm rel}\otimes 
S)(\heta\otimes\eta)(W(b\otimes a))\Vs{5}\\
&=&\lambda^{-i/4}W(\heta_L\otimes \eta)((b^*\otimes a^*)W^*)\Vs{5}\\ \Vs{5}
&=&\lambda^{-i/4}(\heta_L\otimes \eta)(W(b^*\otimes a^*)W^*)
=\lambda^{-i/4}\eta^B(\Ad_W(b^*\otimes a^*)) .
\end{array}
\]
Since the algebraic tensor product of $\Dom(\heta)\cap\Dom(\heta)^*$ and
$\Dom(\eta_L)\cap\Dom(\eta)^*$ is a core for the GNS mapping 
$\eta_d^B$, it follows from the above argument that 
$\Ad_W(c)\in\Dom(\eta^B)$ and
\begin{equation}
\label{tetsuya}
U\eta_d^B(c)=\lambda^{-i/4}\eta^B(\Ad_W(c))
\end{equation}
for any $c\in\Dom(\eta_d^B)$.  Hence we see that 
$\Ad_W(\Dom(\eta_d^B))\subset\Dom(\eta^B)$. Passing to the dual 
statement we obtain
$\Ad_{W^*}(\Dom(\eta^B))\subset\Dom(\eta_d^B)$. Combining the two inclusions we
obtain the equality $\Ad_{W^*}(\Dom(\eta^B))=\Dom(\eta_d^B)$. To end the proof
we notice that \rf{tetsuya} is equivalent to \rf{for8.3}.

\end{pf}\vs

Now we shall introduce the unitary antipode, the scaling group and 
the Haar weight related to $(B,\delta^B)$:

\begin{equation}
\label{danecodouble}
\left.\begin{array}{r@{\;}c@{\;}l}
R^B&=&({\hR}\otimes R)\comp{\Ad}_{W^*},\\ \Vs{5}
\tau_t^B&=&{\htau}_t\otimes\tau_t,\\ \Vs{5}
h^B&=&\hh\comp\hR\otimes h,\Vs{5}
\end{array}\right\}
\end{equation}
where $h^B=\hh\comp\hR\otimes h$ is the weight on $B$ associated with the GNS
mapping $\eta^B=\heta_L\otimes\eta$. Using the properties
\[
({\hR}\otimes R)(W)=W, \quad ({\hR}\otimes 
R)\comp{\Ad}_{W^*}={\Ad}_W\comp({\hR}\otimes R),
\quad({\htau}_t\otimes\tau_t)(W)=W, \quad\hlambda\lambda=1
\]
one can easily verify the following:

\begin{Prop}
\label{lm8.5}
$R^B$ is an involutive antiautomorphism of $B$, 
$\seq{\tau_t^B}{t\in\real}$ is a strongly continuous
one parameter group of automorphisms of $B$ commuting with $R^B$ and $h^B$
is a locally finite, strictly faithful, lower semicontinuous weight on $B$ such
that $h^B\comp\tau^B_t=h^B$ for all $t\in\real$. Moreover
\begin{equation}
\label{unimodular}
h^B\comp R^B=h^B \;.
\end{equation}
\end{Prop}

\begin{pf}
Only the last formula needs a justification. Let $h_d^B$ be the weight on $B$
related to the GNS map $\eta^B_d$. By definition
\[
h^B_d=\hh\otimes h\comp R=(\hh\comp\hR\otimes h)\comp(\hR\otimes 
R)=h^B\comp(\hR\otimes R).
\]
Therefore
\[
h^B\comp R^B=h^B_d\comp\Ad_{W^*}=h^B,
\]
where the last equality follows immediately from \rf{for8.3}.
\end{pf}\vs

Our main aim in this section is to show that $(B,\delta^B)$ is a 
weighted Hopf C$^*$-algebra. To this end we shall use
Theorem \reef{newthm} with $A,\delta,R,\tau$ and $h$ replaced by 
$B,\delta^B,R^B,\tau^B$ and $h^B$.
Then the role of $\Hil$ and $\eta$ is played by 
$\Hil^B=\Hil\otimes\Hil$ and $\eta^B$ and the unitary $W$ (denoted in
the present case by $W^B$) acts on 
$\Hil^B\otimes\Hil^B=\Hil\otimes\Hil\otimes\Hil\otimes\Hil$. We shall 
prove that
\begin{equation}
\label{KTOB}
W^B=U_{12}\hW_{13}U_{12}^*W_{24}.
\end{equation}
satisfies all the assumptions of Theorem \reef{newthm}. Remembering 
that $W\in M(\hA\otimes A)\subset
M(\skrK(\Hil)\otimes A)$ and $\hW\in M(A\otimes\hA)\subset 
M(\skrK(\Hil)\otimes\hA)$ we see that
$W^B\in M(\skrK(\Hil)\otimes\skrK(\Hil)\otimes\hA\otimes A)$ and
\begin{equation}
\label{multi}
W^B\in M(\skrK(\Hil^B)\otimes B).
\end{equation}

In the following propositions we shall use a special notation. For 
any $\varphi\in A^*$ and any $a\in A$ we set:
$\varphi_*(a)=\varphi*a$. Then $\varphi_*$ is a linear mapping acting 
on $A$. Similarly for any $\psi\in\hA^*$ and any $b\in\hA$
we set: $\psi_*(b)=\psi*b$. Then $\psi_*$ is a linear mapping acting 
on $\hA$. Clearly
\[
\begin{array}{r@{\;}c@{\;}l}
\varphi_*&=&(\id\otimes\varphi)\comp\delta,\\
\psi_*&=&(\id\otimes\psi)\comp\hdelta.\Vs{4}
\end{array}
\]
One can easily verify that $\varphi_*$ and $\psi_*$ are completely 
bounded. Consequently $\psi_*\otimes\id$ and
$\id\otimes\varphi_*$ are bounded linear mappings acting on $B=\hA\otimes A$.
\begin{Prop}
For any $\psi\in\hA^*$, $\varphi\in A^*$ we have$\,:$
\begin{equation}
\label{forPr8.4}
\begin{array}{r@{\;}c@{\;}l}
(\id\otimes\id\otimes\psi\otimes\varphi)W^B&=&U\left[(\id\otimes\psi)\hW\otimes 
1\right]U^*
\left[1\otimes(\Vs{4}\id\otimes\varphi)W\right],\\ \Vs{5}
(\id\otimes\id\otimes\psi\otimes\varphi)\comp\delta^B&=&\Ad_W\comp
(\psi_*\otimes\id)\comp\Ad_{W^*}\comp (\id\otimes\varphi_*).
\end{array}
\end{equation}
\end{Prop}
\begin{pf}
The first formula follows immediately from \rf{KTOB}. Using the 
second formula of \rf{warkocz} we compute:
\[
\begin{array}{r@{\;}c@{\;}l}
(\sigma_W\otimes\id\otimes\id)\comp\delta^B
&=&(\sigma_W\otimes\id\otimes\id)\comp
(\id\otimes\sigma_W\otimes\id)\comp(\hdelta\otimes\delta)\\ \Vs{5}
&=&(\id\otimes\hdelta\otimes\id)\comp(\sigma_W\otimes\id)\comp(\id\otimes\delta).
\end{array}
\]
Therefore
\[
\delta^B=(\sigma_W^{-1}\otimes\id\otimes\id)\comp
(\id\otimes\hdelta\otimes\id)\comp(\sigma_W\otimes\id)\comp(\id\otimes\delta)
\]
and
\[
(\id\otimes\id\otimes\psi\otimes\varphi)\comp\delta^B=\sigma_W^{-1}\comp
(\id\otimes\psi_*)\comp\sigma_W\comp(\id\otimes\varphi_*).
\]
Remembering that $\sigma_W=\sigma\comp\Ad_{W^*}$ and using the 
relation $\sigma\comp(\id\otimes\psi_*)\comp\sigma=
\psi_*\otimes\id$ we get the second formula of \rf{forPr8.4}.
\end{pf}

\begin{Prop}

\label{pr77}\Vs{2}

$1.$ For any $\varphi\in A_*$ and any $c\in\Dom(\eta^B)$ we have$\,:$ 
$(\id\otimes\varphi_*)c\in\Dom(\eta^B)$ and
\[
\eta^B((\id\otimes\varphi_*)c)=\left[1\otimes(\Vs{4}\id\otimes\varphi)W\right]\eta^B(c) 
\;.
\]

$2.$ For any $\psi\in\hA_*$ and any $c\in\Dom(\eta^B_d)$ we have$\,:$ 
$(\psi_*\otimes\id)c\in\Dom(\eta^B_d)$ and
\[
\eta^B_d((\psi_*\otimes\id)c)=\left[(\id\otimes\psi)\hW\otimes 
1\right]\eta^B_d(c) \;.
\]

$3.$ For any $\mu\in B_*$ and any $c\in\Dom(\eta^B)$ we have$\,:$ 
$(\id\otimes\mu)\delta^B(c)\in\Dom(\eta^B)$ and
\[
\eta^B\left((\id\otimes\mu)\delta^B(c)\right)=\left[\Vs{4}(\id\otimes\mu)W^B\right]\eta^B(c).
\]

\noindent In the last point, $\id$ denotes the identity map acting on 
operators in $\Hil\otimes\Hil$.
\end{Prop}
\begin{pf} Ad 1. For any $\sigma$-weakly continuous functional $\varphi$, the
mapping $\id\otimes\varphi$ maps strongly convergent sequences into
strongly convergent sequences. Formula $\delta(a)=W(1\otimes a)W^*$ 
(cf. \rf{coproduct}) shows that the comultiplication
$\delta$  has the same property. Therefore $\varphi_*$ maps strongly 
convergent sequences into strongly convergent sequences.
  Therefore we may assume that $c=b\otimes a$, where 
$b\in\Dom(\heta_L)$ and $a\in\Dom(\eta)$ (by definition
$\Dom(\heta_L)\otimes_{\rm alg}\Dom(\eta)$ is a core of 
$\eta^B=\heta_L\otimes\eta$). In this case our relation takes the form
\[
\heta_L(b)\otimes\eta(\varphi*a)=\heta_L(b)\otimes\left[\id\otimes\varphi)W\right]\eta(a),
\]
which immediately reduces to \rf{Ftran} (cf. \rf{lslice}).\vs

Ad 2. In this point we may assume that $c=b\otimes a$, where 
$b\in\Dom(\heta)$ and $a\in\Dom(\eta_L)$. In this case our
relation takes the form
\[
\heta\left((\id\otimes\psi)\hdelta(b)\right)\otimes\eta_L(a)
=\left[(\id\otimes\psi)\hW\right]\heta(b)\otimes\eta_L(a),
\]
which immediately reduces to \rf{zalozeniex}.\vs

Ad 3.
By Remark \reef{koniec3} we may assume that $\mu=\psi\otimes\varphi$, 
where $\psi\in\hA^*$ and $\varphi\in A^*$.
(the reader should notice that $\hA_*\otimes_{\rm alg}A_*$ is norm 
dense in $B_*$). Let
$c\in\Dom(\eta^B)$. Using Statement 1, Proposition \reef{pro8.3}, 
Statement 2 and again Proposition \reef{pro8.3} we
obtain: $(\id\otimes\varphi_*)c\in\Dom(\eta^B)$, 
$\Ad_{W^*}(\id\otimes\varphi_*)c\in\Dom(\eta^B_d)$,
$(\psi_*\otimes\id)\Ad_{W^*}(\id\otimes\varphi_*)c\in\Dom(\eta^B_d)$ and
\[
\Ad_W(\psi_*\otimes\id)\Ad_{W^*}(\id\otimes\varphi_*)c\in\Dom(\eta^B) \;.
\]
According to \rf{forPr8.4} the last expression coincides with 
$(\id\otimes\mu)\delta^B(c)$. Taking into account \rf{forPr8.4} and 
\rf{for8.3} we have:
\[
\begin{array}{r@{\;}c@{\;}l}
\left[(\id\otimes\mu)W^B\right]\eta^B(c)&=&U\left[(\id\otimes\psi)\hW\otimes 
1\right]U^*
\left[1\otimes(\Vs{4}\id\otimes\varphi)W\right]\eta^B(c)\\
&=&U\left[(\id\otimes\psi)\hW\otimes 1\right]U^*
\eta^B((\id\otimes\varphi_*)c)\Vs{6}\\
&=&\lambda^{i/4}U\left[(\id\otimes\psi)\hW\otimes 
1\right]\eta^B_d(\Ad_{W^*}(\id\otimes\varphi_*)c)\Vs{6}\\
&=&\lambda^{i/4}U\eta^B_d((\psi_*\otimes\id)\Ad_{W^*}(\id\otimes\varphi_*)c)\Vs{5}\\
&=&\eta^B(\Ad_W(\psi_*\otimes\id)\Ad_{W^*}(\id\otimes\varphi_*)c)\Vs{5}\\&=&
\eta^B((\id\otimes\mu)\delta^B(c)).\Vs{5}
\end{array}
\]
\end{pf}\vs

Let $\kappa^B=R^B\comp\tau^B_{i/2}$. Then 
$\kappa^B=(\hkappa\otimes\kappa)\comp\Ad_{W^*}$. Let $\psi\in \hA^*$ 
and $\varphi\in
A^*$ be linear functionals such that $\psi\comp\hkappa\in \hA^*$ and 
$\varphi\comp\kappa\in A^*$. For any $c\in B$ we set:
\begin{equation}
\label{mu}
\mu(c)=(\psi\otimes\varphi)(W^*cW).
\end{equation}
Then $\mu\in B^*$ and one can easily verify that
\[
\mu\comp\kappa^B=\psi\comp\hkappa\otimes\varphi\comp\kappa.
\]
In particular $\mu\comp\kappa^B\in B^*$.\vs

\begin{Prop}
\label{pr88}
For any $\mu\in B^*$ such that $\mu\comp\kappa^B\in B^*$ we have$\,:$
\begin{equation}
\label{zalozenie1y}
(\id\otimes\mu\comp\kappa^B)W^B=(\id\otimes\mu)({W^B}^*).
\end{equation}
\end{Prop}
\begin{pf}
In the following computation we use leg numbering notation for 
operators acting on $\Hil^{\otimes 4}=
\Hil\otimes\Hil\otimes\Hil\otimes\Hil$. Our first aim is formula 
\rf{UUUU}. We start with the pentagonal equation
$W^*_{14}W^*_{13}=W_{34}W^*_{13}W^*_{34}$. Multiplying the both sides 
from the right by $W^*_{12}$ we get
\[
W^*_{14}W^*_{13}W^*_{12}=W_{34}W^*_{13}W^*_{12}W^*_{34}.
\]
Combining this formula with the pentagonal equation 
$W^*_{13}W^*_{12}=W_{23}W^*_{12}W^*_{23}$ we obtain:
\begin{equation}
\label{UU}
W^*_{14}W_{23}W^*_{12}W^*_{23}=W_{34}W_{23}W^*_{12}W^*_{23}W^*_{34}.
\end{equation}
For any operator $V$ acting on $\Hil^{\otimes 4}$ we set:
\[
{\mathcal Ad}(V)=(J\otimes\hJ\otimes J\otimes\hJ)V(J\otimes\hJ\otimes 
J\otimes\hJ).
\]
Clearly ${\mathcal Ad}$ is an conjugate linear multiplicative 
operation acting on $\skrL(\Hil^{\otimes 4})$. Using  \rf{th5.11b} we 
obtain immediately
the relations: ${\mathcal Ad}(W_{34})=W^*_{34}$, ${\mathcal 
Ad}(W^*_{14})=W_{14}$ and
${\mathcal Ad}(W^*_{12})=W_{12}$. Furthermore, taking into account 
definition \rf{U} we see that
${\mathcal Ad}(W_{23})=U_{23}W_{23}^*$. Applying ${\mathcal Ad}$ to 
the both sides of \rf{UU} we get:
\[
W_{14}U_{23}W_{23}^*W_{12}W_{23}U_{23}^*=W_{34}^*U_{23}W_{23}^*W_{12}W_{23}U_{23}^*W_{34}.
\]
Multiplying both sides by $W_{13}$ from the right and using pentagon equation
$W_{34}W_{13}=W_{13}W_{14}W_{34}$ we obtain
\begin{equation}
\label{UUU}
W_{14}U_{23}W_{23}^*W_{12}W_{23}U_{23}^*W_{13}=W_{34}^*U_{23}W_{23}^*W_{12}W_{23}U_{23}^*W_{13}W_{14}W_{34}.
\end{equation}
The reader should notice that the `third leg' of 
$W_{23}U_{23}^*=\{(\hJ\otimes J)W^*(\hJ\otimes J)\}_{23}$
`belongs' to the commutant $A'$ of $A$, whereas $W_{13}\in 
M(\hA\otimes 1\otimes A\otimes 1)$. Therefore $W_{23}U_{23}$
commutes with $W_{13}$ and
\[
W_{23}^*W_{12}W_{23}U_{23}^*W_{13}=W_{23}^*W_{12}W_{13}W_{23}U_{23}^*=W_{12}U_{23}^*.
\]
Therefore \rf{UUU} takes the form
\[
W_{14}U_{23}W_{12}U_{23}^*=W_{34}^*U_{23}W_{12}U_{23}^*W_{14}W_{34}.
\]
Now we perform the cyclic permutation of the first three $\Hil$ in 
$\Hil^{\otimes 4}=
\Hil\otimes\Hil\otimes\Hil\otimes\Hil$. It results with the following 
replacement of the `leg numbers':
$1\rightarrow 3\rightarrow 2\rightarrow 1$. Then $W_{12}$ goes into 
$W_{31}=\hW^*_{13}$ and our formula takes the form
\[
W_{34}U_{12}\hW^*_{13}U_{12}^*=W_{24}^*U_{12}\hW^*_{13}U_{12}^*W_{34}W_{24}.
\]
Rearranging this formula we obtain
\begin{equation}
\label{UUUU}
{W^B}^*=W_{24}^*U_{12}\hW^*_{13}U_{12}^*=W_{34}U_{12}\hW^*_{13}U_{12}^*W_{24}^*W_{34}^*.
\end{equation}

Let $\mu$ be the functional on $B$ introduced by \rf{mu}. Applying 
$\id\otimes\mu$ to the both sides of \rf{UUUU} and using in
the third step formulae
\rf{SLWx} and
\rf{zalozenie1x} we obtain:
\[
\begin{array}{r@{\;}c@{\;}l}
(\id\otimes\mu){W^B}^*&=&(\id\otimes\id\otimes\psi\otimes\varphi)\left(U_{12}\hW^*_{13}U_{12}^*W_{24}^*\right)\\ 
\Vs{5}
&=&U\left[(\id\otimes\psi)(\hW^*)\otimes 
1\right]U^*\left[\Vs{4}1\otimes(\id\otimes\varphi)(W^*)\right]\\ 
\Vs{5}
&=&U\left[(\id\otimes\psi\comp\hkappa)\hW\otimes 
1\right]U^*\left[\Vs{4}1\otimes(\id\otimes\varphi\comp\kappa)W\right]\\
&=&(\id\otimes\id\otimes\psi\comp\hkappa\otimes\varphi\comp\kappa)\left(U_{12}\hW_{13}U_{12}^*W_{24}\right)\Vs{5}\\
&=&(\id\otimes\mu\comp\kappa^B)(W^B).\Vs{5}
\end{array}
\]
In this way we showed that \rf{zalozenie1y} holds for any $\mu$ of 
the form \rf{mu}. Clearly the set of functionals of
the form \rf{mu} is linearly weakly$^*$ dense in $B^*$. Moreover for any
$t\in\real$ the functional $\mu\comp\tau^B_t$ is again of the form 
\rf{mu}. Indeed, remembering that $W$ commutes
with $Q\otimes Q$ we see that $\tau^B_t=\htau_t\otimes\tau_t$ 
commutes with $\Ad_{W^*}$ and
\[
\mu\comp\tau^B_t(c)=(\psi\comp\htau_t\otimes\varphi\comp\tau_t)(W^*cW).
\]
Now, using Remark \reef{koniec31} we obtain \rf{zalozenie1y} in full 
generality.
\end{pf}\vs

Let us summarize the results obtained so far. We have constructed the 
sequence of objects:
$B,\delta^B,R^B,\tau^B,h^B,\eta^B,W^B$. Taking into account 
Proposition \reef{lm8.3}, Proposition \reef{lm8.5}, formula
\rf{multi}, Statement 3 of Proposition \reef{pr77} and Proposition 
\reef{pr88} we see that these objects satisfy all the
assumptions of Theorem
\reef{newthm}. Making profit of this theorem we obtain

\begin{Thm}
The pair $(B,\delta^B)$ is a unimodular weighted Hopf C$^*$-algebra 
with $\lambda=1$. The unitary antipode, the scaling
group and the Haar weight related to $(B,\delta^B)$ are given by 
\rf{danecodouble}. The corresponding Kac-Takesaki
operator coincides with \rf{KTOB}.
\end{Thm}

The unimodularity of $(B,\delta^B)$ follows from \rf{unimodular}. It 
shows that the left and right Haar weight
coincide.
By the theory developed in Section \reef{sek3}, $W^B$ is a manageable 
multiplicative unitary. Inserting definition \rf{U}
into
\rf{KTOB} we obtain
\[
W^B=\left((\hJ\otimes J)W(\hJ\otimes J)\right)_{12}W_{12}\hW_{13}
W_{12}^*\left((\hJ\otimes J)W^*(\hJ\otimes J)\right)_{12}W_{24}.
\]
By pentagonal equation $W_{12}\hW_{13}W_{12}^*=\hW_{23}\hW_{13}$. Therefore
\[
W^B=\left((\hJ\otimes J)W(\hJ\otimes J)\right)_{12}\hW_{23}\hW_{13}
\left((\hJ\otimes J)W^*(\hJ\otimes J)\right)_{12}W_{24}.
\]
This formula admits further simplification:
\[
W^B=\hW_{23}\hW_{13}\left((J\otimes J)\hW^*(J\otimes J)\right)_{23}W_{24}.
\]
The verification of the last equality is left to the reader as an 
(un-easy) exercise.

\section{Appendix}

\renewcommand{\thesubsection}{\Alph{subsection}}
\renewcommand{\thesection}{\thesubsection}
\setcounter{equation}{0}\setcounter{Thm}{0}
\subsection{Multipliers and morphisms}
\label{A}

In this Appendix, we collect some technical statements concerning the 
multipliers and the morphisms
of the $C^*$-algebras.\vs

For any $\varphi\in A^*$ and $a,b\in A$ we define functionals 
$b\varphi$ and $\varphi a$ in $A^*$ by
\[
(b\varphi)(a)=\varphi(ab)=(\varphi a)(b).
\]
Then $A^*$ is a Banach $A$-bimodule.\vs

Let $A$ be a $C^*$-algebra and $\{e_{\alpha}\}$ an approximate
identity of $A$.  Let $\varphi\in A^*$.  Then $\varphi$ is written in 
the form $\varphi(a)=(x|\pi(a)y)$, where $\{\pi,\Hil\}$ is a 
representation of $A$ and $x,y\in\Hil$.  For any $a\in A$ we have
\[
|(e_{\alpha}\varphi-\varphi)(a)|=|\varphi(ae_{\alpha}-a)|
=|(x|\pi(ae_{\alpha}-a)y)|\leq \norm{x}\,\norm{a}\,\norm{\pi(e_{\alpha})y-y}
\]
and
\[
|(\varphi e_{\alpha}-\varphi)(a)|=|\varphi(e_{\alpha}a-a)|
=|(x|\pi(e_{\alpha}a-a)y)|\leq \norm{\pi(e_{\alpha})x-x}\,\norm{a}\,\norm{y}.
\]
Hence the sets $\{\varphi e_{\alpha}\}$ and $\{e_{\alpha}\varphi\}$ 
converges in norm to $\varphi$.
Thus we find that
\begin{eqnarray}
\overline{\{a\varphi: a\in A, \varphi\in A^*\}}
=\overline{\{\varphi a: a\in A, \varphi\in A^*\}}=A^*.
\label{DW}
\end{eqnarray}

Now we recall the Doran-Wichmann's factorization theorem:

\begin{Thm}[\cite{dw}]
Let $B$ be a Banach algebra and $E$ a left $($resp. right$)$ Banach 
$B$-module. If $B$ has a left $($resp. right$)$ approximate identity 
bounded by $K\geq 1$, then for any $x\in E$ and for any 
$\varepsilon>0$ there exist $b\in B$ and $y\in E$ such that
\[
x=by,\quad \norm{b}\leq K,\quad
y\in\overline{Bx}\ \ (resp.\ y\in\overline{xB}), \quad \norm{y-x}<\varepsilon.
\]
In particular, $\overline{BE}=BE$ $($resp. $\overline{EB}=EB)$.
\end{Thm}

Combining this theorem with (\reef{DW}), we find that

\begin{Prop}
\label{AP2}
Let $A$ be a $C^*$-algebra. Then $A^*=AA^*=A^*A$, i.e. for any 
$\varphi\in A^*$ there exist $\psi\in A^*$ and $a\in A$ such that 
$\varphi=a\psi$ $($or $\varphi=\psi a)$.
\end{Prop}

\begin{Lem}
Let $A$ and $B$ be $C^*$-algebras and $\varphi\in B^*$.
Then the slice mapping $\id_A\otimes\varphi :A\otimes B\rightarrow A$
has an extension
$(\id_A\otimes\varphi)\tilde{\,}$ from $M(A\otimes B)$ to $M(A)$ such that
$\norm{(\id_A\otimes\varphi)\tilde{\,}}=\norm{\varphi}$.
\end{Lem}

\begin{pf} Let $\seq{c_{\alpha}}{}$ be a net in $A\otimes B$ 
converging to zero with respect to the strict topology on $M(A\otimes 
B)$.
By Proposition \reef{AP2} $\varphi$ is of the form $b\psi$ for some 
$b\in A$ and $\psi\in A^*$.
Then, by the boundness of
$\psi$, for any $a\in A$,
$\{(\id_A\otimes \varphi)(c_{\alpha})\}a
=(\id_A\otimes\psi)(c_{\alpha}(a\otimes b))$ converges in norm
to zero.
Similarly, we obtain that the net
$\seq{a\{(\id_A\otimes\varphi
)(c_{\alpha})\}}{}$ converges in norm to zero. This proves that the
net $\seq{(\id_A\otimes\varphi
)(c_{\alpha})}{}$ converges to zero with respect to the strict
topology. Hence the slice mapping
$\id_A\otimes\varphi:A\otimes B\rightarrow A$ has an extension to the
mapping $(\id_A\otimes\varphi)\tilde{\,}$\ from $M(A\otimes B)$ to
$M(A)$.

It is known that the unit ball of a $C^*$-algebra is strictly dense 
in the unit ball of the multiplier algebra.  Therefore the norm of 
$(\id_A\otimes\varphi)\tilde{\,}$\ coincides with the norm of 
$\id\otimes\varphi$. It is clear that norms $\id\otimes\varphi$ and 
$\varphi$ are the same.
\end{pf}\vs

Throughout this paper, the notion of the morphism between the two 
$C^*$-algebras plays an important role. Here we recall the definition 
of the morphism.

\begin{Def}[\cite{w11}]\label{defmor}
Let $A$ and $B$ be $C^*$-algebras. Then $\pi$ is said to be in the 
set $\Mor (A,B)$ of morphisms from $A$ to $B$ if $\pi$ is a 
homomorphism from the $C^*$-algebra $A$ to the multiplier algebra 
$M(B)$ such that the subset $\pi (A)B$ of $B$ is norm dense in $B$.
\end{Def}

Each morphism $\pi$ of $A$ to $B$ can be extended uniquely to a 
morphism $\widetilde{\pi}$ of the multiplier algebra $M(A)$ to the 
multiplier algebra $M(B)$:
\[
\widetilde{\pi}(m)\sum_i\pi(a_i)x_i=\sum_i\pi(ma_i)x_i
\]
for $m\in M(A)$, $a_i\in A$ and $x_i$ in the carrier Hilbert space of $B$.
In this paper this extension of $\pi$ is denoted by the same letter 
(we omit tilde).\vs

The morphisms admit even further extensions to unbounded elements 
affiliated with  $C^*$-algebras (see \cite{w1,w11}).  If $\rho$ is a 
strictly positive self-adjoint operator affiliated with $A$, then 
$\seq{\rho^{it}}{t\in\real}$ is a one parameter group of unitaries in 
$M(A)$ continuous with respect to the strict topology and 
$\pi(\rho)^{it}=\pi(\rho^{it})$. In this formula $\pi(\rho)$ denotes 
the element affiliated with $B$ corresponding (via morphism 
$\pi\in\mathrm{Mor}(A,B)$) to $\rho$ affiliated with $A$.

\setcounter{equation}{0}\setcounter{Thm}{0}
\subsection{GNS maps}
\label{GNSmap}
This Appendix is devoted to the notion of GNS map. Neglecting certain 
details this notion is equivalent to the notion
of generalized Hilbert algebra. In our opinion it is more convenient, 
when we have to deal with many different scalar
products defined on the same operator algebra. \vs

In the following definition operator algebra means either 
C$^*$-algebra or von Neumann algebra of operators acting on a
Hilbert space $\Hil$. In the C$^*$-algebra case we shall assume that 
it acts on $\Hil$ in a non-degenerate way.
Consequently in the case of von Neumann algebra we understand that it 
contains the identity operator $1_\Hil$.

\begin{Def}
Let $B$ be an operator algebra acting on a Hilbert space $\Hil$ and 
$\eta$ be a (unbounded) linear mapping from $B$
into $\Hil$. We say that $\eta$ is a GNS map if the domain 
$\Dom(\eta)$ is a left ideal in $B$ and
\begin{equation}
\label{GNS}
\eta(ab)=a\eta(b)
\end{equation}
for any $a\in B$ and $b\in\Dom(\eta)$. The GNS map $\eta$ is called 
closed if   it is closed with respect to strong
operator topology on $B$ and norm topology on $\Hil$. We say that a 
GNS map defined on a C\,$^*$-algebra is densely
defined if its domain is norm dense in the C\,$^*$-algebra. Similarly 
a GNS map defined on a von Neumann algebra is
said to be densely defined if its domain is strongly dense in the von 
Neumann algebra.
\end{Def}

The reader should notice that using the term 'densely defined` we 
refer to different topologies in C$^*$-algebra and
von Neumann algebra cases. On the contrary to define closedness we 
use in both cases the same strong operator topology
on $B$. The same topology is used when we speak about continuous GNS 
maps, closable GNS maps, the closure of a GNS map
and a core (essential domain) of a GNS map.\vs

Let $B$ be an operator algebra acting on a Hilbert space $\Hil$ and 
$\Omega\in\Hil$. For any $b\in B$ we set
\[
\eta(b)=b\Omega.
\]
Then $\eta$ is a continuous GNS map with $\Dom(\eta)=B$. To obtain 
less trivial example assume that $\Hil$ is a direct
orthogonal sum of a family $\seq{\Hil_i}{i\in I}$ of $B$-invariant 
subspaces and for each $i\in I$ pick up a vector
$\Omega_i\in\Hil_i$. Then the formula
\begin{equation}
\label{postac}
\eta(b)=\sum_{i\in I}b\Omega_i.
\end{equation}
defines a closed GNS map. By definition $\Dom(\eta)$ is the set of 
all $b\in B$ for which the above series is norm
convergent. It turns out (cf. Theorem \reef{Bszesc}) that any closed 
GNS map defined on a von Neumann algebra acting on a
separable Hilbert space is of the form \rf{postac} with a denumerable 
index set $I$. \vs

Now we introduce the notion of the commutant $\eta'$ of a GNS-map 
$\eta$. Let $B$ be an operator algebra acting on a
Hilbert space $\Hil$, $\eta$ be a GNS map defined on  $B$ and $B'$ be 
the commutant of $B$. Assume that the
intersection of kernels of all $b\in\Dom(\eta)$ is trivial. This is 
always the case, when $\eta$ is densely defined.
For any $c'\in B'$ and $x\in\Hil$ we set:
\begin{equation}
\label{komutant}
\left(
\begin{array}{c}
   c'\in\Dom(\eta') \\
  \mbox{and} \\
   x=\eta'(c')
\end{array}
\right) \Longleftrightarrow
  \left(
\begin{array}{c}
   c'\eta(b)=bx \\
   \mbox{for all} \\
   b\in\Dom(\eta)
\end{array}
\right).
\end{equation}
The reader should notice that there is at most one vector $x$ 
satisfying the right hand side of the above equivalence.
One can easily verify that $B'$ is a von Neumann algebra and that 
$\eta'$ is a closed GNS map defined on $B'$. It turns
out (cf. Theorem \reef{podkomutant}) that $\eta'$ is densely defined 
provided $\eta$ is closable.\vs

In the following statements we use the double commutant $\eta''$ of a 
GNS map $\eta$. If $\eta$ is defined on an
operator algebra $B$ then $\eta''$ is defined on the double commutant 
$B''$. By the famous von Neumann density theorem,
$B''$ is the strong closure of $B$. We have the following

\begin{Thm}
\label{podkomutant} Let $B$ be an operator algebra acting on a Hilbert space
$\Hil$ and $\eta$ be a closable GNS mapping densely defined on $B$. Then
$\eta'$ is densely defined on the von Neumann algebra $B'$, $\eta''$ is an
extension of $\eta$ and $\Dom(\eta)$ is a core of $\eta''$.
\end{Thm}

\begin{pf} The inclusion $\eta\subset\eta''$ is obvious. Let $\overline{\eta}$
be the closure of $\eta''$ restricted to $\Dom(\eta)$, $a\in B''$ and
$x\in\Hil$. We have to show that
\begin{equation}
\label{komutant1}
  \left(
\begin{array}{c}
   a\eta'(p)=px \\
   \mbox{for all} \\
   p\in\Dom(\eta')
\end{array}
\right)
  \Longrightarrow
\left( \begin{array}{c}
   a\in\Dom(\overline{\eta}) \\
  \mbox{and} \\
   x=\overline{\eta}(a)
\end{array}
\right).
\end{equation}
The condition on the right hand side of \rf{komutant1} means that the 
pair $(a,x)$ belongs to the closure of the graph
of $\eta$ in the (strong$\,\times\,$norm)-topology. Remembering that 
the neighborhoods in strong topology are
determined by finite sets of vectors it is sufficient to show that 
for any $\Omega_1,\Omega_2,\dots,\Omega_n\in\Hil$
and
\begin{equation}
\label{d111}
\begin{array}{c}
\mbox{for any }\varepsilon>0\mbox{ there exists }b\in\Dom(\eta)\mbox{ 
such that} \\
\Vs{5}\norm{x-\eta(b)}^2+\sum_{k=1}^n\norm{a\Omega_k-b\Omega_k}^2<\varepsilon^2.
\end{array}
\end{equation}

Let $\widetilde{\Hil}$ be the direct sum of $1+n$ copies of
$\Hil$. The copies will be labelled by numbers $0,1,2,\dots,n$.
For any $k=0,1,\dots,n$ we denote by $x_k\in\Hil$ the $k$-th
component of a vector $x\in\widetilde{\Hil}$. Consequently for any
$b\in\skrL(\widetilde{\Hil})$, the matrix elements of $b$ will be
denoted by $b_{kl}\in\skrL(\Hil)$ ($k,l=0,1,2,\dots,n$). In what
follows we shall use the diagonal action of  $B$ on
$\widetilde{\Hil}$: any $b\in B$ defines an operator
$\tilde{b}\in\skrL(\widetilde{\Hil})$ with matrix elements
$\tilde{b}_{kl} =b\delta_{kl}$ ($k,l=0,1,2,\dots,n$).\vs

For any $b\in\Dom(\eta)$ we denote by $\tilde{\eta}(b)$ the element 
of $\widetilde{\Hil}$ such that
\[
\tilde{\eta}(b)_k=\left\{
\begin{array}{cl}
\eta(b)&\mbox{ for }k=0\\ b\Omega_k&\mbox{ for }k=1,2,\dots,n.
\end{array}
\right.
\]
We shall also use the vector $\tilde{x}\in\widetilde{\Hil}$ such that
\[
\tilde{x}_k=\left\{
\begin{array}{cl}
x&\mbox{ for }k=0\\ a\Omega_k&\mbox{ for }k=1,2,\dots,n.
\end{array}
\right.
\]
With this notation, \rf{d111} may be rewritten in the following 
equivalent form:
\begin{equation}
\label{przestrzen} \tilde{x}\in\set{\tilde{\eta}(b)}{\Vs{3.5}b\in 
B}^{\rm closure}
\end{equation}

Let $p\in\skrL(\widetilde{\Hil})$ be the orthogonal projection onto 
orthogonal complement of
$\set{\tilde{\eta}(b)}{\Vs{3}b\in B}$. One can easily show that the 
latter space is invariant under
the diagonal action of $B$. Therefore the matrix elements of $p$ 
belongs to the commutant of $B$:
$p_{kl}\in B'$ for all $k,l=0,1,2,\dots,n$. By definition 
$p\tilde{\eta}(b)=0$. It means that
\[
p_{k0}\eta(b)+\sum_{l=1}^n p_{kl}b\Omega_l=0
\]
and
\[
p_{k0}\eta(b)=-b\sum_{l=1}^n p_{kl}\Omega_l
\]
for $k=0,1,2,\dots,n$. This relation holds for all $b\in\Dom(\eta)$. 
Definition \rf{komutant} shows now that
$p_{k0}\in\Dom(\eta')$ and $\eta'(p_{k0})=-\sum_{l=1}^n p_{kl}\Omega_l$.\vs

Assume now that the condition on the left hand side of \rf{komutant1} 
is satisfied. Then
\[
p_{k0}x=a\eta'(p_{k0})=-a\sum_{l=1}^n p_{kl}\Omega_l
\]
and
\[
p_{k0}x+\sum_{l=1}^n p_{kl}a\Omega_l=0.
\]
This relation holds for $k=0,1,2,\dots,n$. It shows that 
$p\tilde{x}=0$ and \rf{przestrzen} follows. This way the
implication \rf{komutant1} is proved.\vs

  Setting $a=0$ we see that $x=0$ is the only vector killed by all
$p\in\Dom(\eta')$. Remembering that $\Dom(\eta')$ is a left ideal in 
$B'$ and using \cite[Chapitre I, \S 3, Corollaire 3]{Dix} we see that
$\Dom(\eta')$ is dense in $B'$. Comparing now \rf{komutant1} with 
\rf{komutant} we see that $\overline{\eta}$ is an
extension of $\eta''$: $\eta''\subset\overline{\eta}$. The converse 
inclusion is obvious: $\eta''$ is an extension of
$\eta$ and $\eta''$ is closed. Hence $\overline{\eta}\subset\eta''$.

\end{pf}
\vs

Assume now that $\Hil$ is separable. Our nearest aim is to show that 
any closed GNS map is of the form \rf{postac}. We
start with the following

\begin{Prop}
\label{jadro}
  Let $M$ a von Neumann algebra acting on a separable Hilbert space $\Hil$ and
$\eta$ be a closed GNS-mapping on $M$.  Then there exists 
$d\in\Dom(\eta)$ such that $Md$ is a core for $\eta$.  One
may choose $d$ in such a way that $d\geq 0$ and ${\rm 
Sp}(d)\subset\{1,2^{-1},2^{-2},\dots,0\}$.
\end{Prop}

\pf At first we shall prove that there exists elements 
$c_n\in\Dom(\eta)$ ($n=1,2,\dots$) such that the set
\begin{equation}
\label{corex}
  \set{c_n}{n=1,2,\dots\Vs{3}}^{\rm linear\ span}
\end{equation}
is a core for $\eta$. We know that $\eta$ is closed with respect to 
the strong topology on $M$ and norm topology on
$\Hil$. We shall use the graph topology on $\Dom(\eta)$. Let
\[
K=\set{a\in\Dom(\eta)}{\norm{a}\leq 1,\ \norm{\eta(a)}\leq 1}.
\]
Clearly, $K$ is a closed subset of $\Dom(\eta)$. Choosing an 
orthonormal basis $\seq{e_n}{n\in\natu}$ in $\Hil$ and
setting
\[
d(a,b)^2=\sum_{n=1}^{\infty}\frac 1{4^n}\norm{(a-b)e_n}^2+\norm{\eta(a-b)}^2,
\]
we define a metric $d$ on $K$.  One can easily show that the topology 
induced by this metric coincides with the graph
topology on $\Dom(\eta)$ restricted to $K$. Moreover, the mapping
\[
K\ni a\longmapsto \left(\sum_{n=1}^{\infty}\Vs{4}^{\oplus} \frac 
1{2^n}ae_n\right)\oplus\eta(a)
\in\left(\sum_{n=1}^{\infty}\Vs{4}^{\oplus}\Hil\right)\oplus\Hil
\]
is an isometric embedding and the target Hilbert space is separable. 
Therefore there exists a denumerable subset dense
in $K$. The linear span of this subset is the core for $\eta$ of the 
form (\reef{corex}).

  Multiplying if necessary  by a decreasing sequence of positive 
numbers we may assume that $\sum
\norm{c_n}<\infty$ and $\sum \norm{\eta(c_n)}<\infty$.  Using 
sequence $\seq{c_n}{n=1,2,\dots}$, we define a bounded linear
mapping $c$
\[
\Hil\ni x\longmapsto cx=\sum_{n=1}^{\infty}\Vs{4}^{\oplus}c_nx.
\]
Then $c^*c=\sum_{n=1}^{\infty}c_n^*c_n\in M$. Let $d=\sqrt{c^*c}$ and 
$c=ud$ be the polar decomposition of $c$. Then
$d\in M$ and the isometry $u$ is of the form
\[
\Hil\ni x\longmapsto ux=\sum_{n=1}^{\infty}\Vs{4}^{\oplus}u_nx 
\in\sum_{n=1}^{\infty}\Vs{4}^{\oplus}\Hil,
\]
where all $u_n\in M$ and $\norm{u_n}\leq 1$. Relation $c=ud$ means that
\begin{equation}
\label{idea} c_n=u_nd
\end{equation}
for all $n=1,2,\dots$, whereas relation $d=u^*c$ is equivalent to
\[
d=\sum_{n=1}^\infty u_n^*c_n.
\]
The above series is norm converging. Also the series $\sum 
u_n^*\eta(c_n)$ is norm converging. Remembering that $\eta$
is closed we conclude that $d\in\Dom(\eta)$ and $Md\subset\Dom(\eta)$ 
($\Dom(\eta)$ is a left ideal!). Formula
(\reef{idea}) shows now that the set (\reef{corex}) is contained in 
$Md$. Therefore $Md$ is a core for $\eta$.

The element $d$ constructed above is self-adjoint and positive. We 
may assume that $\norm{d}< 1$. Let $f$ be a simple
function on $[0,1)$ defined by $f(0)=0$ and
\[
f(t)=\frac{1}{2^{n+1}} \quad \mbox{for} \quad t\in \left[\frac 
1{2^{n+1}},\frac 1{2^n}\right).
\]
for $n=0,1,2,\dots$. Then $t/2\leq f(t)\leq t$ for $t\in[0,1]$. Using 
this estimate one can easily show that $d$ and
$f(d)$ generate the same left ideal in $M$: $Md=Mf(d)$. Moreover 
$f(d)\geq 0$ and
$\Spec f(d)\subset\{1,2^{-1},2^{-2},\dots,0\}$. This way the last 
Statement of our Proposition is shown.

\qed

\begin{Def}
\label{Bsz} Let $B$ be an operator algebra acting on a separable 
Hilbert space $\Hil$ and $\eta$ be a
closed densely defined GNS-mapping on $B$. Assume that 
$\seq{\Omega_n}{n=1,2,\dots}$ is a
sequence of elements of $\Hil$  such that $B\Omega_n\perp B\Omega_m$ 
for $n\neq m$ and
\begin{equation}
\label{cyklik}
\eta(b)=\sum_n^\infty b\Omega_n.
\end{equation}
for any $b\in \Dom(\eta)$. Then we say that \rf{cyklik} is a
{\rm vector presentation} of the GNS mapping $\eta$.
\end{Def}

\begin{Prop}
\label{Bpiecc}
Let \rf{cyklik} be a vector presentation of a GNS map $\eta$ of a 
closed densely
defined GNS map $\eta$. Then for any $b\in\Dom(\eta'')$ we have
\begin{equation}
\label{cyklikbis}
\eta''(b)=\sum_{n=1}^\infty b\Omega_n.
\end{equation}
\end{Prop}

\begin{pf} Let  $P_n$ be the
orthogonal projection onto $\overline{B\Omega _n}$ ($n=1,2,\dots$). 
Then $P_n$ ($n=1,2,\dots$) are
pairwise orthogonal projections and denoting by $P_0$ the projection 
onto the orthogonal complement
of ${\sum}^\oplus \overline{B\Omega _n}$ we obtain the decomposition of unity:
\[
1=\sum_{n=0}^\infty P_n
\]
Clearly all $P_n$ belong to $B'$. For $n=1,2,\dots$ we have
\[
P_n\eta(b)=P_n\sum_{n=1}^\infty b\Omega_n=b\Omega_n
\]
for all $b\in B$. It shows that $P_n\in\Dom(\eta')$ and 
$\eta'(P_n)=\Omega_n$. Moreover
$P_0\eta(b)=0$ for all $b\in B$. Therefore $P_0\in\Dom(\eta'$ and 
$\eta'(P_0)=0$.\vs

Assume that $b\in\Dom(\eta'')$. Then for any $c\in\Dom(\eta')$ we have
$c\eta''(b)=b\eta'(c)$. Setting $c=P_n$ we obtain 
$P_n\eta''(b)=b\eta'(P_n)=b\Omega_n$ for
$n=1,2,\dots$ and $P_0\eta''(b)=b\eta'(P_0)=0$ for $n=0$. Summing 
over $n$ we get \rf{cyklikbis}.

\end{pf}

\begin{Def}
\label{defexact}
A vector presentation \rf{cyklik} is said to be {\rm exact} if $\Dom(\eta'')$
coincides with the set of
all elements $b\in B''$ such that the sum \rf{cyklikbis} is convergent.
\end{Def}

We know that any bounded GNS-map $\eta$ is of the form 
$\eta(a)=a\Omega$, where $\Omega\in\Hil$. We shall use the
Proposition \reef{jadro} to prove the following nice generalization 
of this fact.

\begin{Thm}
\label{Bszesc}
Let $B$ be an operator algebra acting on a separable Hilbert space 
$\Hil$ and $\eta$ be a closed
densely defined GNS-mapping on $B$. Then $\eta$ admits an exact 
vector presentation.
\end{Thm}

\begin{pf} Let $\eta'$ be the commutant of $\eta$. By Proposition 
\reef{jadro}, there exists
$d\in\Dom(\eta')$ such that $B'd$ is a core for $\eta'$. Assume that 
a vector $x\in\ker d$. Then
$cdx=0$ for any $c\in B'$. Therefore $x$ is killed by all elements of 
$\Dom(\eta')$. Remembering that
$\Dom(\eta')$ is strongly dense in $B'$ we see that $x=0$. This way 
we have shown that $\ker
d=\{0\}$.

By the last sentence of Proposition \reef{jadro}, we may assume that 
$d$ is a self-adjoint positive operator with
spectrum contained in $\{1,2^{-1},2^{-2},\dots,0\}$. Then
\[
d=\sum_{n=1}^\infty 2^{-n}P_n,
\]
where $P_n\in B'$ are orthogonal projections such that $P_nP_m=0$ for 
$n\neq m$. Remembering that the kernel of $d$ is
trivial we obtain
\[
\sum_{n=1}^\infty P_n=I
\]
For any $n=1,2,\dots$, we set:
\[
\Omega_n=2^nP_n\eta'(d).
\]
Then $d\Omega_n=P_n\eta'(d)$ and
\begin{equation}
\label{xx} \sum_{n=1}^\infty d\Omega_n=\sum_{n=1}^\infty P_n\eta'(d)=\eta'(d).
\end{equation}
The reader should notice that the above series is norm converging.
  Since $P_n\in B'$, the subspace $P_n\Hil$ is $B$-invariant and 
$B\Omega_n\subset P_n\Hil$.
Remembering that $P_nP_m=0$ we see that  $B\Omega_n\perp B\Omega_m$ 
for $n\neq m$. Formula $P_n=2^nP_nd$ shows that
$P_n\in\Dom(\eta')$ and $\eta'(P_n)=2^nP_n\eta'(d)=\Omega_n$ for all 
natural $n$.

Let $b\in\Dom(\eta)$. Then, for any natural $n$, 
$P_n\eta(b)=b\eta'(P_n)=b\Omega_n$. Therefore
\[
\eta(b)=\sum_{n=1}^\infty P_n \eta(b)=\sum_{n=1}^\infty b\Omega_n.
\]

Clearly the above formula is a vector presentation of $\eta$. If the series
\[
x=\sum_{n=1}^\infty b\Omega_n
\]
is converging for some $b\in B''$, then, using (\reef{xx}) we obtain:
\[
dx=\sum_{n=1}^\infty db\Omega_n=\sum_{n=1}^\infty 
bd\Omega_n=b\sum_{n=1}^\infty d\Omega_n=b\eta'(d).
\]
Therefore for any $c\in B'$, $cdx=cb\eta'(d)=bc\eta'(d)=b\eta'(cd)$. 
Remembering that $B'd$ is a core for $\eta'$ we
see that $cx=b\eta'(c)$ for all $c\in\Dom(\eta')$. Hence 
$b\in\Dom(\eta'')$. It shows that the vector presentation
constructed above is exact.

\end{pf}

Let us also note the following interesting

\begin{Prop}
\label{18.08.2002c}
Let $B$ an operator algebra acting on a separable Hilbert space
$\Hil$, $\eta$ be a closed densely defined GNS mapping on $B$ with a
vector presentation
\begin{equation}
\label{vectpres}
\eta (b)=\sum_nb\Omega_n
\end{equation}
and $P_n\in\Dom (\eta')$ be the orthogonal projection onto
$\overline{B\Omega_n}$, $n=1,2,\dots$ $($see the proof of Proposition
\reef{Bpiecc}$)$. Then the following statements are equivalent:
\begin{enumerate}
\item The vector presentation \rf{vectpres} is exact.
\item The left ideal of $B'\Vs{5}$ generated by $\{P_1,P_2,\dots\}$ is a
core for $\eta'$.
\end{enumerate}

\end{Prop}
\begin{pf}
Ad 2 $\Rightarrow$ 1.

Let $b$ be an element of $B''$. If the series $\sum_kb\Omega_k$ is
convergent, then for any $a\in B'$ and any $n=1,2,\dots$ we have:
\[
b\eta'(aP_n)=ab\eta'(P_n)=ab\Omega_n=aP_n\sum_kb\Omega_k \;.
\]
If the left ideal in $B'$ generated by $\{P_1,P_2,\dots\}$ is a
core for $\eta'$ then $b\eta'(a)=a\sum_kb\Omega_k$ for all
$a\in\Dom (\eta )$. Hence $b\in\Dom (\eta'')$ and
$\eta''(b)=\sum_kb\Omega_k$. It shows that the vector presentation
\rf{vectpres} is exact.\vs

Ad 1 $\Rightarrow$ 2.

Assume now that the left ideal $I$ of $B'$ generated by
$\{P_1,P_2,\dots\}$ is not a core for $\eta'$. Let $\theta$ be the
closure of $\eta'$ restricted to $I$. Then
$\theta :B'\rightarrow\Hil$ is a GNS mapping, $\theta\subset\eta'$
and $\theta\neq\eta'$. Therefore $\eta''\subset\theta'$ and
$\eta''\neq\theta'$. Hence there exists $b\in\Dom (\theta')$ such
that $b\notin\Dom (\eta'')$. Clearly $P_n\in\Dom (\theta )$ and
$\theta (P_n)=\eta'(P_n)=\Omega_n$. Now we have
\[
b\Omega_n=b\theta (P_n)=P_n\theta'(b)
\]
and the series $\sum b\Omega_n=\sum_nP_n\theta'(b)$ is convergent.
Remember that $b\notin\Dom (\eta'')$ we see that the vector
presentation \rf{vectpres} is not exact.
\end{pf}

\setcounter{equation}{0}\setcounter{Thm}{0}
\subsection{Weights on C$^*$-algebras}
\label{B}

In this Appendix, we collect some basic properties of weights on a separable
$C^*$-algebras and their relations with GNS mappings. Let $h$ be a weight on a
$C^*$-algebra $A$. We shall use the standard notation:
\[
\begin{array}{r@{\;}c@{\;}l}
\Nideal_h&=&\set{a\in A}{h(a^*a)<\infty},\\
\Mideal_h&=&\Vs{5}\set{a\in A_+}{h(a)<\infty}^{\rm linear\ span}\\
&=&\Vs{5}\set{a^*b}{a,b\in\Nideal_h}^{\rm linear\ span}
\end{array}
\]
Then $\Nideal_h$ is a left ideal and $\Mideal_h$ is a hereditary 
subalgebra in $A$.
It is known that $h$ extends uniquely to a linear functional on 
$\Mideal_h$. We say that $h$ is locally finite if
$\Nideal_h$ is dense in $A$.\vs

By the GNS construction one can find a Hilbert space $\Hil$, a representation
$\pi$ of $A$ acting on $\Hil$ and a GNS-map $\eta$ defined on the C$^*$-algebra
$B=\pi(A)$ such that
\[
\Dom(\eta)=\set{\pi(b)}{b\in\Nideal_h},
\]
the range of $\eta$ is dense in $\Hil$ and
\[
\norm{\eta(\pi(a))}^2=h(a^*a)
\]
for $\Vs{4}$all $a\in\Nideal_h$. We say that $(\Hil,\pi,\eta)$ is a 
GNS triple associated with
the weight $h$. It is unique up to a unitary equivalence. If $h$ is 
finite then the
GNS-map $\eta$ is bounded, defined on whole $\pi(A)$ and there exists 
a vector $\Omega\in\Hil$ such
that $\eta(b)=b\Omega$ for all $b\in B$. In the general case we would 
like to have the map $\eta$ closed. This is the
case if the weight $h$ is lower semicontinuous, i.e. if for any 
$t\in\real_+$ the set $\set{a\in A_+}{h(a)\leq t}$ is
closed in $A$ in the sense of the norm topology. By a very important 
result of Haagerup \cite{h}, any lower semicontinuous
weight
$h$ on a C$^*$-algebra $A$ is of the form
\begin{equation}
\label{suma0}
h=\sum_{i\in I}h_i,
\end{equation}
where $\seq{h_i}{i\in I}$ is a family of positive functionals on $A$. 
For separable C$^*$-algebras the
above sum is countable:\vs

\begin{Lem}\label{lem:sls} Let $h$ be a locally finite lower 
semicontinuous weight on a
separable $C^*$-algebra $A$. Then there exists a countable family 
$\seq{h_i}{i\in\natu}$ of
elements of $A^*_+$ such that
\begin{equation}
\label{suma1}
h =\sum_{i\in \natu}h_i \;.
\end{equation}
\end{Lem}

\begin{pf}
We start with formula \rf{suma0}. Replacing if necessary $h_i$ by the sum of
$N$-copies of $h_i/N$, (where $N$ is a natural number larger than 
$\norm{h_i}$) and removing
zero-terms, we may assume that
$\norm{h_i}\leq 1$ and $h_i\neq 0$ for all $i\in I$. Let
$\Lambda$ be the set $\set{\psi\in A^*_+ }{\norm{\psi}\leq 1}$. Since
$A$ is separable, the $\text{weak}^*$ topology on the unit ball of
$A^*$ is metrizable. Denote a metric by $d$. Put
\[
\begin{array}{r@{\;}c@{\;}l}
  I_1&=&\set{i\in I}{\Vs{4}d(h_i ,0)\in [1,\infty )}
\\ I_n&=\Vs{7}&\set{i\in I}{d(h_i ,0)\in\left[ \frac 1{n},\frac 
1{n-1} \right)} \;,\; n\geq 2\ .
\end{array}
  \]
It suffices to show that each $I_n$ is finite. Suppose that $I_n$ is 
infinite for some
$n\in \natu$. Since $\Lambda$ is $\text{weak}^*$ compact, there 
exists a sequence $i_m\in I_n$
($m=1,2,\dots$) and an element $\psi\in\Lambda$ such that 
$h_{i_m}\longrightarrow \psi$ when
$m\rightarrow\infty$. Clearly $d(\psi ,0)\geq n^{-1}$ and $\psi\neq 
0$. For any $a\in\Nideal_h$ we have
\[
\sum_{m=1}^\infty h_{i_m}(a^*a)\leq\sum_{i\in I}h_i(a^*a)
     =h(a^*a)<\infty.
\]
It shows that $\psi (a^*a)=\lim_{m\to\infty}h_{i_m}(a^*a)=0$. Remembering
that $\Nideal_h$ is norm-dense in $A$ we conclude that $\psi=0$. This 
contradiction shows that $I_n$ is finite for all
$n\in\natu $.

\end{pf}

\begin{Thm}\label{thm:omega}
Let $h$ be a locally finite lower
semicontinuous weight on a separable $C^*$-algebra $A$ and
$(\Hil,\pi,\eta)$ be the GNS-triple associated with the weight $h$. 
Then $\Hil$ is
separable and $\eta$ is closed.
\end{Thm}

\begin{pf} We may assume that $h$ is given by \rf{suma1}. Let
$(\Hil_i,\pi_i,\eta_i)$be the GNS triple associated with
$h_i$. In this case $\eta_i$ are continuous. Remembering that
$A$ is separable and that the range of $\eta_i$ are dense in
$\Hil_i$ we conclude that  $\Hil_i$ are separable. By Lemma
\reef{lem:sls}, the index set $I$ is denumerable. Therefore the direct sum
\[
\widetilde{\Hil}=\sum_{i\in I}\Vs{4}^\oplus \Hil_i
\]
is a separable Hilbert space. Using
\rf{suma1} one can easily show that
\[
\left\norm{\eta(\pi(a))\right}^2=\sum_{i\in 
I}h_i(a^*a)=\left\norm{\sum_{i\in I}\Vs{4}^\oplus 
\eta_i(\pi_i(a))\right}^2
\]
for any $a\in\Nideal_h$. Remembering that the range of $\eta$ is 
dense in $\Hil$ we see that
there exists a unique isometry 
$U:\Hil\longrightarrow\widetilde{\Hil}$ such that
\[
U\eta(\pi(a))=\sum_{i\in I}\Vs{4}^\oplus \eta_i(\pi_i(a))\; .
\]
Therefore $\Hil$ is isomorphic to a subspace of a separable Hilbert 
space. It shows that $\Hil$ is
separable. \vs

Let $i\in I$. Then $h_i(a^*a)\leq h(a^*a)$ for any $a\in\Nideal_h$. 
Combining this
relation with the Schwartz inequality one can easily verify that $|h_i(a)|\leq
\sqrt{\norm{h_i}}\norm{\eta(\pi(a))}$. Therefore there exists a 
vector $\Omega_i\in\Hil$ such that
\begin{equation}
\label{slw1}
h_i(a)=\IS{\Omega_i}{\eta(\pi(a))}
\end{equation}
for any $a\in\Nideal_h$. One can also easily show that
$\left(\eta(\pi(c))\Vs{3},\eta(\pi(b))\right)\longrightarrow 
h_i(c^*b)$ is a continuous
sesquilinear form on the range of $\eta$. Therefore there exists a 
bounded positive operator $T_i$
acting on $\Hil$ such that
\begin{equation}
\label{slw2}
h_i(c^*b)=\ITS{\eta(\pi(c))\Vs{3}}{T_i}{\eta(\pi(b))}
\end{equation}
for any $c,b\in\Nideal_h$.

We shall prove that for any  $x\in\Hil$ and $b\in B=\pi(A)$ we have 
the following equivalence:
\begin{equation}
\label{equix}
\left(\begin{array}{c} b\in\Dom(\eta)\\ \mbox{and} \\ x=\eta(b)
\end{array}\right)
\Longleftrightarrow \left(\begin{array}{c} T_ix=b\Omega_i\\ \mbox{for 
all }  i\in I
\end{array}\right).
\end{equation}

Replacing in \rf{slw1} $a$ by $c^*b$ and comparing with \rf{slw2} we obtain
$T_i\eta(\pi(c))=\pi(c)\Omega_i$ for any $c\in\Nideal_\varphi$. It shows that
\begin{equation}
\label{slw3} T_i\eta(b)=b\Omega_i
\end{equation}
for any $b\in\Dom(\eta)$. This way we showed the `$\Rightarrow$'- 
part of \rf{equix}.

We shall prove the converse. Setting in \rf{slw2}, $c=b$ and summing 
over $i\in I$ we see that
\begin{equation}
\label{slw4}
\sum_{i\in I}T_i=1
\end{equation}
  Assume now that $x\in\Hil$, $a\in A$ and $T_ix=\pi(a)\Omega_i$ for 
all $i\in I$ and take a sequence
$\seq{a_n}{n=1,2,\dots}$ of elements of $\Nideal_h$ converging to $a$ 
in norm. Using formula
\rf{slw2} one can easily show that the sequence
$\seq{T_i^{\frac{1}{2}}\eta(\pi(a_n))}{n=1,2,\dots}$ is convergent:
\[
\lim_{n\rightarrow\infty} T_i^{\frac{1}{2}}\eta(\pi(a_n))=y_i,
\]
where $y_i\in\Hil$ and $\norm{y_i}^2=h_i(a^*a)$. Moreover, by formula 
\rf{slw3},
\[
T_i^{\frac{1}{2}}y_i=\lim_{n\rightarrow\infty}
T_i\eta(\pi(a_n))=\lim_{n\rightarrow\infty}\pi(a_n)\Omega_i=\pi(a)\Omega_i=T_ix
\]
It shows that $y_i-T_i^{\frac{1}{2}}x\in\ker T_i^{\frac{1}{2}}$. On 
the other hand $y_i$ and
$T_i^{\frac{1}{2}}x$ clearly belong to the closure of the range of 
$T_i^{\frac{1}{2}}$. Therefore
$y_i-T_i^{\frac{1}{2}}x=0$ and $y_i=T_i^{\frac{1}{2}}x$. Now we have:
\[
h_i(a^*a)=\norm{y_i}^2=\ITS{x}{T_i}{x}.
\]
Summing over $i\in I$ and using \rf{slw4} we obtain
\[
h(a^*a)=\sum_{i\in I}h_i(a^*a)=\IS{x}{x}<\infty.
\]
It shows that $\pi(a)\in\Dom(\eta)$. Now we may use formula 
\rf{slw3}. It shows that
$T_i(\eta(\pi(a))-x)=0$ and taking into account \rf{slw4} we 
conclude: $x=\eta(\pi(a))$. The
`$\Leftarrow$' of \rf{equix} is shown.

The closedness of $\eta$ is an easy consequence of \rf{equix}. Let
$\seq{b_\lambda}{\lambda\in\Lambda}$ be a net of elements of 
$\Dom(\eta)$ converging strongly
to an element $b_\infty\in B$, such that 
$x_\lambda=\eta\left(b_\lambda\right)$ converges to a
vector $x_\infty\in\Hil$. Then $T_ix_\lambda=b_\lambda\Omega_i$. 
Passing to the limit we get
$T_ix_\infty=b_\infty\Omega_i$. It means that $b_\infty\in\Dom(\eta)$ 
and $\eta(b_\infty)=x_\infty$.

\end{pf}

\begin{Thm}
\label{wagaaa}
Let $B$ a C$^*$-algebra acting on a separable Hilbert space $\Hil$ and
$\eta:B\to\Hil$ be a closed densely defined GNS map.  Then the formula$\,:$
\begin{equation}
\label{waga}
h(b^*b)=\left\{
\begin{array}{cl}
\IS{\eta''(b)}{\eta''(b)}&\hspace{5mm}\mbox{\rm if 
}b\in\Dom(\eta'')\\ +\infty &\hspace{6mm} \mbox{otherwise}\Vs{5}
\end{array}
\right. ,
\end{equation}
where $b$ runs over $B''$, defines a normal semifinite, weight $h$ on 
$B''$. Its restriction to $B$ is a locally finite
lower semicontinuous weight on $B$.
\end{Thm}

\begin{pf} Let $M=B''$. We shall use Theorem \reef{Bszesc}. Let 
\rf{cyklik} be an exact vector presentation of $\eta$.
Combining
\rf{cyklikbis} with \rf{waga} we obtain:
\[
h(b^*b)=\sum_n \IS{\Omega_n}{b^*b\,\Omega_n}
\]
for all $b\in\Dom(\eta'')$. Due to the exactness, this formula holds 
for all $b\in M$. Therefore for any $b\in M_+$ we
have:
\[
h(b)=\sum_n \IS{\Omega_n}{b\Omega_n}.
\]
It shows that $h$ is a normal weight on $M$ (cf. \cite{h}). Clearly 
$h(b^*b)<\infty$ for $b\in\Dom(\eta)$. The latter set
is weakly dense in $M$. Hence $h$ is semifinite. The last assertion is obvious.

\end{pf}

\begin{Def} A weight $h$ on a $C^*$-algebra $A$ is said to be
\begin{enumerate}
\item {\rm faithful} if for any $a\in A$, the relation $h(a^*a)=0$ implies
$a=0$.
\item {\rm strictly faithful} if for any sequence $\seq{a_n}{n\in\natu }$ in
$A$ such that the sequence $\seq{h(a_n^*a_n)}{n\in\natu }$ is bounded and
$\lim_{n\to\infty}h(a_na_n^*)=0$, $b=0$ is the only element of $A$ satisfying
the relation
\begin{equation}
\label{strfaith}
h((b-a_n)^*(b-a_n))\leq h(a_n^*a_n)
\end{equation}
for all $n\in\natu$.
\end{enumerate}
\end{Def}

Assume that $h$ is strictly faithful. Inserting $a_n=0$ in \rf{strfaith} we see
that $b=0$ is the only element of $A$ such that $h(b^*b)=0$. It shows that any
strictly faithful weight is faithful.

\begin{Lem}\label{lem:strictfaithful}
Let $h$ be a locally finite lower semicontinuous weight on a 
C$^*$-algebra $A$ and $(\Hil,\pi,\eta)$ be the GNS triple
related to $h$.Then the following conditions are equivalent$\,:$
\begin{enumerate}
\item The weight $h$ is strictly faithful.
\item The weight $h$ is faithful and the mapping
$S_0:\eta(a)\mapsto\eta(a^*)$ {\rm (}where $a$ runs over 
$\Dom(\eta)\cap\Dom(\eta)^*${\rm )} is a closable conjugate
linear operator in $\Hil$.
\item The double commutant $\eta''$ is injective.
\end{enumerate}
\end{Lem}

\begin{pf} Assume that the mapping $\eta(a)\mapsto\eta(a^*)$
  is not closable. It is the same as to assume that the mapping
$\eta(a^*)\mapsto\eta(a)$ is not closable. Therefore there
exists a sequence
$\seq{a_n}{n\in \natu }$ in $\Dom(\eta)\cap\Dom(\eta)^*$ such that
$\seq{\eta(a_n^*)}{n\in \natu }$ in the Hilbert space $\Hil$ 
converges to zero in norm topology (i.e.
the sequence $\seq{h(a_na_n^*)}{n\in \natu }$ of non-negative numbers 
converges to zero) and
$\seq{\eta(a_n)}{n\in\natu }$ in $\Hil$ converges to some non-zero 
element $x\in\Hil$ (in particular the sequence
$\seq{h(a_n^*a_n)}{n\in \natu }$ of non-negative numbers is bounded). 
Remembering that
the image of the GNS map $\eta$ is dense in $\Hil$ we find $b\in 
\Dom(\eta)$ such that
$\norm{\eta(b)-x}<\norm{x}$. Clearly $b\neq 0$. By definition 
$x=\lim\eta(a_n)$. Therefore
$\norm{\eta(b)-\eta(a_n)}<\norm{\eta(a_n)}$ for almost all $n$. This 
equation is equivalent to \rf{strfaith} and shows that the
weight $h$ is not strictly faithful. In this way we showed that 
Assertion 1 implies Assertion 2.\vs

We shall prove the converse. We now assume that the mapping 
$\eta(a)\mapsto\eta(a^*)$ is closable.
Take a sequence $\seq{a_n}{n\in\natu }$ in $\Dom(\eta)\cap 
\Dom(\eta)^*$ such that the
sequence $\seq{\eta(a_n^*)}{n\in\natu }$ is norm convergent to zero 
and the sequence
$\seq{\eta(a_n)}{n\in\natu }$ is bounded. Passing to a suitable
subsequence, we may assume that the sequence 
$\seq{\eta(a_n^*)}{n\in\natu }$ is norm convergent to
zero and the sequence $\seq{\eta(a_n)}{n\in\natu }$ is weakly 
convergent to some element $x\in\Hil$.
Then the pair $(0,x)$ belongs to the (norm$\,\times\,$weak)-closure 
of the $\Graph S_0=\set{(\eta(a^*),\eta(a))}
{a\in\Dom(\eta)\cap\Dom(\eta)^*}$. In Hilbert spaces the weak and the 
norm closures of any convex set coincide.
Therefore denoting by $S$ the closure of $S_0$ we obtain $x=S0=0$. 
Now we compute as follows. Let $b$ be an element in the
$C^*$-algebra $A$ such that the inequality \rf{strfaith} holds. Then 
$h(b^*b)\leq
2\text{Re}(\eta(b)|\eta(a_n))$ and taking the limit $n\to\infty$, we 
obtain $h(b^*b)\leq 0$.  Now the faithfulness of
$h$ shows that $b=0$. This proves that Assertion 2 implies Assertion 1.\vs

The equivalence of Assertions 2 and 3 is the well known part of the 
Tomita-Takesaki theory.

\end{pf}

\setcounter{equation}{0}\setcounter{Thm}{0}
\subsection{Tomita-Takesaki theory}
\label{TomTak}\Vs{2}

In this section we collect the main results on Tomita-Takesaki theory 
used in this paper. The less
known results will be proven. Throughout the Section, $M$ will be a 
von Neumann algebra of
operators acting on a Hilbert space $\Hil$. We shall assume that the embedding
$M\subset\skrL(\Hil)$ is standard. Then any normal semifinite 
faithful weight $h$ on $M$ can be
represented by a closed GNS map $\eta$ densely defined on $M$ with 
range dense in $\Hil$, related
to $h$ by the formula
\[
h(a^*a)=\left\{
\begin{array}{c@{\ }c}
\IS{\eta(a)}{\eta(a)}&\mbox{ if }a\in\Dom(\eta)\\+\infty&\mbox{ 
otherwise}\Vs{4}
\end{array}\right.
\]
for any $a\in M$. Similarly any normal semifinite faithful weight $k$ 
on $M'$ can be represented by a closed GNS map
$\theta'$ densely defined on $M'$ with range dense in $\Hil$. In what 
follows, $\theta$ will denote the commutant of $\theta'$.
With this notation we have:

\begin{Thm}\Vs{2}
\label{19.08.2002a}

{\rm 1}. There exists unique closed conjugate linear operator $S_{\rm 
rel}$ such that $\set{\theta(a)}{a\in\Dom(\theta)\cap\Dom(\eta)^*}$ is
a core for $S_{\rm rel}$ and
$
S_{\rm rel}\theta(a)=\eta(a^*)
$
for any $a\in\Dom(\theta)\cap\Dom(\eta)^*$.\vs

{\rm 2}. The adjoint operator $F_{\rm rel}=S_{\rm rel}^*$ has a core 
$\set{\theta'(b)}{b\in\Dom(\theta')\cap\Dom(\eta')^*}$ and
$
F_{\rm rel}\theta'(b)=\eta'(b^*)
$
for any $b\in\Dom(\theta')\cap\Dom(\eta')^*$.\vs

{\rm 3}. The operator $\Delta_{\rm rel}=F_{\rm rel}S_{\rm rel}$ 
depends only on the weights $h$ and $k$ $($it is independent of the 
choice of the GNS
maps $\eta$ and $\theta'$ representing $h$ and $k)$.\vs

{\rm 4}. The operator $J_{\rm rel}$ appearing in the polar 
decomposition $S_{\rm rel}=J_{\rm rel}\Delta^{1/2}_{\rm rel}$ is 
antiunitary and $J_{\rm
rel}MJ_{\rm rel}^*=M'$.\vs
\end{Thm}

To see that the domain of $S_{\rm rel}$ is dense in $\Hil$ we notice that
$a^*b\in\Dom (\theta )$ and $(a^*b)^*=b^*a\in\Dom (\eta )$ for any
$a\in\Dom (\eta )$ and $b\in\Dom (\theta )$. Therefore
$a^*\theta (b)=\theta (a^*b)$ belongs to the domain of $S_{\rm rel}$ and
$S_{\rm rel}a^*\theta (b)=\eta (b^*a)=b^*\eta (a)$. It turns out that
$
\{a^*\theta (b):a\in\Dom (\eta ),\;b\in\Dom (\theta )\}
$
is a core for $S_{\rm rel}$. In fact one can prove a little more:
\begin{Prop}
\label{19.08.2002b}

Let $\Dom_1\subset\Dom (\eta )$ be a core for $\eta$ and
$\Dom_2\subset\Dom (\theta )$ be a core for $\theta$. Then
\[
\{a^*\theta (b):a\in\Dom_1,\;b\in\Dom_2\}
\] is a core for $S_{\rm rel}$.
\end{Prop}
The similar statement holds for $F_{\rm rel}$.\vs

According to \cite{Connes}, the relative modular operator 
$\Delta_{\rm rel}$ is called the Radon-Nikodym derivative of $h$ with 
respect to
$k$ and denoted by $\frac{dh}{dk}$. Clearly $\Delta_{\rm 
rel}=\frac{dh}{dk}$ is a positive self-adjoint operator acting on 
$\Hil$.

\begin{Thm}
Let $h$ be a weight on $M$ and $k$ be a normal semifinite faithful
weight on $M'$. Then for any $a\in M$and $t\in\real$ the element
\[
\left(\frac{dh}{dk}\right)^{it}a\left(\frac{dh}{dk}\right)^{-it}
\]
belongs to $M$ and is independent of the choice of the weight $k$.
\end{Thm}
We set
\[
\sigma_t^h(a)=
\left(\frac{dh}{dk}\right)^{it}a\left(\frac{dh}{dk}\right)^{-it} .
\]
Clearly $\seq{\sigma_t^h}{t\in\real}$ is a strongly continuous one
parameter group of automorphisms of $M$. It is called the modular 
automorphism group related to the weight $h$.
\begin{Thm}
Let $h$ and $h_1$ be normal semifinite faithful weights on $M$ and
$k$ be a normal semifinite faithful weight on $M'$. Then for any
$t\in\real$ the element
\[
\left(\frac{dh}{dk}\right)^{it}\left(\frac{dh_1}{dk}\right)^{-it}
\]
belongs to $M$ and is independent of the choice of the weight $k$.
\end{Thm}

We set
\[
(Dh:Dh_1)_t=
\left(\frac{dh}{dk}\right)^{it}\left(\frac{dh_1}{dk}\right)^{-it} .
\]
Clearly $((Dh:Dh_1)_t)_{t\in\real}$ is a strongly continuous one
parameter family of unitaries in $M$. It is called the Radon-Nikodym
cocycle. One can easily verify that
\[
\begin{array}{r@{\;}c@{\;}l}
(Dh:Dh_1)_{t+\tau}
    &=& \sigma_{\tau}^h((Dh:Dh_1)_t)(Dh:Dh_1)_{\tau} \\
    &=& (Dh:Dh_1)_t \sigma_t^{h_1}((Dh:Dh_1)_{\tau}), \Vs{5} \\
\sigma_t^h(a)&=&(Dh:Dh_1)_t \sigma_t^{h_1}(a)(Dh:Dh_1)_t^* \Vs{6}
\end{array}
\]
for any $t,\;\tau\in\real$ and $a\in M$.
\begin{Thm}
\label{19.08.2002c}
Let $\eta$ and $\theta$ be closed GNS mappings densely defined on $M$
with the ranges dense in $\Hil$, $S_{\rm rel}$ be the corresponding
closed conjugate linear operator introduced in Statement $1$ of 
Theorem \reef{19.08.2002a} and
$a\in M$. Assume that $a\in\Dom(\theta )$ and
$\theta (a)\in\Dom (S_{\rm rel})$. Then $a^*\in\Dom (\eta )$ and
$($obviously$)$ $S_{\rm rel}\theta (a)=\eta (a^*)$.
\end{Thm}
\begin{pf}
By the Theorem \reef{podkomutant} $\eta''=\eta$. Therefore it is
sufficient to show that
\begin{equation}
\label{hoshi}
a^*\eta'(b)=bS_{\rm rel}\theta (a)
\end{equation}
for any $b\in\Dom (\eta')$. For any $c\in\Dom (\theta')$ we have
\[
\begin{array}{r@{\;}c@{\;}l}
(\theta'(c)|bS_{\rm rel}\theta (a))
  &=& (b^*\theta'(c)|S_{\rm rel}\theta (a)) \\
  &=& (\theta'(b^*c)|S_{\rm rel}\theta (a)) \Vs{5}.
\end{array}
\]
The reader should notice that $b^*c\in\Dom (\theta')$ and
$(b^*c)^*=c^*b\in\Dom (\eta')$ (because $\Dom (\theta')$ and
$\Dom (\eta')$ are left ideals in $M'$). Therefore (cf. Statement 2
of Theorem D1) $\theta'(b^*c)\in\Dom (F_{\rm rel})$ and
$F_{\rm rel}\theta'(b^*c)=\eta'(c^*b)=c^*\eta'(b)$. So we have
\[
\begin{array}{r@{\;}c@{\;}l}
(\theta'(c)|bS_{\rm rel}\theta (a))
   &=& (\theta'(b^*c)|S_{\rm rel}\theta (a)) \\  \Vs{5}
   &=& (\theta (a)|F_{\rm rel}\theta'(b^*c)) \\  \Vs{5}
   &=& (\theta (a)|c^*\eta'(b))=(c\theta (a)|\eta'(b)) \\  \Vs{5}
   &=& (a\theta'(c)|\eta'(b))=(\theta'(c)|a^*\eta'(b))
\end{array}
\]
and \rf{hoshi} follows.
\end{pf}

In the next theorem we set $\theta =\eta$. In this case we write
$S,\;F,\;\Delta$ and $J$ instead of $S_{\rm rel},\;F_{\rm rel},\;
\Delta_{\rm rel}$ and $J_{\rm rel}$. One can easily verify that
$S$ and $F$ are involutive: $S^{-1}=S$ and $F^{-1}=F$. It implies
that $J^2=1,\;J^*=J$ and $J\Delta J=\Delta^{-1}$. With this notation
we have
\begin{Thm}
\label{18.08.2002a}
Let $h$ be a weight associated with the closed densely defined GNS
mapping with dense range. Then \vs

$1.$ $\Dom (\eta )$ is $\sigma^h$-invariant and
\[
\eta (\sigma^h_t(a))=\Delta^{it}\eta (a)
\]
for any $a\in\Dom (\eta )$ and $t\in\real$. \vs

$2.$ $\Dom (\eta')=J\Dom (\eta )J$ and
\[
\eta'(JaJ)=J\eta (a)
\]
for any $a\in\Dom (\eta )$.
\end{Thm}

We shall also use

\begin{Thm}
\label{19.08.2002d}
Let $\eta$ be a closed GNS mapping densely defined on $M$ with range
dense in $\Hil$, $h$ be the weight on $M$ related to $\eta$ and
$\sigma_{i/2}^h$ be the analytic generator of the corresponding
modular automorphism group. Assume that $a\in\Dom (\eta )$ and
$b\in\Dom (\sigma_{i/2}^h)$. Then $ab\in\Dom (\eta )$ and
\begin{equation}
\label{hoshitrois}
\eta (ab)=J\sigma_{i/2}^h(b)^*J\eta (a) .
\end{equation}
\end{Thm}
\begin{pf}
We know that $\sigma_t^h(b)=\Delta^{it}b\Delta^{-it}$. Let
$c=\sigma_{i/2}^h(b)$. Then $b\Delta^{1/2}\subset\Delta^{1/2}c$ and
by passing to the adjoint operators
\begin{equation}
\label{hoshideux}
c^*\Delta^{1/2}\subset \Delta^{1/2}b^* .
\end{equation}
Let $\seq{u_i}{i\in I}$ be a net of positive elements of $\Dom (\eta )$
converging strongly to $I$. Then $u_ia$ and $(u_ia)^*=a^*u_i$ belong
to $\Dom (\eta )$. Therefore $\eta (a^*u_i)=S\eta (u_ia)=
J\Delta^{1/2}\eta (u_ia)=\Delta^{-1/2}J\eta (u_ia)\in
\Dom (\Delta^{1/2})$. By \rf{hoshideux} $b^*$ maps
$\Dom (\Delta^{1/2})$ into $\Dom (\Delta^{1/2})$. Hence
$\eta (b^*a^*u_i)=b^*\eta (a^*u_i)\in\Dom (\Delta^{1/2})=\Dom (S)$.
Using Theorem *** we see that $u_iab=(b^*a^*u_i)^*\in\Dom (\eta )$.
Remember that $Jc^*J\in M'$ commute with $u_i\in M$ we obtain
\[
\begin{array}{r@{\;}c@{\;}l}
u_iJc^*J\eta (a)
    &=& Jc^*Ju_i\eta (a) \\  \Vs{5}
    &=& Jc^*J\eta (u_ia) = Jc^*JS\eta (a^*u_i)  \\  \Vs{5}
    &=& Jc^*\Delta^{1/2}\eta (a^*u_i)
           = J\Delta^{1/2}b^*\eta (a^*u_i) \\  \Vs{5}
    &=& S\eta (b^*a^*u_i)  =  \eta (u_iab) .
\end{array}
\]
We take the limit with respect to $i$. Since $\eta$ is closed, we 
conclude that $ab\in\Dom (\eta )$ and
\begin{equation}
\eta (ab)=Jc^*J\eta (a) .
\end{equation}
Clearly this formula coincides with \rf{hoshitrois}.
\end{pf}

The modular group and Connes' Radon-Nikodym cocycle are canonical 
constructions. If $\tau$ is a normal
automorphism of a von Neumann algebra $M$ then for any faithful 
semifinite normal weight $h$ on
$M$, $h\comp\tau$ is a faithful semifinite normal weight on $M$ and
$\sigma^{h\comp\tau}_t=\tau^{-1}\comp\sigma^h_t\comp\tau$ for all 
$t\in\real$. Similarly for any
faithful semifinite normal weights $h$ and $h_1$ on $M$ we have: 
$(D(h\comp\tau):D(h_1\comp\tau))_t=
\tau^{-1}\left((Dh:Dh_1)_t\right)$ for all $t\in\real$. For 
antiautomorphisms the formulae are
slightly modified:
\begin{Prop}
\label{acanon}
Let $M$ be a von Neumann algebra and $R$ be a normal antiautomorphism 
of $M$. Then
\begin{enumerate}
\item For any faithful semifinite normal weight $h$ on
$M$ and any $t\in\real$ we have
\[
\sigma^{h\comp R}_t=R^{-1}\comp\sigma^h_{-t}\comp R.
\]
\item For any faithful semifinite normal weights $h$ and $h_1$ on $M$ 
and any $t\in\real$ we have:
\[
\left(D(h\comp R):D(h_1\comp 
R)\Vs{3}\right)_t=R^{-1}\left((Dh_1:Dh)_{-t}\Vs{3}\right).
\]
\end{enumerate}
\end{Prop}

\setcounter{equation}{0}\setcounter{Thm}{0}
\subsection{Tensor Product}
\label{Tensor}\Vs{2}

In this section we discuss the tensor product of GNS-mappings. Let 
$A$ be a non-degenerated
$C^*$-algebra of operators acting on a Hilbert space $\Hil$ and $\eta 
:A\rightarrow\Hil$ be a
densely defined closed GNS mapping. Similarly let $B$ be a 
non-degenerated $C^*$-algebra of
operators acting on a Hilbert space ${\mathcal K}$ and $\theta 
:B\rightarrow {\mathcal K}$ be a
densely defined closed GNS mapping. Then $A\otimes B$ is a 
non-degenerated $C^*$-algebra of
operators acting on $\Hil\otimes {\mathcal K}$. We shall consider the 
algebraic tensor product
$\eta\otimes_{\rm alg}\theta$. By definition $\Dom (\eta\otimes_{\rm 
alg}\theta )=\Dom (\eta
)\otimes_{\rm alg}\Dom (\theta )$. Clearly
\[
\label{tensorcomm1}
(a\otimes b)(\eta\otimes_{\rm alg}\theta )(c)=c(\eta'(a)\otimes\theta'(b))
\]
for any $c\in\Dom (\eta )\otimes_{\rm alg}\Dom (\theta )$, $a\in\Dom (\eta')$
and $b\in\Dom (\theta')$. Using this formula one can easily show that the
mapping $\eta\otimes_{\rm alg}\theta$ is closable. Indeed if $(c_{\alpha})$
is a net of elements of $\Dom (\eta )\otimes_{\rm alg}\Dom (\theta )$
converging strongly to $0$ such that
$(\eta\otimes_{\rm alg}\theta )(c_{\alpha})\rightarrow 
x\in\Hil\otimes{\mathcal K}$ then using \rf{tensorcomm1} and passing
to the limit we obtain $(a\otimes b)x=0$ for all $a\in\Dom (\eta')$ and
$b\in\Dom (\theta')$. Hence $x=0$ and $\eta\otimes_{\rm alg}\theta$ 
is closable.\vs

By definition $\eta\otimes\theta$ is the closure of
$\eta\otimes_{\rm alg}\theta$. We shall prove that $\eta\otimes\theta$ is a GNS
mapping acting from $A\otimes B$ into $\Hil\otimes {\mathcal K}$.\vs

Let $a\in A$, $b\in B$ and $c\in\Dom (\eta\otimes\theta )$. Then there exists a
net $(c_{\alpha})$ of elements of $\Dom (\eta )\otimes_{\rm alg}\Dom (\theta )$
such that $c_{\alpha}\rightarrow c$ strongly and
$(\eta\otimes_{\rm alg}\theta )(c_{\alpha})\rightarrow (\eta\otimes\theta)(c)$
in norm. Passing to the limit in the obvious equality
\[
(\eta\otimes_{\rm alg}\theta )((a\otimes b)c_{\alpha})
=(a\otimes b)(\eta\otimes_{\rm alg}\theta )(c_{\alpha})
\]
we conclude that $(a\otimes b)c\in\Dom (\eta\otimes\theta )$ and
\[
(\eta\otimes_{\rm alg}\theta )((a\otimes b)c)
=(a\otimes b)(\eta\otimes_{\rm alg}\theta )(c) \;.
\]
Using the same method one may replace in the above statement 
$a\otimes b$ by any
element of $A\otimes B$. It shows that $\eta\otimes\theta$ is a GNS-mapping.

\begin{Rem}
By our definition $\Dom (\eta )\otimes_{\rm alg}\Dom (\theta )$ is a core for
$\eta\otimes\theta$.
\end{Rem}

In the same way one introduces the tensor product of GNS mappings defined on
von Neumann algebras.

\begin{Thm}
Let $\eta :A\rightarrow\Hil$ and $\theta :B\rightarrow {\mathcal K}$ be closed
densely defined GNS mappings. Then
\begin{equation}
\label{18.08.2002}
(\eta\otimes\theta )'=\eta'\otimes\theta' \;.
\end{equation}
\end{Thm}

\begin{pf}
For simplicity assume that $\eta''$ and $\theta''$ are faithful and have dense
ranges. Then $(\eta\otimes\theta)''$ is also faithful and has dense range. We
may use Tomita-Takesaki theory. Using Proposition \reef{19.08.2002b} 
with $\theta =\eta$ replaced by $\eta\otimes\theta$ and $D_1=D_2$ 
replaced by
$\Dom (\eta )\otimes_{\rm alg}\Dom (\theta )$ one can easily show that the
operator $S$ related to $(\eta\otimes\theta )''$ is the tensor product of
operators $S$ related to $\eta$ and $\theta$. Consequently we have the same
statements for operators $\Delta$ and $J$. Now \rf{18.08.2002} 
follows from Statement 2 of Theorem \reef{18.08.2002a}.
\end{pf}

\begin{Cor}
Let
\[
\eta (a)=\sum_{\ell}a\Omega_{\ell}
\mbox{\hspace{1cm}and\hspace{1cm}}
\theta (b)=\sum_kb\Psi_k
\]
be exact vector presentation of GNS mappings $\eta$ and $\theta$. Then
\begin{equation}
\label{tensortensor}
(\eta\otimes\theta )(c)=\sum_{k,\ell}c(\Omega_{\ell}\otimes\Omega_k)
\end{equation}
is an exact vector presentation of $\eta\otimes\theta$.
\end{Cor}

\begin{pf}
Let $P_{\ell}\in\skrL (\Hil )$ be the orthogonal projection onto
$\overline{A\Omega_{\ell}}$ and $Q_k\in\skrL ({\mathcal K})$ be the orthogonal
projection onto $\overline{B\Psi_k}$. Then $P_{\ell}\otimes Q_k$ is the
orthogonal projection onto
$\overline{(A\otimes B)(\Omega_{\ell}\otimes\Psi_k)}$ According to 
Proposition \reef{18.08.2002c}
the left ideal in $A'$ (in $B'$ resp.) generated by $\{P_1,P_2,\dots\}$
(by $\{Q_1,Q_2,\dots\}$ resp.) is a core for $\eta'$ (for $\theta'$ 
resp.). Therefore
the left ideal generated by $\{P_{\ell}\otimes Q_k: \ell , k=1,2,\dots\}$ is a
core of $\eta'\otimes\theta'=(\eta\otimes\theta )'$. Using again 
Proposition \reef{18.08.2002c} we see
that \rf{tensortensor} is an exact vector presentation.
\end{pf}

\setcounter{equation}{0}\setcounter{Thm}{0}
\subsection{Analytic generators}
\label{AnaGene}

In this Section we give necessary information on analytic generators 
of one parameter groups of
isometries of Banach spaces. Let $B$ be a Banach space and 
$\seq{\alpha_t}{t\in\real}$ be a one
parameter group of isometries of $B$: $\norm{\alpha_t(a)}= \norm{a}$,
$\alpha_t(\alpha_s(a))=\alpha_{t+s}(a)$ and $\alpha_0(a)=a$ for any 
$a\in B$ and $t,s\in\real$. We
shall always assume that the group is pointwise continuous: for any $a\in B$,
$\norm{\alpha_t(a)-a}\rightarrow 0$ when $t\rightarrow 0$.\vs

  Let us recall \cite{Zsi} that the analytical generator $\alpha_i$ of the group
$\seq{\alpha_t}{t\in\real}$ is the linear operator acting on $B$ in 
the following way:

  For any $a,b\in B$: $a\in\Dom(\alpha_{i})$ and $b=\alpha_{i}(a)$  if 
and only if there exists a
mapping $z\mapsto a_z\in B$ continuous on the strip 
$\set{z\in\compl}{\ \Im z\in[0,1]}$ and holomorphic
in the interior of this strip such that $a_t=\alpha_t(a)$ for all 
$t\in\real$ and
$a_{i}=b$.\vs

It is known that $\alpha_i$ is a closed densely defined linear mapping.   \vs

\begin{Prop}
\label{pierwsze} Let $B_k$  be Banach spaces, 
$\seq{\alpha^k_t}{t\in\real}$ be pointwise continuous one
parameter groups of isometries of $B_k$ and $\alpha_i^k$ be the 
corresponding analytical generators
$(k=0,1,\dots,K-1)$.

Moreover let
\[
\Psi:B_0\times B_1\times\dots\times B_{K-1}\longrightarrow \compl
\]
be a continuous K-linear form such that
\begin{equation}
\label{inva} 
\Psi(\alpha^0_i(a_0),\alpha^1_i(a_1),\dots,\alpha^{K-1}_i(a_{K-1}))=
\Psi(a_0,a_1,\dots,a_{K-1})
\end{equation}
for any $a_k\in \Dom(\alpha_i^k)$ $(k=0,1,\dots, K-1)$. Then
\begin{equation}
\label{teza} 
\Psi(\alpha^0_t(a_0),\alpha^1_t(a_1),\dots,\alpha^{K-1}_t(a_{K-1}))=
\Psi(a_0,a_1,\dots,a_{K-1})
\end{equation}
for any $t\in\real$ and $a_k\in B_k$ $(k=0,1,\dots,K-1)$.
\end{Prop}\vs

To prove this proposition we shall use the following result known as 
Carlson's Lemma (cf. \cite[page 228]{wpure}):

\begin{Prop} Let $f(\cdot)$ be a holomorphic function on the half-plane
$\compl_+=\set{z\in\compl}{\Im z>0}$. Assume that for some constant 
$M\in\real$, $f(z)e^{-M\Im z}$ is
bounded on $\compl_+$ and that $f(in)=0$ for $n\in\natu$. Then $f$ 
vanishes identically.
\end{Prop}

{\it Proof of Proposition \reef{pierwsze}}: Let
$\overline{\compl}_+=\compl_+\cup\real$ be the closure of
$\compl_+$. For any $k=1,2,\dots,K$ we set
\[
\skrD_k=\left\{a\in B_k:
\begin{array}{c}
\mbox{there exists continuous map:}\hspace{8mm} 
\overline{\compl}_+\ni z\longrightarrow a(z)\in B_k\\
\mbox{holomorphic on }\compl_+\mbox{ such that for some constant 
}M\in\real,\\ a(z) e^{-M\Im z}
\mbox{ is bounded on }\overline{\compl}_+  \mbox{ and } 
a(t)=\alpha^k_t(a) \mbox{ for all } t\in\real
\end{array}
\right\}.
\]

Let $a_k\in\skrD_k$ ($k=0,1,\dots,K-1$), $a_k(z)$ be the holomorphic 
extension of $\alpha^k_t(a_k)$ (cf. the above definition of 
$\skrD_k$) and
\[
g(z)=\Psi(a_0(z),a_1(z),\dots,a_{K-1}(z)).
\]
Then, for any $n\in\natu$, $a_k\in \Dom(\left(\alpha^k_i\right)^n)$ and
$\left(\alpha^k_i\right)^na_k=a_k(ni)$. Taking into account 
Assumption \rf{inva} we see that
$g(in)=g(0)$. Now one can easily show that the function 
$f(z)=g(z)-g(0)$ satisfies the assumptions of
the Carlson's lemma. Therefore $g(z)=g(0)$ for all $z\in 
\overline{\compl}_+$. For real $z$ this
formula coincides with \rf{teza}.\vs

To end the proof we have to show that for each $k$, $\skrD_k$ is 
dense in $B_k$. Let $a\in B_k$. For
any $\varepsilon>0$ we set
\[
a_\varepsilon=\frac{1}{2\pi}\int_\real
\frac{e^{1+i\tau}}{{(1+i\tau)}^2}\,\alpha^k_{\varepsilon\tau}(a)\,d\tau.
\]
It turns out that $a_\varepsilon\in\skrD_k$. Indeed one can easily 
verify that the map
\[
\overline{\compl}_+\ni z\longrightarrow 
a_\varepsilon(z)=\frac{e^{-i\frac{z}{\varepsilon}}}{2\pi}
\int_\real\frac{e^{1+i\tau}}{{(1-i\frac{z}{\varepsilon}+i\tau)}^2}\,
\alpha^k_{\varepsilon\tau}(a)\,d\tau\in B_k
\]
is the holomorphic continuation of $\alpha^k_t(a_\varepsilon)$ 
satisfying all the requirements of
definition of $\skrD_k$. Using the Lebesgue dominating  convergence 
theorem, one can easily
show that
\[
\lim_{\varepsilon\rightarrow
0}a_\varepsilon=\left(\frac{1}{2\pi}\int_\real\frac{e^{1+i\tau}}{{(1+i\tau)}^2}\,d\tau\right)\,a=a.
\]
It shows that any element of $B_k$ is a norm-limit of elements of $\skrD_k$.
\hfill $\framebox[2.2mm]{}$

\begin{Thm}
\label{drugie} Let $B_k$  be Banach spaces, 
$\seq{\alpha^k_t}{t\in\real}$ be pointwise continuous one
parameter groups of isometries of $B_k$ and $\alpha_i^k$ be the 
corresponding analytical generators
$(k=1,2,\dots,K)$.

Moreover let
\begin{equation}
\label{Psi}
\Psi:B_1\times B_2\times\dots\times B_{K-1}\longrightarrow B_K
\end{equation}
be a continuous $(K-1)$-linear mapping such that
\begin{equation}
\label{gwi1} \hspace*{-40mm}\Psi(a_1,a_2,\dots,a_{K-1})\in 
\Dom(\alpha^K_i)\vspace{-2mm}
\end{equation}
and
\begin{equation}
\label{A101} \alpha^K_i(\Psi(a_1,a_2,\dots,a_{K-1}))=
\Psi(\alpha^1_i(a_1),\alpha^2_i(a_2),\dots,\alpha^{K-1}_i(a_{K-1}))\vs
\end{equation}
for any $a_k\in \Dom(\alpha_i^k)$ $(k=1,2,\dots ,K-1)$. Then
\begin{equation}
\label{A103} \alpha^K_t(\Psi(a_1,a_2,\dots,a_{K-1}))=
\Psi(\alpha^1_t(a_1),\alpha^2_t(a_2),\dots,\alpha^{K-1}_t(a_{K-1}))
\end{equation}
for any $t\in\real$ and $a_k\in B_k$ $(k=1,2,\dots,{K-1})$.
\end{Thm}

\begin{pf} 
Let $B_0$ be the set of all $a\in B_K^*$ such that the mapping $t\in{\Bbb R}\to a\circ\alpha_t^K\in B_K^*$ is norm continuous.  Then $B_0$ is a closed subspace of $B_K^*$.  For each $a\in B_K^*$ we set
$$
{\mathcal R}_{\varepsilon}(a)=\frac 1{\sqrt{\varepsilon\pi}}\int_{\Bbb R}a\circ\alpha_t^K\ e^{-t^2/\varepsilon}dt.
$$
Then 
$$
{\mathcal R}_{\varepsilon}(a)\circ\alpha_s^K-{\mathcal R}_{\varepsilon}(a)=\frac 1{\sqrt{\varepsilon\pi}}\int_{\Bbb R}a\circ\alpha_t^K(e^{-(t-s)^2/\varepsilon}-e^{-t^2/\varepsilon})dt
$$
and hence ${\mathcal R}_{\varepsilon}(a)\in B_0$.  Since 
$$
{\mathcal R}_{\varepsilon}(a)-a=\frac 1{\sqrt{\varepsilon\pi}}\int_{\Bbb R}(a\circ\alpha_t^K-a)e^{-t^2/\varepsilon}dt,
$$
$B_0$ is weakly$^*$ dense in $B_K^*$.  
The canonical bilinear form on
$B_0\times B_K$, will be denoted by
$\langle\cdot,\cdot\rangle$. 
Since ${\mathcal R}_{\varepsilon}(a)\circ\alpha_t^K={\mathcal R}_{\varepsilon}(a\circ\alpha_t^K)$ for all $t\in{\Bbb R}$, the formula
\begin{equation}
\label{A102} \langle\alpha^0_t(a_0),\alpha^K_t(a_K)\rangle=\langle 
a_0,a_K\rangle,
\end{equation}
where $a_0\in B_0$, $a_K\in B_K$, $t\in\real$, defines a one parameter group
$\seq{\alpha^0_t}{t\in\real}$ of isometries of $B_0$. Let 
$\alpha^0_i$ be the corresponding
analytical generator. By the holomorphic continuation we have:
\begin{equation}
\label{A100} \langle\alpha^0_i(a_0),\alpha^K_i(a_K)\rangle=\langle 
a_0,a_K\rangle
\end{equation}
for any $a_0\in \Dom(\alpha^0_i)$ and $a_K\in \Dom(\alpha^K_i)$.\vs

For any $a_k\in \Dom(\alpha_i^k)$ $(k=0,1,\dots, K-1)$ we set:
\[
\Psi'(a_0,a_1,\dots,a_{K-1})=\langle a_0,\Psi(a_1,\dots,a_{K-1})\rangle.
\]
Clearly $\Psi'$ is a continuous K-linear form defined on $B_0\times 
B_1\times\dots\times B_{K-1}$.
Using \rf{A100} and \rf{A101} one can easily check that
\[
\Psi'(\alpha^0_i(a_0),\alpha^1_i(a_1),\dots,\alpha^{K-1}_i(a_{K-1}))= 
\Psi'(a_0,a_1,\dots,a_{K-1})
\]
for any $a_k\in \Dom(\alpha_i^k)$ $(k=0,1,\dots, K-1)$. Now, 
Proposition \reef{pierwsze} shows that
\[
\langle\alpha^0_t(a_0),\Psi(\alpha^1_t(a_1),\dots,\alpha^{K-1}_t(a_{K-1}))\rangle= 
\langle
a_0,\Psi(a_1,\dots,a_{K-1})\rangle
\]
for any $a_k\in B_k$ $(k=0,1,\dots, K-1)$ and $t\in\real$. Comparing 
this result with \rf{A102} we
get \rf{A103}.

\end{pf}\vs

In particular, inserting in Theorem \reef{drugie}, $K=2$, $B_1=B_2=B$ 
and $\Psi=\id_B$ we get

\begin{Cor}
\label{uniquegr}
Let $B$ be a Banach space and $\seq{\alpha^1_t}{t\in\real}$ and $\seq{\alpha^2_t}{t\in\real}$
be one parameter groups of isometries acting on $B$. Assume that the 
analytic generators
$\alpha^1_i\subset\alpha^2_i$. Then $\alpha^1_i=\alpha^2_i$ and 
$\alpha^1_t=\alpha^2_t$ for all
$t\in\real$.
\end{Cor}

Let us go back to Theorem \reef{drugie}. Assume for the moment that 
the image of \rf{Psi} is linearly
dense in $B_K$. Then the linear span of elements of the form 
\rf{gwi1} is a core for
$\alpha_i^K$. This fact follows from the following

\begin{Thm}
\label{core}
  Let $\alpha=\seq{\alpha_t}{t\in\real}$ be a one parameter group of 
isometries acting on a Banach space
$B$ and $\Dom$ be an $\alpha$-invariant linear subset of 
$\Dom(\alpha_i)$. Assume that $\Dom$ is dense in $B$.
Then $\Dom$ is a core for $\alpha_i$.
\end{Thm}

To prove this Theorem we shall use the regularizing operator 
$\skrR_\alpha$ acting on $B$. By
definition
\begin{equation}
\label{reg} 
\skrR_\alpha=\frac{1}{\sqrt{\pi}}\int_{-\infty}^{\infty}e^{-\tau^2}\alpha_\tau\,d\tau.
\end{equation}
One can easily verify that $\skrR_\alpha$ is a norm one operator 
acting on $B$, the range
$\skrR_\alpha(B)\subset \Dom(\alpha_i)$ and
\[
\alpha_i\comp\skrR_\alpha=\frac{1}{\sqrt{\pi}}\int_{-\infty}^{\infty}e^{-(\tau-i)^2}\alpha_\tau\,d\tau.
\]
is a bounded operator acting on $B$ with norm less or equal $e$. 
Endowing $\Dom(\alpha_i)$ with the
graph norm ($\norm{b}_{\rm Graph}=\norm{b}+\norm{\alpha_i(b)}$ for 
any $b\in \Dom(\alpha_i)$) we see that
$\skrR_\alpha$ maps $B$ into $\Dom(\alpha_i)$ in a continuous way. If 
$\Dom$ is a dense subset of $B$, then
$\skrR_\alpha(\Dom)$ is a dense subset of $\skrR_\alpha(B)$ in the 
sense of the graph topology of
$\Dom(\alpha_i)$. It means that the restrictions of $\alpha_i$ to 
$\skrR_\alpha(\Dom)$ and
$\skrR_\alpha(B)$ have the same closures:
\begin{equation}
\label{jeden}
\overline{\left.\alpha_i\right|_{\skrR_\alpha(\Dom)}}=
\overline{\left.\alpha_i\right|_{\skrR_\alpha(B)}}.
\end{equation}

\begin{Prop}
\label{pomoc}
Let $\alpha=\seq{\alpha_t}{t\in\real}$ be a one parameter group of 
isometries acting on a
Banach space $B$, $b\in B$ and $\Dom\subset B$ be an 
$\alpha$-invariant linear subset. Then \vs

1. $\skrR_\alpha(b)$ belongs to the smallest closed 
$\alpha$-invariant linear subspace of $B$ containing $b$.\vs

2. $b$ belongs to the smallest closed $\alpha$-invariant linear 
subspace of $B$ containing $\skrR_\alpha(b)$.\vs

3. $\skrR_\alpha(\Dom)$ and $\Dom$ have the same closures:
\[
\overline{\skrR_\alpha(\Dom)}=\overline{\Dom}.
\]
\end{Prop}

\begin{pf} Ad 1 and 2. It is sufficient to show that for any 
continuous linear functional $\varphi$ on $B$ we
have:
\begin{equation}
\label{equi}
  \left(
\begin{array}{c}
\varphi(\alpha_t(b))=0\\{\rm for\ all\ }{t\in\real}
\end{array}
\right) \Longleftrightarrow \left(
\begin{array}{c}
\varphi(\alpha_t(\skrR_\alpha(b)))=0\\{\rm for\ all\ }{t\in\real}
\end{array}
\right).
\end{equation}
For any $\varphi\in B^d$ we set: $f_\varphi(t)=\varphi(\alpha_t(b))$ 
and $g_\varphi(t)=\varphi(\alpha_t(\skrR_\alpha(b)))$. Then 
$f_\varphi$
and $g_\varphi$ are bounded continuous functions on $\real$ and may 
be considered as tempered distributions on $\real$. In what follows,
$\widetilde{f}_\varphi$ and $\widetilde{g}_\varphi$ will denote the 
Fourier transform of $f_\varphi$ and $g_\varphi$. Therefore
\rf{equi} is equivalent to
\begin{equation}
\label{equi2}
  \left(
\widetilde{f}_\varphi=0
\right) \Longleftrightarrow \left(\Vs{4}
\widetilde{g}_\varphi=0
\right).
\end{equation}
Taking into account \rf{reg} we obtain:
\[
  \varphi\left(\skrR_\alpha(b)\right)=
\frac{1}{\sqrt{\pi}}\int_{-\infty}^{\infty}e^{-\tau^2}\varphi\left(\alpha_\tau(b)\right)\,d\tau.
\]
Replacing $b$ by $\alpha_t(b)$ we obtain
\[
g_\varphi(t)= 
\frac{1}{\sqrt{\pi}}\int_{-\infty}^{\infty}e^{-\tau^2}f_\varphi(t+\tau)\,d\tau
= 
\frac{1}{\sqrt{\pi}}\int_{-\infty}^{\infty}e^{-(t-\tau)^2}f_\varphi(\tau)\,d\tau.
\]
On the right hand side we have the convolution product of the Gaussian function
$\frac{1}{\pi}e^{-\tau^2}$ with $f_\varphi$. Passing to the Fourier 
transform we see that
$\widetilde{g}_\varphi$ equals to $\widetilde{f}_\varphi$ multiplied 
by the Fourier transform of
the Gaussian function. The latter has no zero points in $\real$ and 
\rf{equi2} follows.\vs

Ad 3. It follows immediately from Statements 1 and 2.
\end{pf}\vs

Now we are able to prove Theorem \reef{core}.\vs

\begin{pf} We shall use Statement 3 of Proposition \reef{pomoc} in 
the following context. Instead of $B$ we take
$\Dom(\alpha_i)$ endowed with the graph norm. Clearly $\alpha$ is a 
group of isometries of $\Dom(\alpha_i)$.
Therefore for any $\alpha$-invariant subset $\Dom\subset 
\Dom(\alpha_i)$, the closures (in the sense of the
graph topology) of
$\skrR_\alpha(\Dom)$ and $\Dom$ coincide. In other words
\begin{equation}
\label{dwa}
\overline{\left.\alpha_i\right|_{\skrR_\alpha(\Dom)}}= 
\overline{\left.\alpha_i\right|_\Dom}.
\end{equation}
Inserting in this formula $\Dom=\Dom(\alpha_i)$ and using the 
inclusion $\Dom(\alpha_i)\subset B$ we obtain
\begin{equation}
\label{trzy} \overline{\left.\alpha_i\right|_{\skrR_\alpha(B)}}= \alpha_i.
\end{equation}
Combining now \rf{jeden}, \rf{dwa} and \rf{trzy} we get:
$\overline{\left.\alpha_i\right|_\Dom}=\alpha_i$. It means that 
$\Dom$ is a core for $\alpha_i$.
\end{pf}

\begin{Thm}
\label{removeReg}
Let $\alpha=\seq{\alpha_t}{t\in\real}$ be a one parameter group of 
isometries acting on a
Banach space $B$ and $a,b\in B$. Assume that 
$\alpha_i\comp\skrR_\alpha(a)=\skrR_\alpha(b)$. Then 
$a\in\Dom(\alpha_i)$ and $\alpha_i(a)=b$.
\end{Thm}

\begin{pf} We shall use Statement 2 of Proposition \reef{pomoc} in 
the following context. Instead of $B$ we take
$\Btilde=B\oplus B$. Clearly $\seq{\alphatilde_t}{t\in\real}$, where 
$\alphatilde_t=\alpha_t\oplus\alpha_t$ for all $t\in\real$, is
a group of isometries of $\Btilde$. 
$\Graph(\alpha_i)=\set{(a,\alpha_i(a))}{a\in\Dom(\alpha_i)}$ is a 
closed, $\alphatilde$-invariant
subspace of $\Btilde$. We assumed that 
$\skrR_{\tilde{\alpha}}(a,b)=(\skrR_\alpha(a),\skrR_\alpha(b))\in\Graph(\alpha_i)$. 
Using
Statement 2 of Proposition \reef{pomoc} and remembering that 
$\Graph(\alpha_i)$ is $\alphatilde$ invariant we see that
$(a,b)\in\Graph(\alpha_i)$. It means that $a\in\Dom(\alpha_i)$ and 
$b=\alpha_i(a)$.
\end{pf}

\begin{Rem} {\rm One can easily generalize Proposition 
\reef{pierwsze} and Theorem \reef{drugie} for
mappings
$\Psi$ that are conjugate linear with respect to some variable (and 
linear with respect to remaining
ones). If $\Psi(a_1,a_2,\dots,a_K)$ is conjugate linear with respect 
to $a_k$, then in all entries
the $\alpha^k_i$ should be replaced by $\alpha^k_{-i}$ (the analytic 
generator of
$\seq{\alpha^k_{-t}}{t\in\real}$).}
\end{Rem}

\end{document}